\documentclass[11pt,a4paper]{article}

\usepackage[utf8]{inputenc}
\usepackage{amsmath}
\usepackage{amsfonts}
\usepackage{amssymb}
\usepackage{graphicx}
\usepackage{xcolor}

\usepackage[left=3cm,right=3cm,top=3cm,bottom=3cm]{geometry}

\newtheorem{theorem}{Theorem}[section]

\newtheorem{corollary}[theorem]{Corollary}
\newtheorem{definition}[theorem]{Definition}
\newtheorem{example}[theorem]{Example}

\newtheorem{remark}[theorem]{Remark}

\newenvironment{proof}[1][Proof]{\noindent\textbf{#1.} }{\ \rule{0.5em}{0.5em}}
\usepackage{algorithm,algorithmic}

\usepackage{natbib} 

\usepackage{hyperref}
\hypersetup{
    colorlinks=true,
    linkcolor=black,
    filecolor=magenta,   
    citecolor=blue,   
    urlcolor=cyan,
    pdfpagemode=FullScreen,
}

\author{}
\date{}
\title{ \textbf{A new robust approach for multinomial logistic regression with complex design model}}
\author{Elena Castilla$^{1}$ and Pedro J. Chocano$^{2}$
\\$^{1}${\small Department of Statistics and O.R., Complutense University of Madrid,  Spain}
\\$^{2}${\small Department of Algebra, Geometry and Topology, Complutense University of Madrid,  Spain}
}
\begin{document}
\maketitle
\begin{abstract}
Robust estimators and Wald-type tests are developed for the multinomial logistic regression based on $\phi$-divergence measures. The robustness of the proposed estimators and tests is proved through the study of their influence functions and it is also illustrated with two numerical examples and an extensive simulation study.
\end{abstract}
\section{Introduction}

Multinomial logistic regression model, also known as polytomous logistic regression model, is widely used in health and life sciences for analyzing nominal qualitative response variables and their relationship with respect to their corresponding explanatory variables or covariates.  \cite{daniels1997hierarchical} used hieratical multinomial logistic regresion models  to examine how the rates of cardiac procedures depend on patient-level characteristics, including age, gender and race. \cite{dreassi2007polytomous}  also used multinomial logistic regression to  detect uncommon risk factors related to oral cavity, larynx and lung cancers.   Recently, \cite{ke2016semi} proposed a risk prediction model using semi-varying coefficient multinomial logistic regression  to assess  correct prediction rates  when classifying the patients with early rheumatoid arthritis.  Further examples of application of these methods can be found in \cite{blizzard2007log}, \cite{bull2007confidence} and  \cite{bertens2016nomogram}, among others.  Although most of classical literature deals with the cases of simple random sampling scheme, the application of multinomial logistic regression model under complex survey setting (with stratification, clustering or unequal selection probabilities, for example) can be found, for example, in \cite{binder1983variances, roberts1987logistic, morel1989logistic} and \cite{morel2012overdispersion}.

Most of the results mentioned above are based on (pseudo) maximum likelihood estimators (PMLEs), which are well-known to be efficient, but also non-robust. Therefore, testing procedures based on MLEs may face serious robustness problems. \cite{castilla2018_BIOMETRICS} developed density power divergence (DPD) based robust estimators (MDPDEs) and Wald-type tests for multinomial logistic regression model under simple random sampling. This approach was extended to complex design in \cite{castilla2019_ABHIK}.  In \cite{Castilla2018_ASTA}, pseudo minimum $\phi$-divergence estimators (PM$\phi$Es), as well as new estimators for the intra-cluster correlation coefficient were developed.  Estimators in the Cressie Read subfamily with tuning parameter $\lambda>0$ were shown to be an efficient alternative to classical PMLE ($\lambda=0$) for small samples sizes. However, the robustness issue was not considered. In the cited paper of \cite{castilla2019_ABHIK}, some simulation studies showed that  Cressie Read estimators with negative tuning parameter were even more robust than MDPDEs,   in terms of efficiency, for low-moderate intra-cluster correlations. However, this problem was not theoretically studied and hypothesis testing  was not considered. In this paper, we prove, through the study of the influence function, the robustness of PM$\phi$Es with $-1<\lambda<0$ and we develop robust PM$\phi$Es based Wald-type tests for testing composite null hypothesis.  In Section \ref{sec:model}, we present the multinomial logistic regression model as well as the framework necessary to define the PM$\phi$Es. Based on their asymptotic distribution, robust Wald-type tests are developed in Section \ref{Sec:Wald}. An extensive simulation study and two numerical examples  illustrate the robustness of the proposed estimators and Wald-type tests, in Section \ref{sec:mc} and Section \ref{sec:num} respectively. In Appendix \ref{sec:modelrobust}, the study of the influence function of the proposed test statistics is detailed, while in Appendix \ref{app:proofs} we present the proofs of  the main  results. Finally, in Appendix \ref{sec:app_sim}, some extensions of the Monte Carlo simulation study are presented.

\section{Multinomial logistic regression model with complex design \label{sec:model}}

We consider a population $\Omega$ partitioned into $H$ strata and the data consist of $n_h$ clusters in stratum $h$. In the $i$-th cluster $(i=1,...,n_h)$ within the $h$-th stratum $(h=1,...,H)$ we have observed for the $j$-th unit $(j=1,...,m_{hi})$ the values of a categorical response variable $Y$ with $d+1$ categories. Note that we assume there are $H$ strata, $n_h$ clusters in stratum $h$ and $m_{hi}$ units in cluster $i$ of stratum $h$.  The observed responses of the $(d+1)$-dimensional variable $Y$ are denoted by the $(d+1)$-dimensional classification vector.
\begin{equation*}
\boldsymbol{y}_{hij}=\left( y_{hij1},....,y_{hij,d+1}\right) ^{T},\ h=1,...,H,i=1,...,n_{h},\text{ }j=1,...,m_{hi},  
\end{equation*}
with $y_{hijr}$ $=1$  if the $j$-th unit selected from the $i$-th cluster of the $h$-th stratum falls in the $r$-th category and $y_{hijl}$ $=0$ for $l\neq r$. It is very common when working with dummy or qualitative explanatory variables to consider that the $k+1$ explanatory variables are common for all the individuals in the $i$-th cluster of the $h$-th stratum, being denoted as $\boldsymbol{x}_{hi}=\left( x_{hi0},x_{hi1},....,x_{hik}\right) ^{T}$, with the first one, $x_{hi0}=1$,  associated with the intercept.

Let us denote the sampling weight from the $i$-th cluster of the $h$-th stratum by $w_{hi}$. For each $i$, $h$ and $j$, the expectation of the $r$-th element of the random variable $\boldsymbol{Y}_{hij}=(Y_{hij1},...,Y_{hij,d+1})^{T}$, corresponding to the realization $\boldsymbol{y}_{hij}$, is determined by 
\begin{equation}
\pi _{hir}\left( \boldsymbol{\beta }\right) =\mathrm{E}\left[ Y_{hijr}|\boldsymbol{x}_{hi}\right] =\Pr \left( Y_{hijr}=1|\boldsymbol{x}_{hi}\right)=\dfrac{\exp \{\boldsymbol{x}_{hi}^{T}\boldsymbol{\beta }_{r}\}}{1+{\sum_{l=1}^{d}}\exp \{\boldsymbol{x}_{hi}^{T}\boldsymbol{\beta }_{l}\}}, \quad r=1,...,d
,  \label{eq:pi}
\end{equation}%
with \ $\boldsymbol{\beta }_{r}=\left( \beta _{r0},\beta _{r1},...,\beta_{rk}\right) ^{T}\in \mathbb{R}^{k+1}$, $r=1,...,d$ and the associated parameter space given by $\Theta =\mathbb{R}^{d(k+1)}$. It is clear that 
\begin{equation}
\pi _{hid+1}\left( \boldsymbol{\beta }\right) =\dfrac{1}{1+{\sum_{l=1}^{d}}\exp \{\boldsymbol{x}_{hi}^{T}\boldsymbol{\beta }_{l}\}}.  
\end{equation}%

Note that, under homogeneity, the expectation of $\boldsymbol{Y}_{hij}$ does not depends on the unit number $j$. Therefore, from now on, we will denote by 
\begin{equation*}
\widehat{\boldsymbol{Y}}_{hi}=\sum \limits_{j=1}^{m_{hi}}\boldsymbol{Y}_{hij}=\left( \sum \limits_{j=1}^{m_{hi}}Y_{hij1},...,\sum\limits_{j=1}^{m_{hi}}Y_{hij,d+1}\right) ^{T}=(\widehat{Y}_{hi1},...,\widehat{Y}_{hi,d+1})^{T} 
\end{equation*} 
the random vector of counts in the $i$-th cluster of the $h$-th stratum and by $\boldsymbol{\pi }_{hi}\left( \boldsymbol{\beta }\right)$ the $(d+1)$-dimensional probability vector with the elements given in (\ref{eq:pi}), $\boldsymbol{\pi }_{hi}\left( \boldsymbol{\beta }\right) =\left( \pi_{hi1}\left( \boldsymbol{\beta }\right) ,...,\pi _{hi,d+1}\left( \boldsymbol{\beta }\right) \right) ^{T}$.

Following this notation we can define  the empirical and the theoretical probability vectors of the model as
\begin{align}
\widehat{\boldsymbol{p}} &  =\frac{1}{\tau}(w_{11}\widehat{\boldsymbol{y}}_{11}^{T},...,w_{1n_{1}}\widehat{\boldsymbol{y}}_{1n_{1}}^{T},...,w_{H1}\widehat{\boldsymbol{y}}_{H1}^{T},...,w_{Hn_{H}}\widehat{\boldsymbol{y}}_{Hn_{H}}^{T})^{T}  \quad \text{and} \label{eq:emp_prob_vector}\\
\boldsymbol{\pi}(\boldsymbol{\beta})&=\frac{1}{\tau}(w_{11}m_{11}\boldsymbol{\pi}_{11}^{T}(\boldsymbol{\beta}),...,w_{1n_{1}}m_{1n_{1}}
\boldsymbol{\pi}_{Hn_{H}}^{T}(\boldsymbol{\beta}))^{T}, \label{eq:theo_prob_vector}
\end{align}
respectively, where $\tau=\sum \limits_{h=1}^{H}\sum \limits_{i=1}^{n_{h}}w_{hi}m_{hi}$. Probability vectors (\ref{eq:emp_prob_vector}) and (\ref{eq:theo_prob_vector}), both of dimension $(d+1)\sum_{h=1}^H n_h$, will play a basic role in the definition of PM$\phi$Es.

\begin{definition}
Under homogeneity assumption within the clusters and  taking into account the weights $w_{hi}$, the (weighted) pseudo-maximum likelihood estimator (PMLE), $\widehat{\boldsymbol{\beta}}_P$, of $\boldsymbol{\beta}$ is obtained by maximizing
\begin{align}
\mathcal{L}\left( \boldsymbol{\beta}\right) &   =\sum \limits_{h=1}^{H}\sum\limits_{i=1}^{n_{h}}w_{hi}\log\boldsymbol{\pi}_{hi}^{T}\left( \boldsymbol{\beta}\right)\widehat{\boldsymbol{y}}_{hi}, \label{eq:likelihood}
\end{align}
where $\log  \boldsymbol{\pi}_{hi}^{T}\left( \boldsymbol{\beta}\right)=(\log \pi_{hi1}\left( \boldsymbol{\beta}\right), \dots,\log \pi_{hid+1}\left( \boldsymbol{\beta}\right) )$.
\end{definition}

The PMLE can be obtained as the solution of the system of equations $\boldsymbol{u}\left( \boldsymbol{\beta}\right) =\boldsymbol{0}_{d(k+1)}$, where $\boldsymbol{0}_{d(k+1)}$ is the null vector of dimension $d(k+1)$, and
\begin{align}
\boldsymbol{u}\left( \boldsymbol{\beta}\right) & =\sum\limits_{h=1}^{H}\sum \limits_{i=1}^{n_{h}}\boldsymbol{u}_{hi}\left( \boldsymbol{\beta}\right) ,  \quad  \boldsymbol{u}_{hi}\left( \boldsymbol{\beta}\right)  =w_{hi}\boldsymbol{r}_{hi}^{\ast}\left( \boldsymbol{\beta}\right) \otimes \boldsymbol{x}_{hi}, \label{eq:uhi}
\end{align}
where $\otimes$ is the Kronecker product  and $\boldsymbol{r}_{hi}^{\ast}\left( \boldsymbol{\beta}\right)  =\widehat{\boldsymbol{y}}_{hi}^{\ast}-m_{hi}\boldsymbol{\pi}_{hi}^{\ast }\left( \boldsymbol{\beta}\right) $, denoting with superscript $^{\ast}$ the vector obtained deleting the last component from the initial vector.

\begin{remark}\label{remark:over}
The distribution of \ $\widehat{\boldsymbol{Y}}_{hi}$, might be  unknown, as their components jointly, might be correlated. The most common assumption is to consider that \ $\widehat{\boldsymbol{Y}}_{hi}$ has a multinomial sampling scheme, which means that \ $\boldsymbol{Y}_{hij}$, $j=1,...,m_{hi}$ are independent random variables with covariance matrix
\begin{equation*}
\boldsymbol{\Sigma }_{hi} =m_{hi}\boldsymbol{\Delta }(\boldsymbol{\pi }_{hi}\left( \boldsymbol{\beta }\right) ),
\end{equation*}
with \ $\boldsymbol{\Delta }(\boldsymbol{\pi }_{hi}\left( \boldsymbol{\beta }\right)) =\mathrm{diag}(\boldsymbol{\pi }_{hi}\left( \boldsymbol{\beta }\right) )-\boldsymbol{\pi }_{hi}\left( \boldsymbol{\beta }\right) \boldsymbol{\pi }_{hi}^{T}\left( \boldsymbol{\beta }\right)$; and since (\ref{eq:likelihood}) is not an approximation, the term \textquotedblleft pseudo\textquotedblright\ should be dropped. A weaker assumption is to consider that \  $\widehat{\boldsymbol{Y}}_{hi}$ has a multinomial sampling scheme with a overdispersion parameter \ $\nu _{hi}=1+\rho _{hi}^{2}(m_{hi}-1)$,  and%
\begin{equation*}
\boldsymbol{\Sigma }_{hi}=\nu _{hi}m_{hi}\boldsymbol{\Delta }(\boldsymbol{\pi }_{hi}\left(\boldsymbol{\beta }\right) ),  
\end{equation*}
but the distribution of \ $\widehat{\boldsymbol{Y}}_{hi}$ is not in principle used for the estimators. Distributions such as Dirichlet Multinomial, Random Clumped and $m$-inflated belong to this family (see \cite{morel2012overdispersion, alonso2017}  and \cite{ Castilla2018_GIL} for details). In Appendix \ref{sec:algorithms}, the algorithms needed to compute these distributions in the context of PLR model with complex design are presented.  
\end{remark}

\begin{figure}[h!]
\centering
\begin{tabular}{ll}
\includegraphics[scale=0.2875]{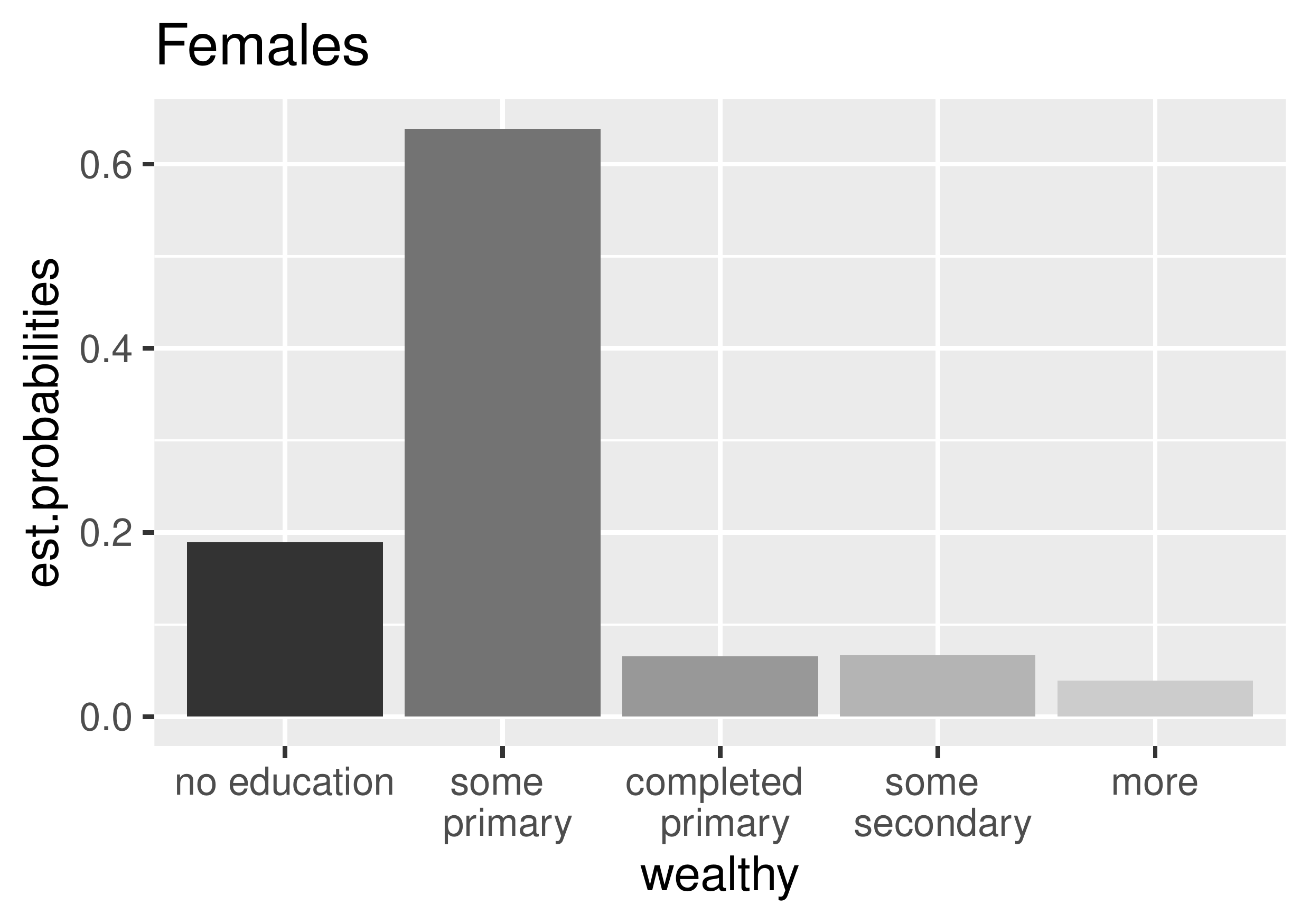}&
\includegraphics[scale=0.2875]{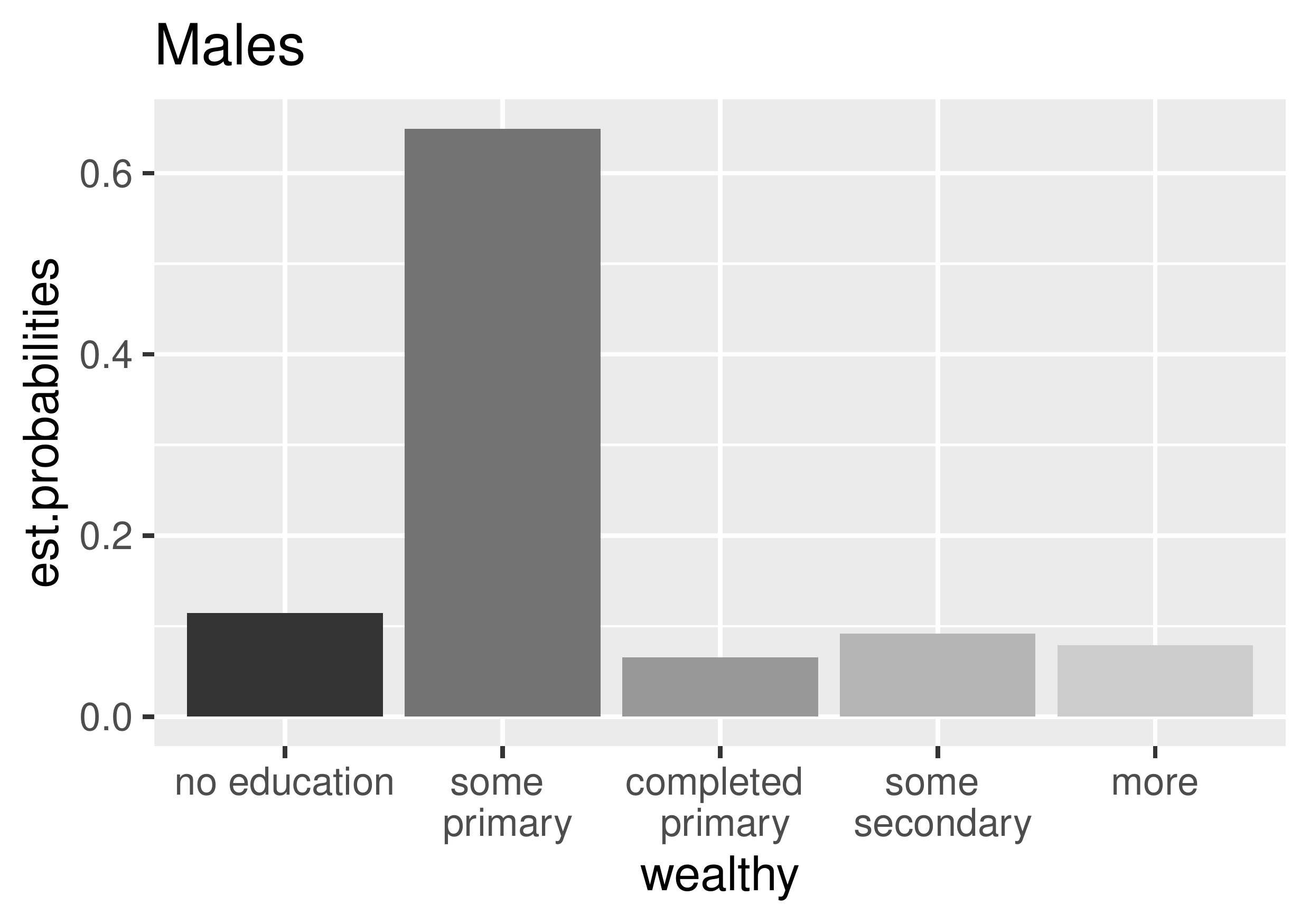}
\end{tabular}
\caption{Estimated probabilities for wealth quintiles by gender. PMLE. \label{fig:malawi_mle_ggplot}}
\end{figure}

\begin{example}[Education in Malawi] \label{example:malawi_mle}
The  2010 Malawi Demographic and Health Survey (2010 MDHS, \cite{Malawi2010})  was implemented by the National Statistical Office (NSO) from June through November 2010, with a nationally representative sample of more than $27,000$ households.  The sample for the 2010 MDHS was designed to provide population and health indicator estimation at the national, regional, and district levels. Let us focus on Tables 2.3.1 and 2.3.2 of the cited study, that present data on educational attainment for female and male household members age $6$ and older, divided in five wealth quintile levels, which are considered as strata. We consider here a response variable with $d+1=5$ categories: ``no education'', ``some primary'', ``completed primary'', ``some secondary'', ``completed secondary or more''. For simplicity, the missing observations are not taken into account. Figure \ref{fig:malawi_mle_ggplot} shows the estimated probabilities by the PMLE of each one of the response categories for each gender. As expected, the proportion of women who have never attended any formal schooling is greater than the proportion of men and the proportion of the population that has attained education declines with its level. In the ensuing work, we will present alternative estimators to the PMLE, which are seen to provide better performance in terms of robustness. 
\end{example}

\subsection{PM$\phi$Es: definition, estimation and asymptotic distribution.}

\begin{definition}
Given the probability vectors (\ref{eq:emp_prob_vector}) and (\ref{eq:theo_prob_vector}),  the family of phi-divergence measures between these probability vectors is given by
\begin{equation}
d_{\phi}\left( \widehat{\boldsymbol{p}},\boldsymbol{\pi}\left( \boldsymbol{\beta}\right) \right) =\frac{1}{\tau}\sum\limits_{h=1}^{H}\sum\limits_{i=1}^{n_{h}}w_{hi}m_{hi}\sum\limits_{s=1}^{d+1}\pi_{his}\left(\boldsymbol{\beta}\right) \phi\left( \frac{\widehat{y}_{his}}{m_{hi}\pi_{his}\left( \boldsymbol{\beta}\right) }\right) ,   \label{eq:phi_div_def}
\end{equation}
where $\phi\in\Phi^{\ast}$ and $\Phi^{\ast}$ denotes the class of all convex functions $\phi:[0,\infty)\rightarrow \mathbb{R}\cup \{\infty \}$, such that $\phi\left( 1\right) =0$, $\phi^{\prime\prime}\left( 1\right) >0$ and we define $0\phi\left(0/0\right) =0$ and $0\phi\left( p/0\right)=p\lim_{u\rightarrow\infty}\phi\left( u\right)/u$. 
\end{definition}

Notice that, for $\phi(x)=x\log x-x+1$ in (\ref{eq:phi_div_def}) , we have the  so-called Kullback Leibler divergence. For more details about phi-divergence measures see \cite{pardo2005}.

Let us now consider the PMLE in (\ref{eq:likelihood}). It can be shown that it is related to the Kullback-Leibler divergence measure as follows
\begin{align}
d_{KL}\left( \widehat{\boldsymbol{p}},\boldsymbol{\pi}\left(\boldsymbol{\beta}\right) \right) & =K-\frac{1}{\tau} \mathcal{L}\left( \boldsymbol{\beta}\right),  \notag
\end{align} 
with $K$  being a constant not dependent on $\boldsymbol{\beta}$ (\cite{Castilla2018_ASTA}). Therefore, the maximization of $\mathcal{L}\left( \boldsymbol{\beta}\right)$ is equivalent to the minimization of $d_{KL}\left( \widehat{\boldsymbol{p}},\boldsymbol{\pi}\left(\boldsymbol{\beta}\right) \right)$, i.e., PMLE is the one which minimizes the Kullback-Leibler divergence between the empirical and theoretical probability vectors, \ref{eq:emp_prob_vector} and \ref{eq:theo_prob_vector},
\begin{equation}
\widehat{\boldsymbol{\beta}}_P=\arg\min_{\boldsymbol{\beta\in\Theta}}d_{KL}\left( \widehat{\boldsymbol{p}},\boldsymbol{\pi}\left(\boldsymbol{\beta}\right) \right) .   \label{eq:kullback_estimators}
\end{equation}

The definition of PM$\phi$E arises from the idea of generalize definition (\ref{eq:kullback_estimators}) to other phi-divergence measures. 

\begin{definition}
\label{def:PMphiE} We consider the multinomial logistic regression model with
complex survey defined in (\ref{eq:pi}). The PM$\phi$E of $\boldsymbol{\beta}$ is defined as
\begin{equation*}
\widehat{\boldsymbol{\beta}}_{\phi,P}=\arg\min_{\boldsymbol{\beta}\in\Theta}d_{\phi}\left( \widehat{\boldsymbol{p}},\boldsymbol{\pi}\left( \boldsymbol{\beta}\right) \right).
\end{equation*}
\end{definition}

\bigskip

Once the PM$\phi$Es are defined, it is necessary to provide the equations needed to obtain them. From equation (\ref{eq:phi_div_def}), it is clear that  the PM$\phi$E of $\boldsymbol{\beta}$, $\widehat{\boldsymbol{\beta }}_{\phi,P}$, is obtained by solving the system of equations $\boldsymbol{u}_{\phi }\left( \boldsymbol{\beta }\right)=\boldsymbol{0}_{d(k+1)}$, where

\begin{equation}
\boldsymbol{u}_{\phi }\left( \boldsymbol{\beta }\right)=\sum\limits_{h=1}^{H}\sum\limits_{i=1}^{n_{h}}\boldsymbol{u}_{\phi,hi}\left( \boldsymbol{\beta }\right) ,  
\end{equation}%
with
\begin{align*}
\boldsymbol{u}_{\phi ,hi}\left( \boldsymbol{\beta }\right) 
& =\frac{w_{hi}m_{hi}}{\phi ^{\prime \prime }(1)}\frac{\partial \boldsymbol{\pi }_{hi}^{T}(\boldsymbol{\beta })}{\partial \boldsymbol{\beta }}\boldsymbol{f}_{\phi ,hi}(\tfrac{\widehat{\boldsymbol{y}}_{hi}}{m_{hi}},\boldsymbol{\beta }), \quad \frac{\partial \boldsymbol{\pi }_{hi}^{T}(\boldsymbol{\beta })}{\partial\boldsymbol{\beta }} =\left( \boldsymbol{I}_{d\times d},\boldsymbol{0}_{d\times 1}\right) \boldsymbol{\Delta }(\boldsymbol{\pi }_{hi}\left(\boldsymbol{\beta }\right) )\otimes \boldsymbol{x}_{hi},
\end{align*}
and

\begin{align*}
\boldsymbol{f}_{\phi ,hi}(\tfrac{\widehat{\boldsymbol{y}}_{hi}}{m_{hi}},\boldsymbol{\beta })& =(f_{\phi ,hi1}(\tfrac{\widehat{y}_{hi1}}{m_{hi}},\boldsymbol{\beta }),...,f_{\phi ,hi(d+1)}(\tfrac{\widehat{y}_{hi(d+1)}}{m_{hi}},\boldsymbol{\beta }))^{T},  \notag \\
f_{\phi ,his}(x,\boldsymbol{\beta })& =\frac{x}{\pi _{his}(\boldsymbol{\beta})}\phi^{\prime }\left( \frac{x}{\pi _{his}(\boldsymbol{\beta })}\right)-\phi \left( \frac{x}{\pi _{his}(\boldsymbol{\beta })}\right) .
\end{align*}

Theorem  \ref{th:Th1} establishes the asymptotic distribution of the PM$\phi$Es, which will be the basis of the definition of the family of Wald-type tests in Section \ref{Sec:Wald}. The proof of this theorem can be found in \cite{Castilla2018_ASTA}.

\begin{theorem}
\label{th:Th1}Let $\widehat{\boldsymbol{\beta}}_{\phi,P}$ the PM$\phi$E of parameter $\boldsymbol{\beta}$ for a multinomial logistic regression model with complex survey, $n$ the total of clusters in all the strata of the sample and $\eta _{h}^{\ast}$ an unknown proportion obtained as $\lim_{n\rightarrow\infty }\frac{n_{h}}{n}=\eta_{h}^{\ast}$, $h=1,...,H$.
Then, we have
\begin{equation}\label{eq:asymEst}
\sqrt{n}(\widehat{\boldsymbol{\beta}}_{\phi,P}-\boldsymbol{\beta}^0)\overset{\mathcal{L}}{\underset{n\mathcal{\rightarrow}\infty }{\longrightarrow}}\mathcal{N}\left( \boldsymbol{0}_{d(k+1)},\mathbf{J}^{-1}\left(\boldsymbol{\beta}^0\right) \mathbf{G}\left( \boldsymbol{\beta}^0\right)\mathbf{J}^{-1}\left( \boldsymbol{\beta}^0\right) \right) ,
\end{equation}
where $\boldsymbol{\beta}^0$ is the true parameter value  and
\begin{align*}
\mathbf{J}\left( \boldsymbol{\beta}\right) &=\lim_{n\rightarrow\infty }\mathbf{J}_{n}\left( \boldsymbol{\beta}\right) =\sum\limits_{h=1}^{H}\eta_{h}^{\ast}\lim_{n_{h}\rightarrow\infty}\mathbf{J}_{n_{h}}^{(h)}\left( \boldsymbol{\beta}\right), \\
\mathbf{G}\left(\boldsymbol{\beta }\right) &=\lim_{n\rightarrow\infty}\mathbf{G}_{n}\left(\boldsymbol{\beta }\right) =\sum\limits_{h=1}^{H}\eta_{h}^{\ast}\lim_{n_{h}\rightarrow\infty}\mathbf{G}_{n_{h}}^{(h)}\left( \boldsymbol{\beta}\right) ,
\end{align*}
with%
\begin{align*}
\mathbf{J}_{n}\left( \boldsymbol{\beta}\right) &=\frac{1}{n}\sum\limits_{h=1}^{H}\sum \limits_{i=1}^{n_{h}}w_{hi}m_{hi}\boldsymbol{\Delta}(\boldsymbol{\pi}_{hi}^{\ast}\left( \boldsymbol{\beta}\right) )\otimes\boldsymbol{x}_{hi}\boldsymbol{x}_{hi}^{T},\\
\mathbf{J}_{n_{h}}^{(h)}\left( \boldsymbol{\beta}\right) &=\frac{1}{n_{h}}\sum\limits_{i=1}^{n_{h}}w_{hi}m_{hi}\boldsymbol{\Delta}(\boldsymbol{\pi}_{hi}^{\ast}\left( \boldsymbol{\beta}\right) )\otimes\boldsymbol{x}_{hi}\boldsymbol{x}_{hi}^{T},\\
\mathbf{G}_{n}\left( \boldsymbol{\beta}\right) &=\frac{1}{n}\sum\limits_{h=1}^{H}\sum \limits_{i=1}^{n_{h}}\mathbf{V}[\boldsymbol{U}_{hi}\left( \boldsymbol{\beta}\right) ],\\
\mathbf{G}_{n_{h}}^{(h)}\left( \boldsymbol{\beta}\right) &=\frac{1}{n_{h}}\sum\limits_{i=1}^{n_{h}}\mathbf{V}[\boldsymbol{U}_{hi}\left( \boldsymbol{\beta}\right) ],\text{\quad}\mathbf{V}[\boldsymbol{U}_{hi}\left(\boldsymbol{\beta }\right) ]=w_{hi}^{2}\mathbf{V}[\widehat{\boldsymbol{Y}}_{hi}^{\ast }]\otimes\boldsymbol{x}_{hi}\boldsymbol{x}_{hi}^{T},
\end{align*}
$\mathbf{J}\left( \boldsymbol{\beta}\right) $ is the Fisher information matrix, $\mathbf{V}[\boldsymbol{\cdot}]$ denotes the variance-covariance matrix of a random vector and $\boldsymbol{U}_{hi}\left( \boldsymbol{\beta }\right) $ is the random variable generator of $\boldsymbol{u}_{hi}\left(\boldsymbol{\beta}\right) $, given by (\ref{eq:uhi}).
\end{theorem}

\begin{remark}
Matrices $\mathbf{J}(\boldsymbol{\beta}^0)$ and $\mathbf{G}(\boldsymbol{\beta}^0)$ of Theorem \ref{th:Th1} can be consistently estimated as
\begin{align*}
\widehat{\mathbf{J}}_n(\widehat{\boldsymbol{\beta}}_{\phi,P})&=\frac{1}{n}\sum\limits_{h=1}^{H}\sum \limits_{i=1}^{n_{h}}w_{hi}m_{hi}\boldsymbol{\Delta}(\boldsymbol{\pi}_{hi}^{\ast}(\widehat{\boldsymbol{\beta}}_{\phi,P})) \otimes\boldsymbol{x}_{hi}\boldsymbol{x}_{hi}^{T}\\
\widehat{\mathbf{G}}_n(\widehat{\boldsymbol{\beta}}_{\phi,P})&=\frac{1}{n}\sum\limits_{h=1}^{H}\sum \limits_{i=1}^{n_{h}}\left(\boldsymbol{u}_{hi}(\widehat{\boldsymbol{\beta}}_{\phi,P})-\frac{1}{n}\boldsymbol{u}(\widehat{\boldsymbol{\beta}}_{\phi,P}) \right)\left(\boldsymbol{u}_{hi}(\widehat{\boldsymbol{\beta}}_{\phi,P})-\frac{1}{n}\boldsymbol{u}(\widehat{\boldsymbol{\beta}}_{\phi,P}) \right)^T.
\end{align*}
\end{remark}

\bigskip
An important family of phi-divergence measures is obtained by restricting $\phi$ from the family of convex functions $\Phi^*$ to the Cressie-Read subfamily, that is to say, $\phi$ is of the following form:
\begin{equation*}
\phi_{\lambda}(x)=\left\{
\begin{array}{ll}
\frac{1}{\lambda(1+\lambda)}\left[ x^{\lambda+1}-x-\lambda(x-1)\right] , & \lambda\in\mathbb{R}-\{-1,0\} \\
\lim_{\upsilon\rightarrow\lambda}\frac{1}{\upsilon(1+\upsilon)}\left[x^{\upsilon+1}-x-\upsilon(x-1)\right] , & \lambda\in\{-1,0\}
\end{array}
\right. .
\end{equation*}
For the Cressie-Read subfamily, it is established that for $\lambda\neq-1$, $$\boldsymbol{u}_{\phi_{\lambda}}\left( \boldsymbol{\beta}\right) =\sum \nolimits_{h=1}^{H}\sum \nolimits_{i=1}^{n_{i}}\boldsymbol{u}_{\phi_{\lambda},hi}\left( \boldsymbol{\beta}\right), $$ where%
\begin{equation}
\boldsymbol{u}_{\phi_{\lambda},hi}\left( \boldsymbol{\beta}\right) =\frac{w_{hi}}{(\lambda+1)m_{hi}^{\lambda}}\boldsymbol{\Delta}^*(\boldsymbol{\pi}_{hi}(\boldsymbol{\beta}))\mathrm{diag}^{-(\lambda+1)}(\boldsymbol{\pi}_{hi}(\boldsymbol{\beta}))\widehat{\boldsymbol{y}}_{hi}^{\lambda+1} \otimes \boldsymbol{x}_{hi}, \label{eq:uphi}
\end{equation}
where  $$\boldsymbol{\Delta}^*(\boldsymbol{\pi}_{hi}(\boldsymbol{\beta}))=(\boldsymbol{I}_d,\boldsymbol{0}_d)\boldsymbol{\Delta}(\boldsymbol{\pi}_{hi}(\boldsymbol{\beta})).$$
We can observe that for $\lambda=0$, we have
\begin{equation*}
\phi_{\lambda=0}(x)=\lim_{\upsilon\rightarrow0}\frac{1}{\upsilon(1+\upsilon )}\left[ x^{\upsilon+1}-x-\upsilon(x-1)\right] =x\log x-x+1,
\end{equation*}
and the associated phi-divergence, coincides with the Kullback divergence. Therefore, the PM$\phi$E of $\boldsymbol{\beta}$ based on $\phi_{\lambda}(x)$ contains as special case the PMLE and $\boldsymbol{u}_{hi}\left( \boldsymbol{\beta}\right) $ given in (\ref{eq:uhi}) matches $\boldsymbol{u}_{\phi,hi}\left( \boldsymbol{\beta}\right) $ given in (\ref{eq:uphi}).  Other important divergences are obtained inside this family: for $\lambda=1$ the  chi-square divergence, for $\lambda=2/3$ the Cressie-Read divergence and for $\lambda=-0.5$ the  Hellinger distance. Note that the Hellinger distance is well-known in statistical theory for its robustness (\cite{lindsay1994}). The robustness of (Cressie-Read) PM$\phi$Es  for $-1<\lambda<0$ is proved, through the study of their influence function, in Appendix \ref{sec:modelrobust}. 


\begin{remark}\label{remark:simple_design}
Along this paper, we are referring to the case of complex sample survey. These all procedures  can be easily  simplified to the case of simple sample survey by considering a single stratum and considering ``observations'' instead of clusters. Some work has been done within the phi-divergence measures and multinomial logistic regression (see \cite{gupta2008, martin2014}) but, to the best of our knowledge, the robustness issue was not previously considered.  In Section \ref{sec:mammo}, an example is provided to illustrate the application of the  proposed methods also in this context. 
\end{remark}

\section{Robust Wald-type tests \label{Sec:Wald}}
In the last years, it has been very common in the statistical literature to consider Wald-type tests based on the minimum distance estimators instead of the MLE. The resulting tests have an excellent behaviour in relation to the robustness with a non-significant loss of efficiency, see for instance, \cite{basu2017_EJS, basu2018_METRIKA}. In this section, we will introduce Wald-type test statistics based on  the PM$\phi$Es  as a generalization of the classical Wald test based on PMLE. As it happened with  PM$\phi$Es, the robustness of these proposed Wald-type tests can be proved for $-1<\lambda<0$ (see Appendix \ref{sec:modelrobust}).

In this context, we are interested in testing
\begin{equation}
H_0:\boldsymbol{M}^T\boldsymbol{\beta}=\boldsymbol{m} \quad \text{against} \quad \boldsymbol{M}^T\boldsymbol{\beta}\neq \boldsymbol{m},
\label{eq:contrast}
\end{equation}
where  $\boldsymbol{M}$ is $d(k+1)\times r$ full rank matrix with $r \leq d(k + 1)$ and $\boldsymbol{m}$ an $r$-vector.

\begin{definition}
 Let $\widehat{\boldsymbol{\beta}}_{\phi,P}$ the PM$\phi$E of $\boldsymbol{\beta}$ and denote $$\widehat{\mathbf{V}}_{n}(\widehat{\boldsymbol{\beta}}_{\phi,P})=\widehat{\mathbf{J}}_{n}^{-1}(\widehat{\boldsymbol{\beta}}_{\phi,P})\widehat{\mathbf{G}}_{n}(\widehat{\boldsymbol{\beta}}_{\phi,P})\widehat{\mathbf{J}}_{n}^{-1}(\widehat{\boldsymbol{\beta}}_{\phi,P}).$$ Then, the
family of Wald-type test statistics for testing the null hypothesis given in (\ref{eq:contrast}) is defined as
\begin{equation}
W_{n}(\widehat{\boldsymbol{\beta}}_{\phi,P})=n\left( \boldsymbol{M}^T\widehat{\boldsymbol{\beta}}_{\phi,P}-\boldsymbol{m}\right)^T\left[\boldsymbol{M}^T\mathbf{V}_{n}(\widehat{\boldsymbol{\beta}}_{\phi,P})\boldsymbol{M} \right]^{-1}\left( \boldsymbol{M}^T\widehat{\boldsymbol{\beta}}_{\phi,P}-\boldsymbol{m}\right).
\label{eq:Wald_phi}
\end{equation}
\end{definition}

\begin{theorem}\label{th:wald_test_phi}
The asymptotic distribution of the Wald-type
test statistics, $W_{n}(\widehat{\boldsymbol{\beta}}_{\phi,P})$, under the null hypothesis in (\ref{eq:contrast}), is a chi-square distribution with $r$ degrees of freedom.
\end{theorem}

\begin{corollary}
Based on Theorem \ref{th:wald_test_phi}, the null hypothesis in (\ref{eq:contrast}) will be rejected if 
\begin{align}\label{eq:wald_test_phi}
W_{n}(\widehat{\boldsymbol{\beta}}_{\phi,P})>\chi^2_{r,\alpha}
\end{align}
being $\chi^2_{r,\alpha}$ the upper $\alpha$-th quantile of $\chi^2_{r}$.
\end{corollary}

\bigskip

The following theorem may be used to approximate the power function for the Wald-
type test statistics given in (\ref{eq:wald_test_phi}).

\begin{theorem}\label{th:power1}
Let $\boldsymbol{\beta}^{0}$ be the true value of the parameter and let us denote
\begin{equation*}
\ell ^{\ast }\left( \boldsymbol{\beta }_{1}\mathbf{,}\boldsymbol{\beta }_{2}\right) =\left( \boldsymbol{M}^{T}\boldsymbol{\beta }_{1}-\boldsymbol{m}\right) ^{T}\left( \boldsymbol{M}^{T}\mathbf{V}\left( \boldsymbol{\beta }_{2}\right) \boldsymbol{M}\right) ^{-1}\left( \boldsymbol{M}^{T}\boldsymbol{\beta }_{1}-\boldsymbol{m}\right) .
\end{equation*}%
Then it holds 
\begin{equation*}
\sqrt{n}\left( \ell ^{\ast }(\widehat{\boldsymbol{\beta }}_{\phi ,P},\boldsymbol{\beta ^{0 }})-\ell ^{\ast }\left( \boldsymbol{\beta ^{0},\beta ^{0 }}\right) \right) \underset{n\rightarrow \infty }{\overset{L}{\longrightarrow }}\mathcal{N}(0,\sigma _{W}^{2}\left( \boldsymbol{\beta^{0}}\right) ),
\end{equation*}
where 
\begin{equation*}
\sigma _{W}^{2}\left( \boldsymbol{\beta ^{0}}\right) =4\left( \boldsymbol{M}^{T}\boldsymbol{\beta }^{0}-\boldsymbol{m}\right) ^{T}\left( \boldsymbol{M}^{T}\mathbf{V}\left( \boldsymbol{\beta }^{0}\right) \boldsymbol{M}\right) ^{-1}\left( \boldsymbol{M}^{T}\boldsymbol{\beta }^{0}-\boldsymbol{m}\right).
\end{equation*}
\end{theorem}

\begin{theorem} \label{Th:PowerApprox_Wald} 
Let $\boldsymbol{\beta }^{0}\in \Theta$, with $\boldsymbol{M}^{T}\boldsymbol{\beta }^{0}\neq \boldsymbol{m}$, be the true value of the parameter such that $\widehat{\boldsymbol{\beta }}_{\phi,P }\underset{n\rightarrow \infty }{\overset{P}{\longrightarrow }}\boldsymbol{\beta }^{0}$. The power function of the Wald-type test  given in (\ref{eq:wald_test_phi}),  is given by 
\begin{equation}
\Pi _{W_{n}(\widehat{\boldsymbol{\beta }}_{\phi ,P})}\left( \boldsymbol{\beta }^{0}\right) =1-\Phi_{n}\left( \frac{1}{\sigma _{W}\left( \boldsymbol{\beta }^{0}\right) }\left( \frac{\chi _{r,\alpha }^{2}}{\sqrt{n}}-\sqrt{n} \ \ell ^{\ast }(\widehat{\boldsymbol{\beta }}_{\phi ,P},\boldsymbol{\beta^{0}})\right) \right)
\end{equation}%
where $\Phi _{n}\left( x\right) $ uniformly tends to the standard normal distribution as $n\rightarrow \infty$. 
\end{theorem}

\begin{corollary}
It is clear that 
\begin{equation*}
\lim_{n\rightarrow \infty }\Pi_{W_{n}(\widehat{\boldsymbol{\beta }}_{\phi ,P})}\left( \boldsymbol{\beta}^{0}\right)=1
\end{equation*}%
for all $\alpha \in \left(0,1\right)$. Therefore, the Wald-type tests are consistent in the sense of Fraser.
\end{corollary}

\begin{remark}
Theorem \ref{Th:PowerApprox_Wald} can be applied in the sense of getting the necessary sample size in order to get that the Wald-type tests have a determinate fix power, i.e., $\Pi_{W_{n}(\widehat{\boldsymbol{\beta }}_{\phi ,P})}\left( \boldsymbol{\beta}^{0}\right) \equiv \pi^{0}$ and size $\alpha$. The necessary sample size is given by 
\begin{equation*}
n=\left[ \frac{A+B+\sqrt{A(A+2B)}}{2\ell ^{\ast }\left( \boldsymbol{\beta}^{0},\boldsymbol{\beta}^{0}\right) ^{2}}\right] +1,
\end{equation*}%
where $\left[ x\right] $ denotes the largest integer less than or equal to $x$, $A=\sigma _{W}^{2}\left( \boldsymbol{\beta }^{0}\right) \left( \Phi^{-1}(1-\pi ^{0})\right) ^{2}$ and $B=2\ell ^{\ast }\left( \boldsymbol{\beta }^{0},\boldsymbol{\beta }^{0}\right) \chi _{r,\alpha }^{2}$.
\end{remark}

\bigskip
We may also find  approximations of the power function of the Wald-type tests given in (\ref{eq:Wald_phi}) at an alternative hypothesis close to the null hypothesis.  Let $\boldsymbol{\beta}_n \in \Theta - \Theta_0$ be a given
alternative, and let  $\boldsymbol{\beta}_0 \in \Theta_0$ (null hypothesis) the element closest to $\boldsymbol{\beta}_n$ in terms of the Euclidean
distance. We may introduce contiguous alternative hypotheses by considering a fixed $\boldsymbol{d}\in \mathbb{R}^{d(k+1)}$ and to permit $\boldsymbol{\beta }_{n}$ moving towards $\boldsymbol{\beta }_{0}$ as $n$ increases through the relation

\begin{equation}
H_{1,n}:\boldsymbol{\beta }_{n}=\boldsymbol{\beta }_{0}+n^{-1/2}\boldsymbol{d}. 
\label{eq:alt1}
\end{equation}
Let us now relax the condition $\boldsymbol{M}^{T}\boldsymbol{\beta }_{0}=\boldsymbol{m}$ defining the null hypothesis. Let $\boldsymbol{\delta }\in \mathbb{R}^{r}$ and consider the following sequence, $\boldsymbol{\beta }_{n}$,  of parameters moving towards $\boldsymbol{\beta }_{0}$ according to 

\begin{equation}
H_{1,n}^{\ast }:\boldsymbol{M}^{T}\boldsymbol{\beta }_{n}-\boldsymbol{m}=n^{-1/2}\boldsymbol{\delta }. \label{eq:alt2}
\end{equation}


\begin{theorem}\label{th:respower_ult}
The asymptotic distribution of \ $W_{n}(\widehat{\boldsymbol{\beta }}_{\phi ,P})$ is given by:
\begin{itemize}
\item [(a)] Under $H_{1,n}$, $W_{n}(\widehat{\boldsymbol{\beta }}_{\phi ,P})\underset{n\rightarrow \infty }{\overset{\mathcal{L}}{\longrightarrow }}\chi_{r}^{2}\left( \Delta \right) $, where $\Delta$ is the parameter of non-centrality given by
$$\Delta=\boldsymbol{d}^{T}\boldsymbol{M}\left[ \boldsymbol{M}^{T}\mathbf{V}(\boldsymbol{\beta }_{0})\boldsymbol{M}\right] ^{-1}\boldsymbol{M}^{T}\boldsymbol{d},$$

\item [(b)] Under $H^*_{1,n}$,  $W_{n}(\widehat{\boldsymbol{\beta }}_{\phi ,P})\underset{n\rightarrow \infty }{\overset{\mathcal{L}}{\longrightarrow }}\chi_{r}^{2}\left( \Delta^* \right) $, where $\Delta^*$ is the parameter of non-centrality given by

$$\Delta^*=\boldsymbol{d}^{T}\boldsymbol{M}\left[ \boldsymbol{M}^{T}\mathbf{V}(\boldsymbol{\beta }_{0})\boldsymbol{M}\right] ^{-1}\boldsymbol{M}^{T}\boldsymbol{d}.$$
\end{itemize}

\end{theorem}

Proofs of the results given in this section can be found in Appendix \ref{app:proofs}.

\clearpage
\section{Monte Carlo Simulation Study \label{sec:mc}}

In this section, we develop a simulation study in order to illustrate the robustness of the proposed estimators and Wald-type tests based on them. Following the simulation studies proposed in \cite{castilla2018_BIOMETRICS} and \cite{castilla2019_ABHIK}, we consider $H=4$ strata  with $n_h$ clusters of the same size $m$, with $m=20$ and $n_h \in \{ 10,20,..,60\}$ for $h=1,\dots,H$. We consider  $d+1 = 3$ categories on the response variable, depending on $k = 2$ explanatory variables. The response variable $\widehat{\boldsymbol{Y}}_{hi}$, described as%

\begin{align*}
\boldsymbol{E}[\widehat{\boldsymbol{Y}}_{hi}]  &  =m\boldsymbol{\pi}_{hi}\left(\boldsymbol{\beta}^{0}\right)  \quad\text{and}\quad\boldsymbol{V}[\widehat{\boldsymbol{Y}}_{hi}]=\nu_{m}m\boldsymbol{\Delta}(\boldsymbol{\pi}_{hi}\left(  \boldsymbol{\beta}^{0}\right)  ),\\
\nu_{m}  &  =1+\rho^{2}(m-1), \quad i=1,\dots,n_h, \ h=1,\dots,H,
\end{align*}
is considered to follow the m-Inflated  multinomial distribution (see Remark \ref{remark:over}), with  parameters $\rho^2=0.5$ and $\boldsymbol{\pi}_{hi}\left(  \boldsymbol{\beta}^{0}\right)$, given by the logistic relationship (\ref{eq:pi}) with 

$$\boldsymbol{\beta}^0=(\beta_{01},\beta_{11},\beta_{21},\beta_{02},\beta_{12},\beta_{22})^T=(0,-0.9,0.1,0.6,-1.2,0.8)^T$$ and $\boldsymbol{x}_{hi}\overset{iid}{\sim}\mathcal{N}(\boldsymbol{0},\boldsymbol{I})$ for all $i=1, \ldots,n_h$, $h=1,\dots,H$. In order to study the robustness issue, these simulations are repeated under contaminated
data having $10\%$ outliers. These outliers are generated by permuting the elements of the outcome variable, such that categories 1, 2, 3 are classified as categories 3, 1, 2 for the outlying observations. Note that this view of considering outliers as classification errors  in the PLR model is, in fact, in line with the general literature on robust analysis of categorical data (\cite{johnson1985influence,croux2003a}) and is covered with the theory developed in Appendix \ref{sec:modelrobust}, where  our ``outlier producing measure''  indeed provides  classification error if the outlier point  yields its mass in a wrong category (see \cite{castilla2019_ABHIK} for more details). 

\begin{figure}[p]
\centering
\begin{tabular}{rc}
 \includegraphics[scale=0.4]{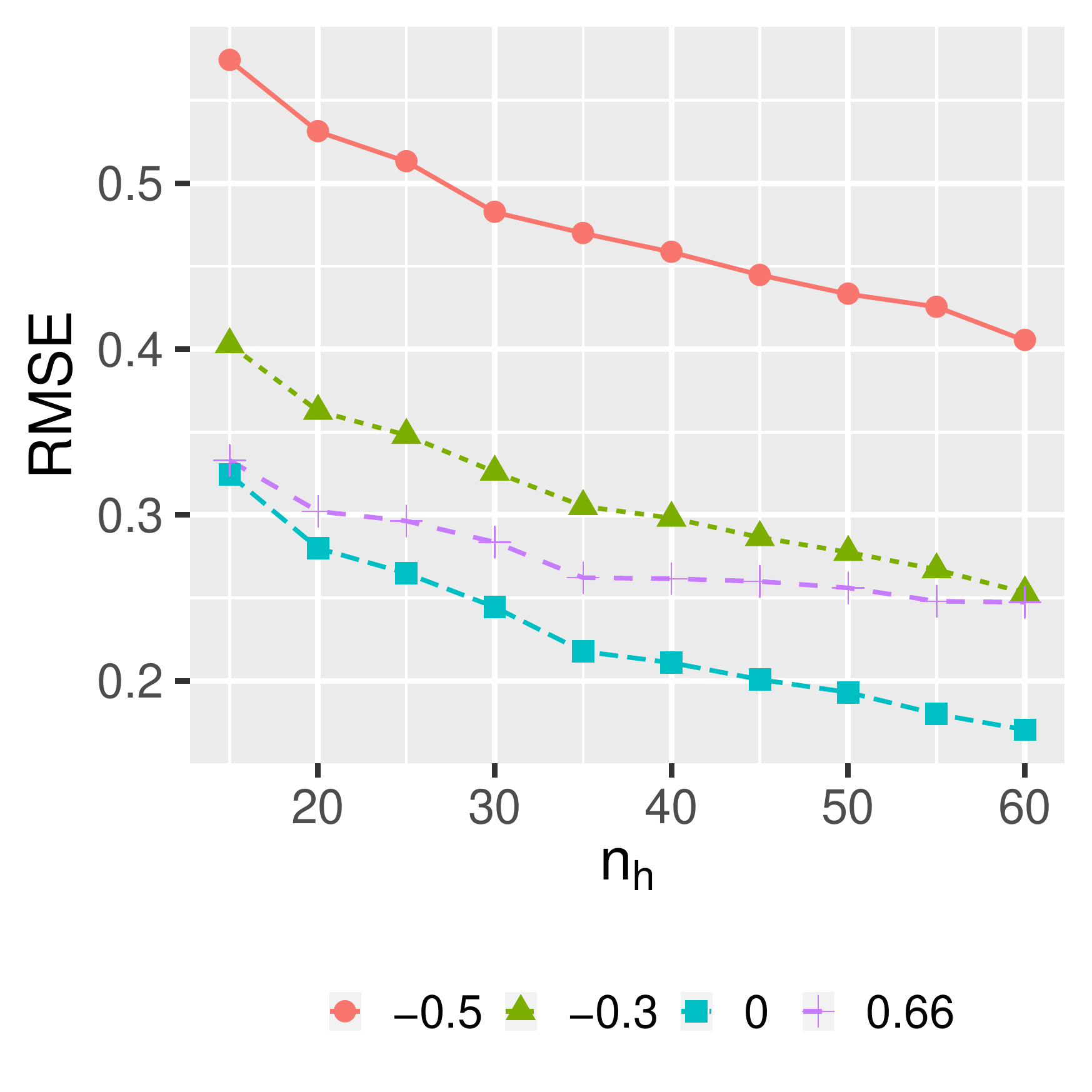}&
\includegraphics[scale=0.4]{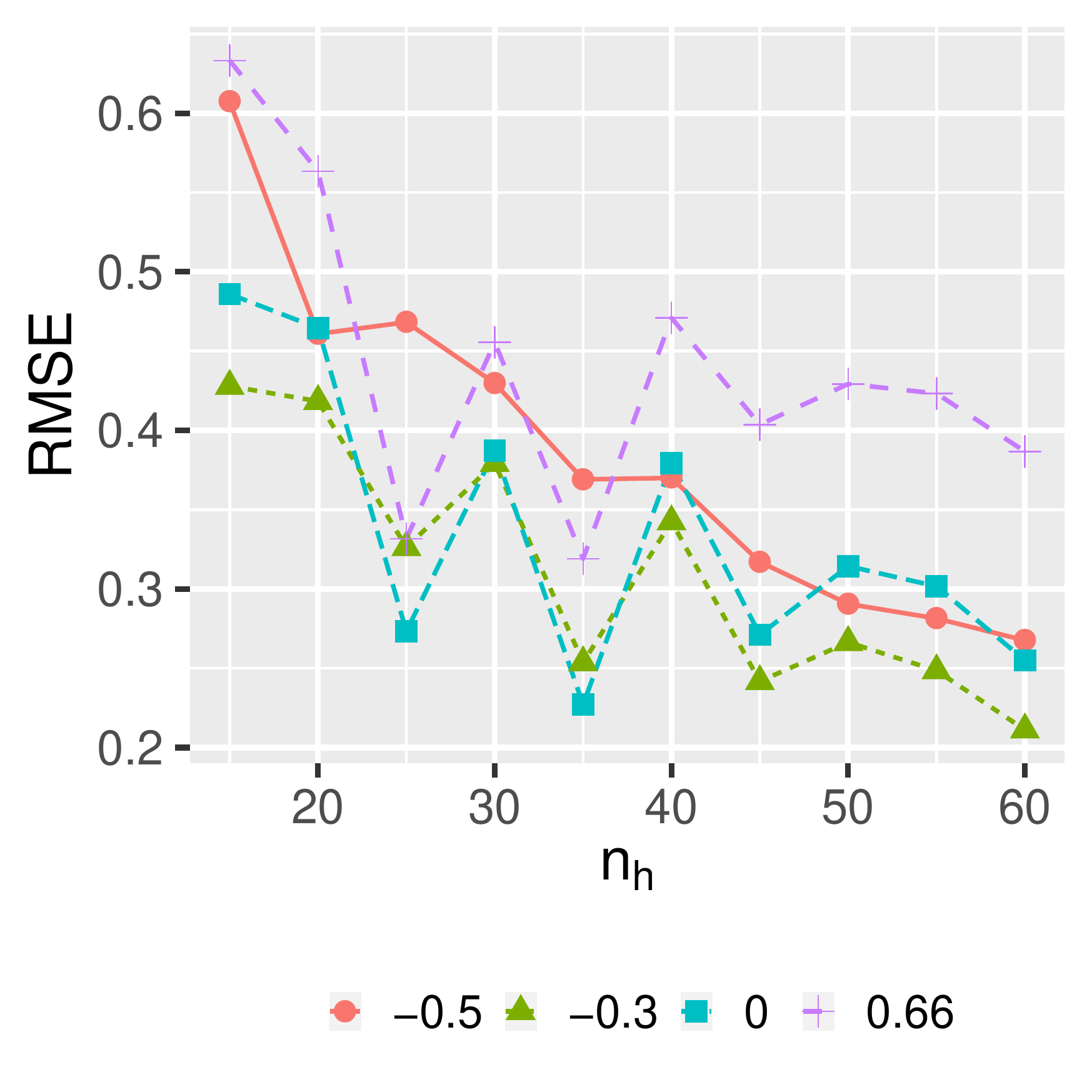}\\
 \includegraphics[scale=0.4]{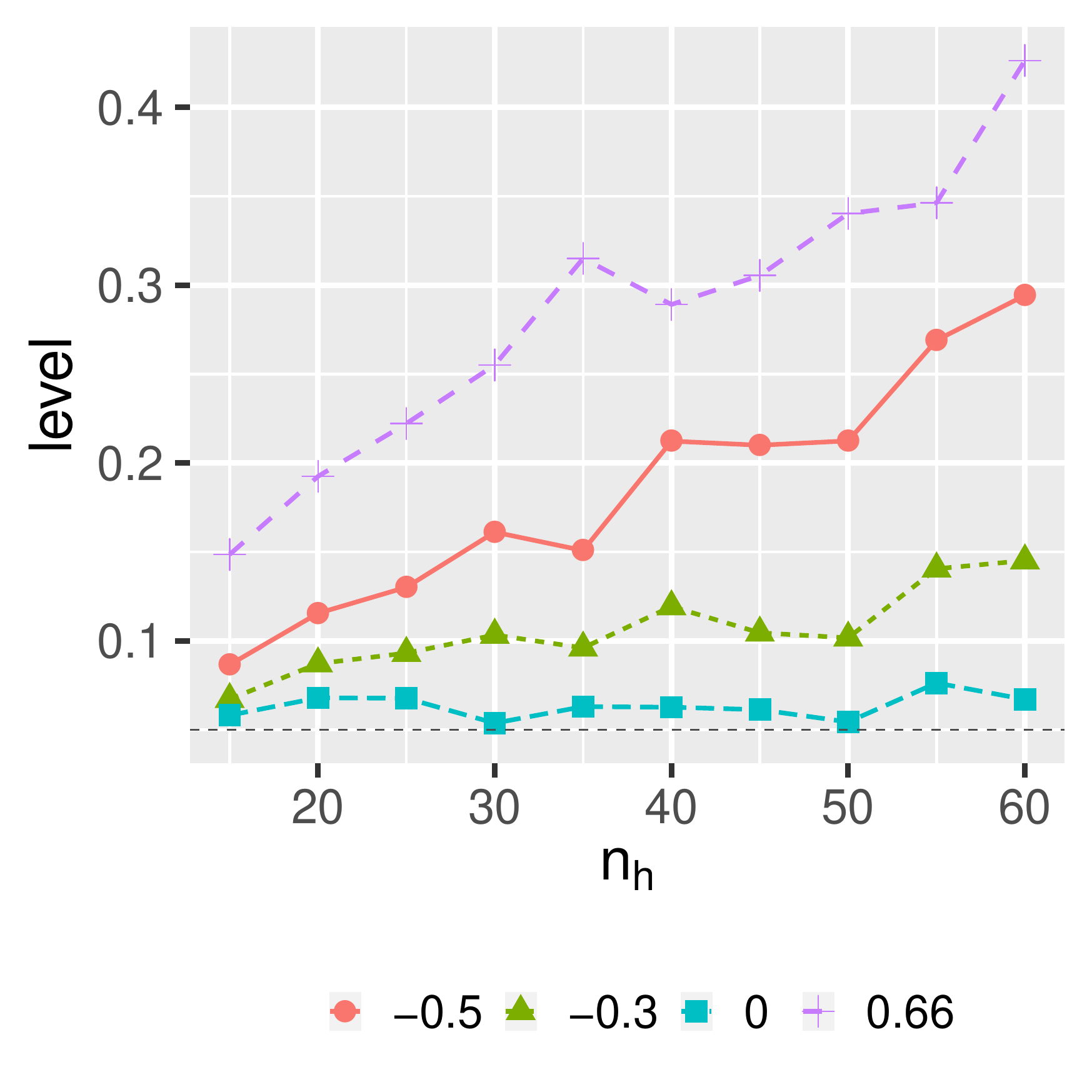}&
\includegraphics[scale=0.4]{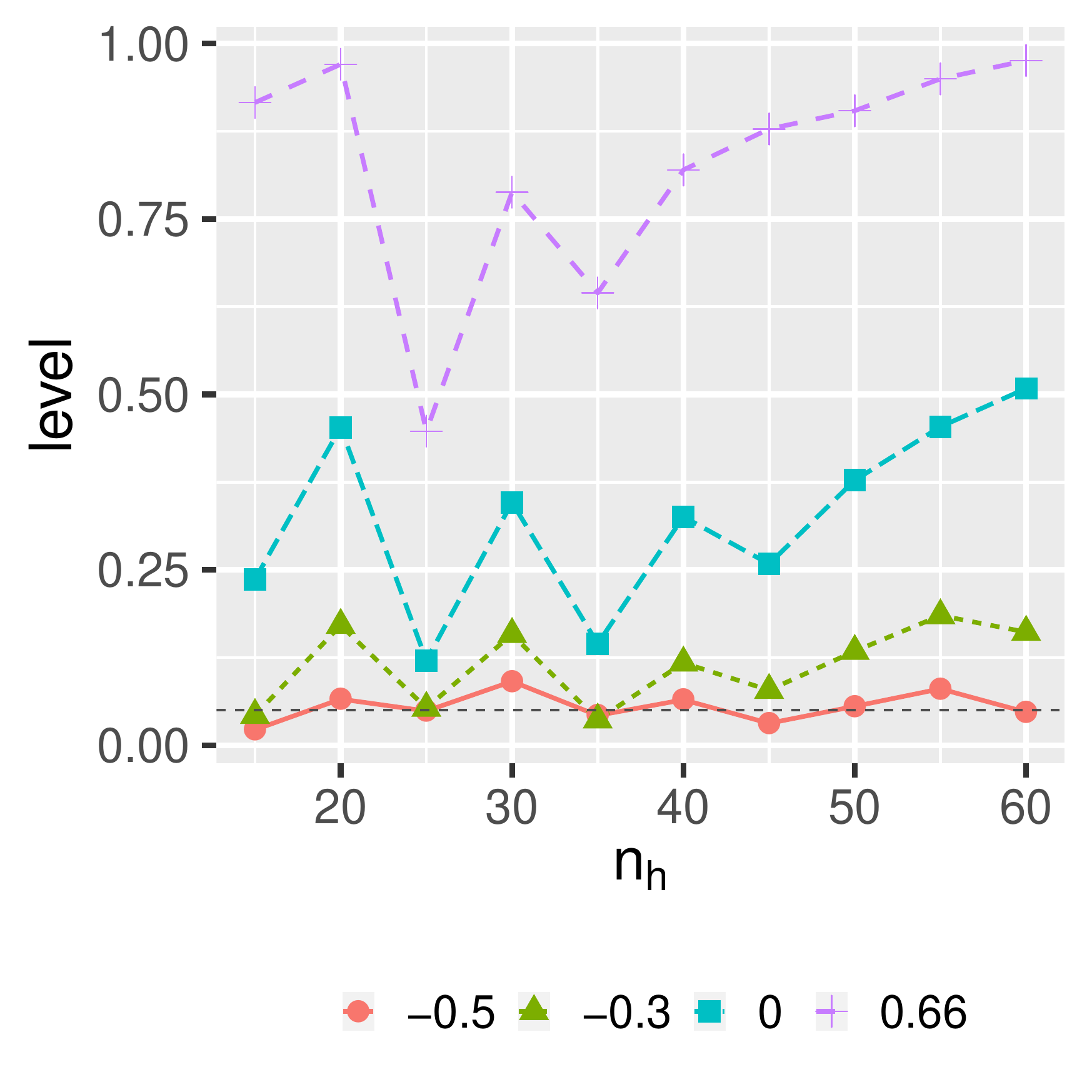}\\
 \includegraphics[scale=0.4]{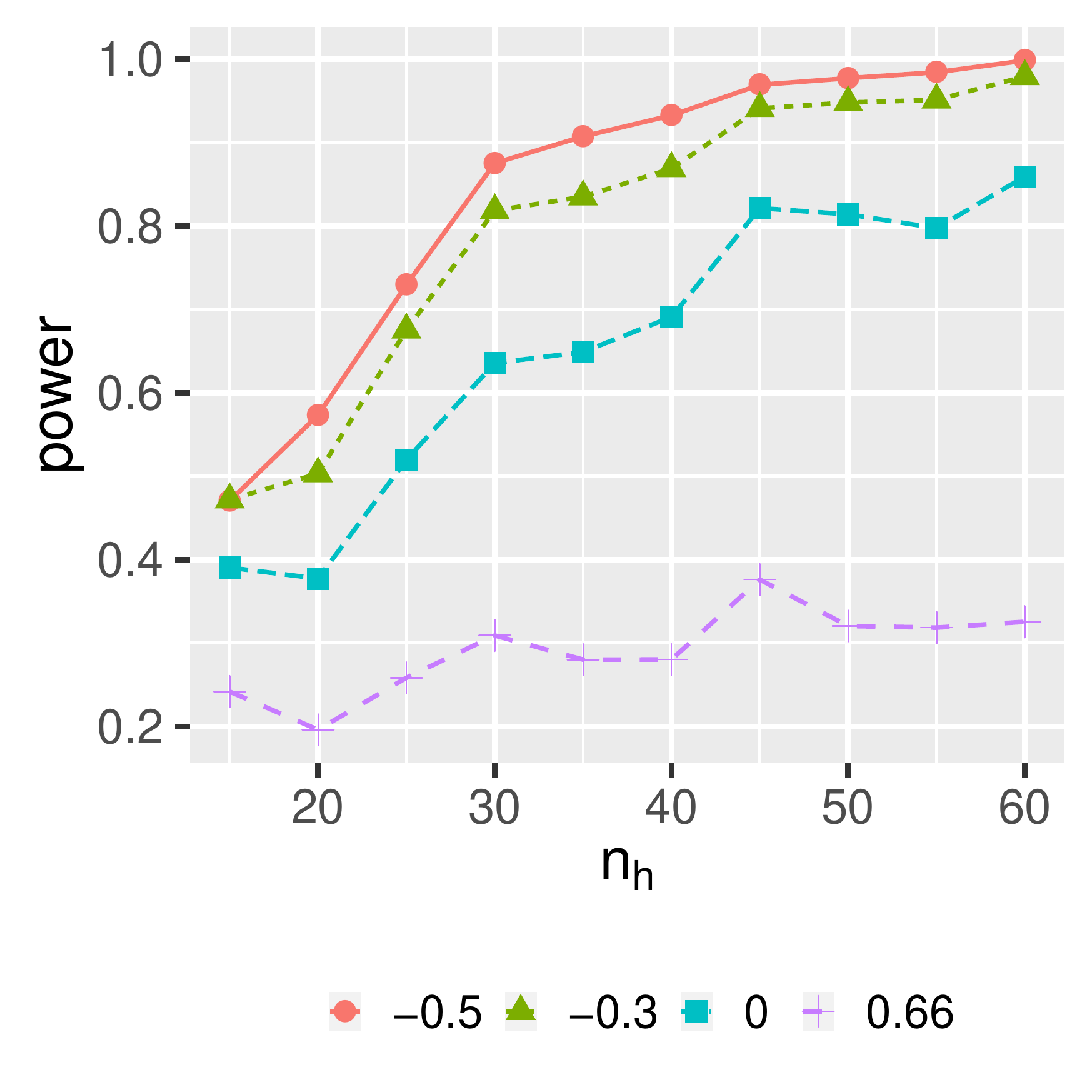}&
\includegraphics[scale=0.4]{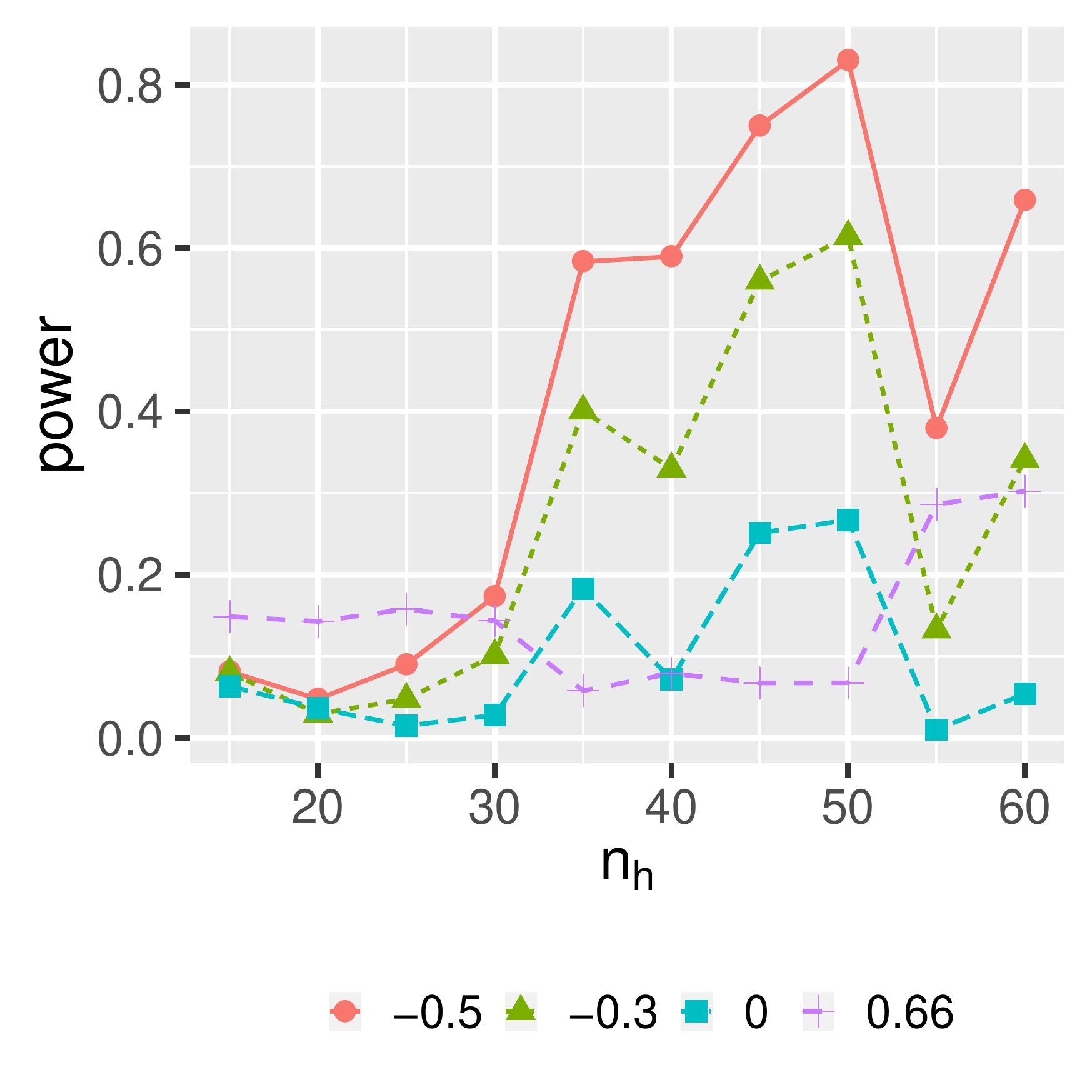}
\end{tabular}
\caption{RMSEs (top), emprirical levels (middle) and empirical powers (bottom).  Non-contaminated and contaminated settings (left and right, respectively). m-Inflated distribution. \label{fig:MC_mI}}
\end{figure}

In this scenario, the root of mean square error (RMSE) for the Cressie-Read PM$\phi$Es of $\boldsymbol{\beta}$ with $\lambda \in \{-0.5,-0.3,0,2/3 \}$ is studied, both for the contaminated and not-contaminated cases (see top of Figure \ref{fig:MC_mI}). To compute the accuracy in terms of contrast, we consider the testing problem

$$H_0: \beta_{11} = -0.9 \quad  \text{vs.} \quad H_1 : \beta_{11} \neq -0.9. $$ 
For computing the empirical test level, we measured the proportion of Wald-type test statistics ex-
ceeding the corresponding chi-square critical value. The simulated test powers were also
obtained under $H_1$ in  in a similar manner (here we consider $\beta_{11}=-1.5$).  We used a nominal level of $0.05$. Both levels and powers are presented in the middle and bottom of Figure \ref{fig:MC_mI}.

It is observed that PMLE ($\lambda=0$) presents the best behavior in terms of efficiency for the non-contaminated setting. In addition,  PM$\phi$E with $\lambda=0.66$ presents a RMSE lower than PM$\phi$E with negative values of $\lambda$.  On the other hand, when it is considered the contaminated setting, PM$\phi$E with $\lambda=-0.3$ presents a better behavior as it can be seen for high values of $n_h$. The other negative value considered for PM$\phi$E improves clearly the RMSE regarding to $\lambda=0.66$, again, for high values of $n_h$. However, the greatest difference is observed when studying empirical levels and powers. Although PMLE remains the best estimator for testing in a pure scenario, for the contaminated setting, better empirical levels are observed for negative values of the tuning parameter $\lambda$, in particular, for $\lambda=-0.5$. In terms of powers, negative values of $\lambda$ present better behavior in both settings, non-contaminated and contaminated. Positive values of $\lambda$ are presented as a good alternative only in terms of efficiency for small sample sizes, in concordance with \cite{Castilla2018_ASTA}. Other alternative scenarios are considered in Appendix \ref{sec:app_sim}.

\section{Numerical Examples \label{sec:num}}
\subsection{Education in Malawi (continuation)}
\begin{figure}[h!]
\centering
\begin{tabular}{rc}
\hspace{-0.4cm} \includegraphics[scale=0.3]{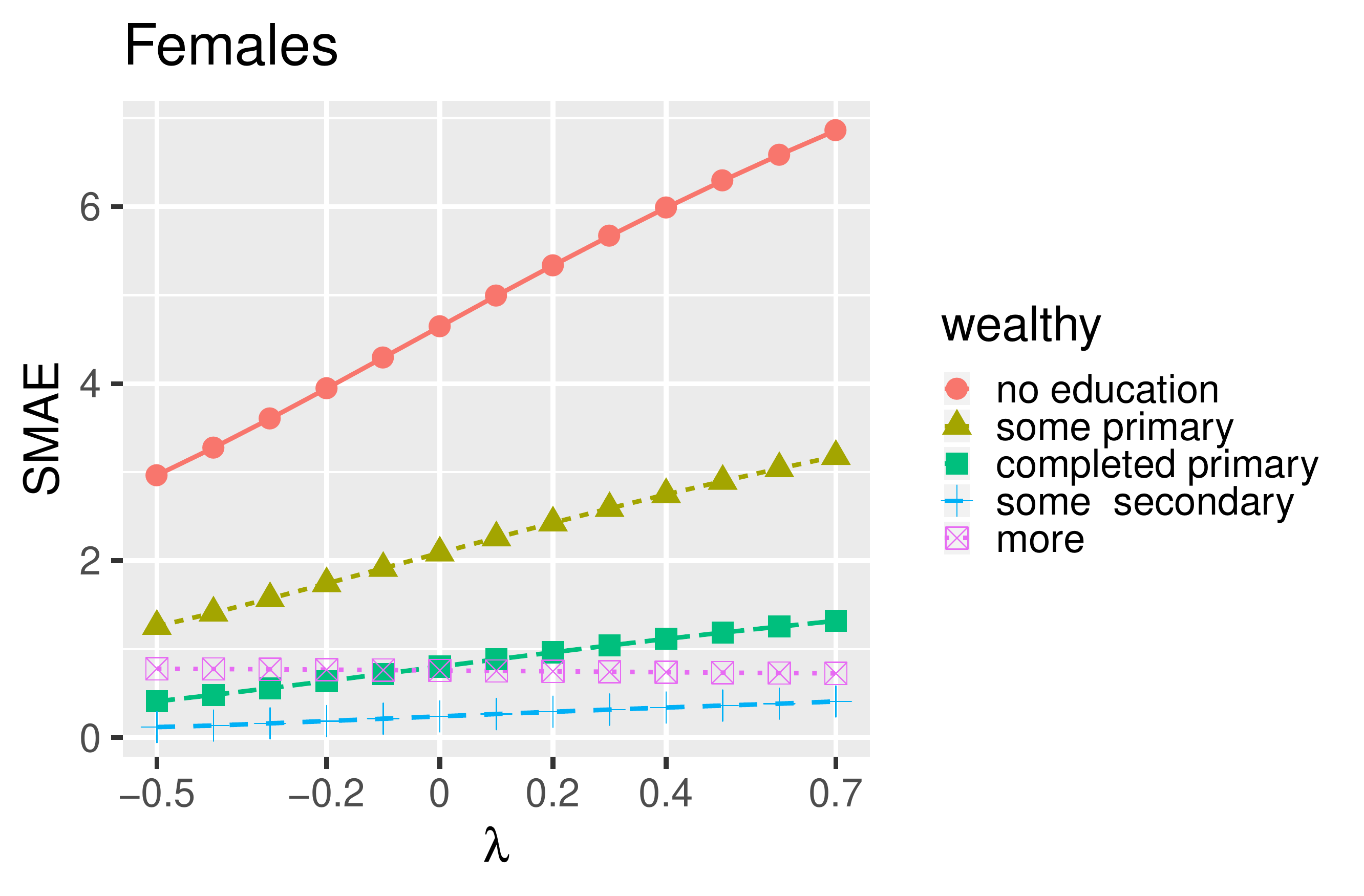}&
\includegraphics[scale=0.3]{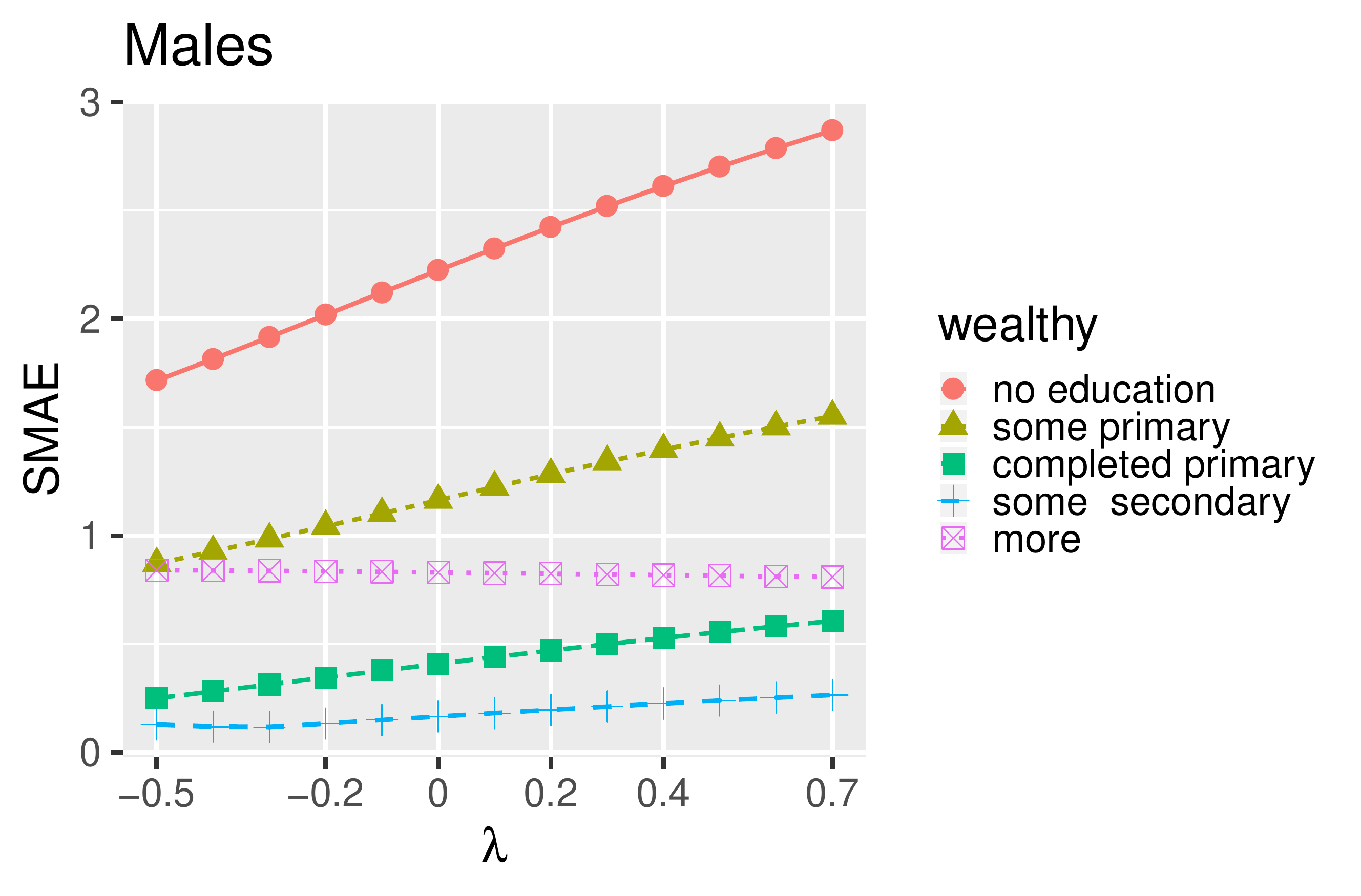}\\
\hspace{-0.4cm} \includegraphics[scale=0.3]{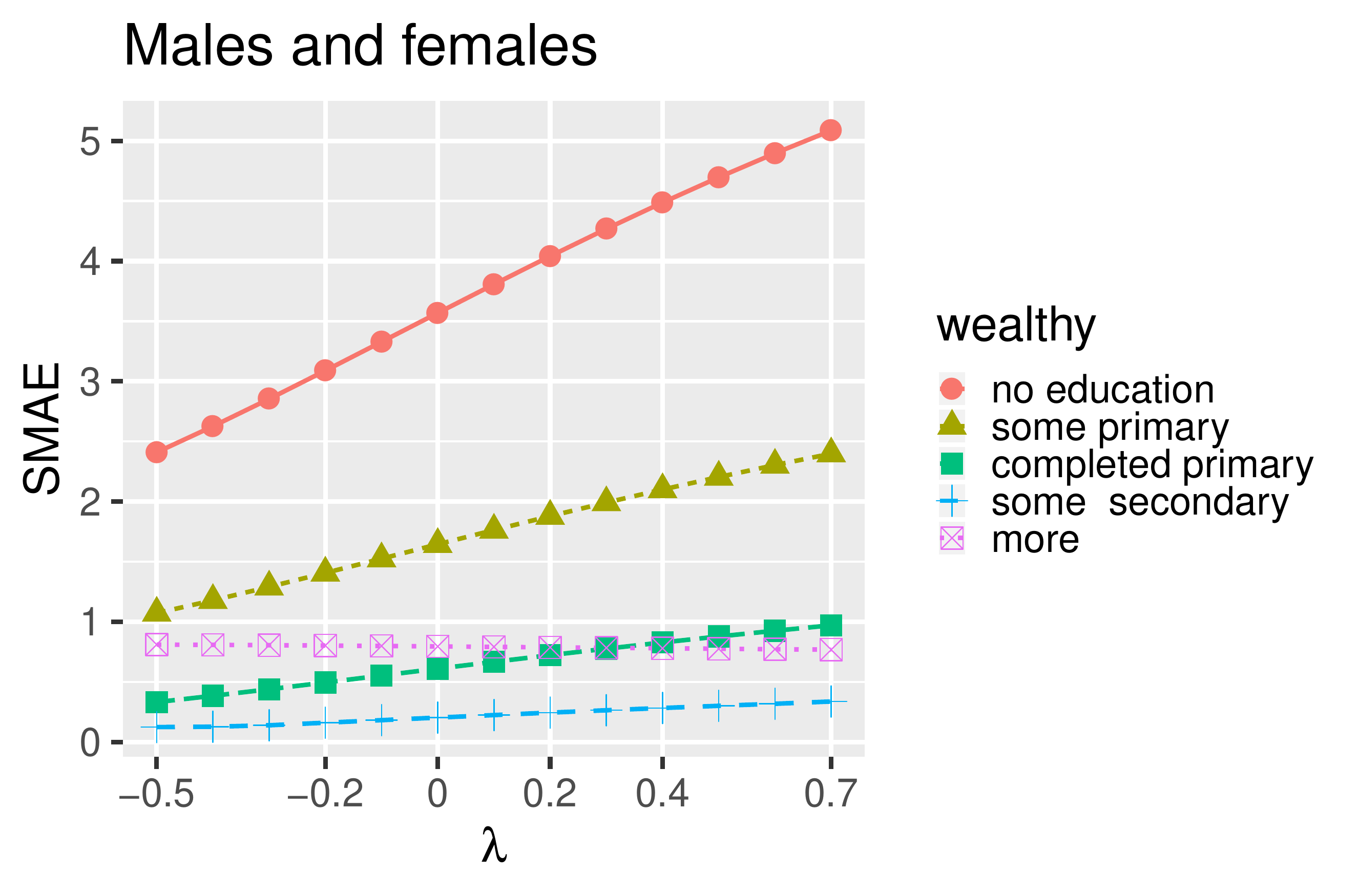}&
\hspace{-1.4cm} \includegraphics[scale=0.293]{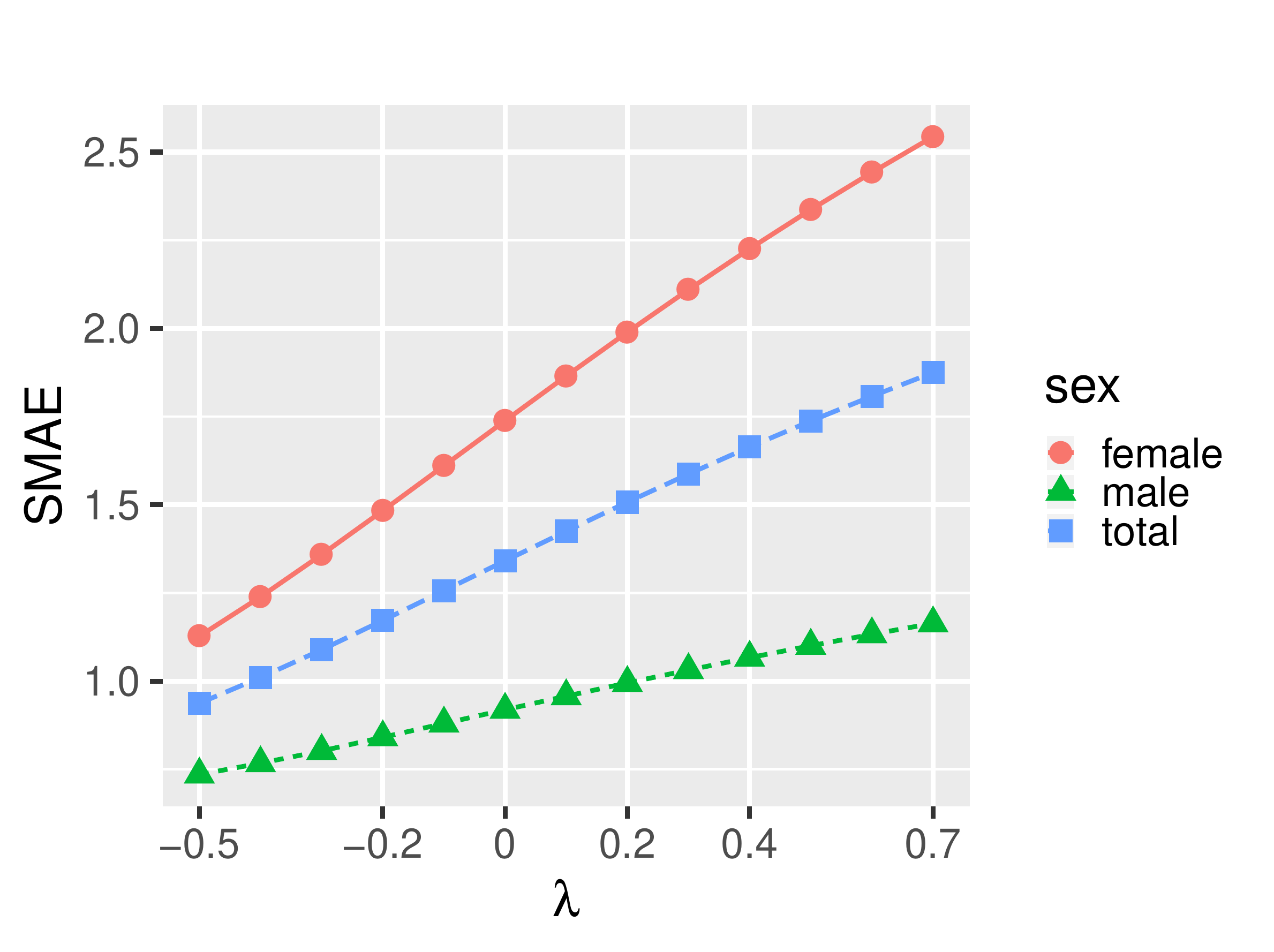}
\end{tabular}
\caption{Education in Malawi: estimated SMAEs for different values of the tuning parameter. \label{fig:malawi_SMAE}}
\end{figure}

Let us continue with the 2010 MDHS presented in Example \ref{example:malawi_mle}. As pointed out there, the 2010 MDHS presents data on educational attainment for female and male by its wealth quintile  level. In this section, we make a comparison of the behaviour between different PM$\phi$Es, when estimating the probabilities of the response categories. For this purpose, and after estimating the PM$\phi$Es in a grid of tuning parameters $\lambda \in \{-0.5, 0.7 \}$ in the Cressie-Read subfamily,  we measure a pondered standardized mean absolute error (SMAE) of the estimated probabilities against the observed probabilities. This is done by distinguishing the strata (wealth quintiles), the clusters (female and male) and jointly, as it can be seen in Figure \ref{fig:malawi_SMAE}.  The lowest SMAEs are obtained for negative values of $\lambda$. Then, they seem to offer a better behavior than the classical PMLE.\\


%

\begin{figure}[p]
\centering
\begin{tabular}{rc}
 \includegraphics[scale=0.4]{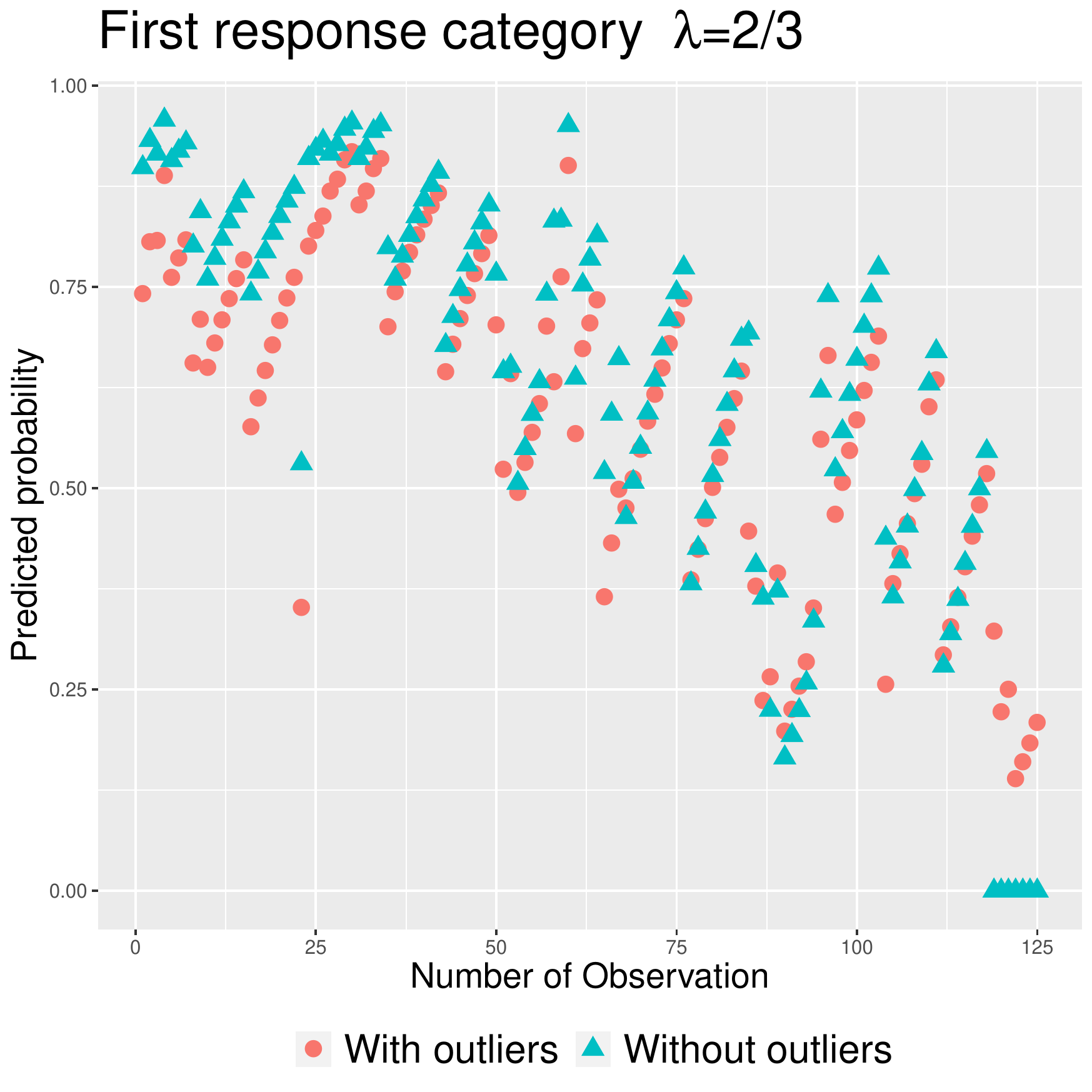}&
\includegraphics[scale=0.4]{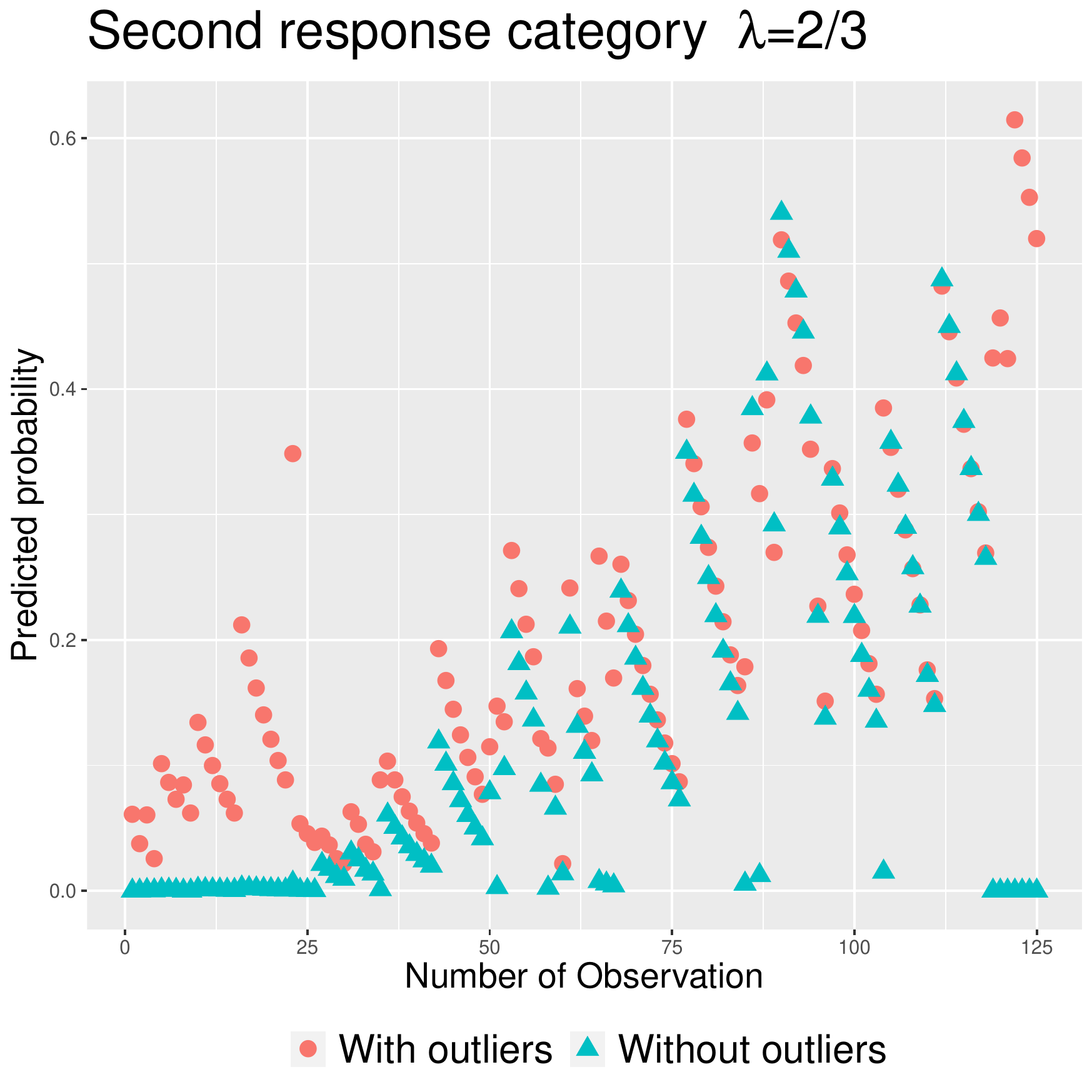}\\
 \includegraphics[scale=0.4]{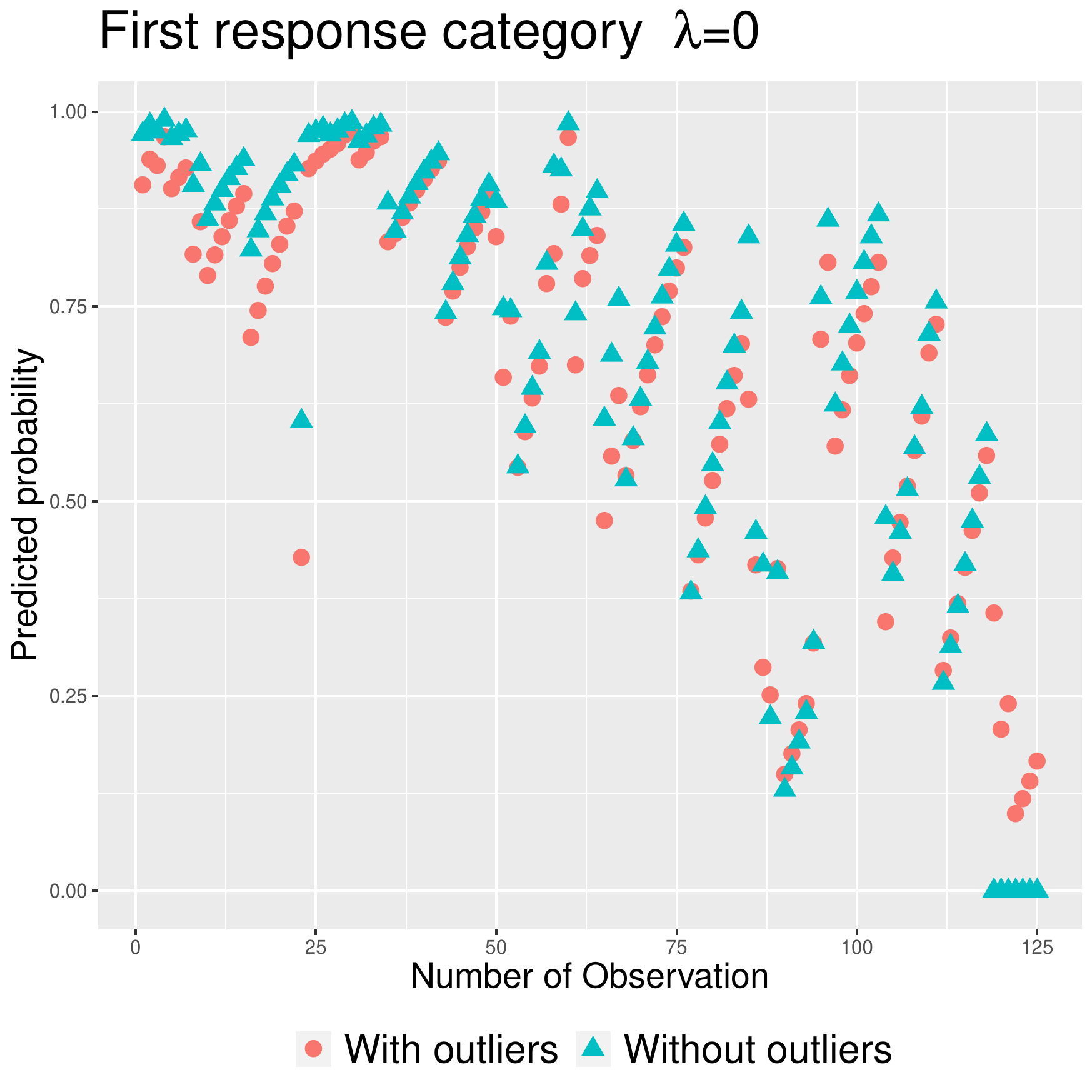}&
\includegraphics[scale=0.4]{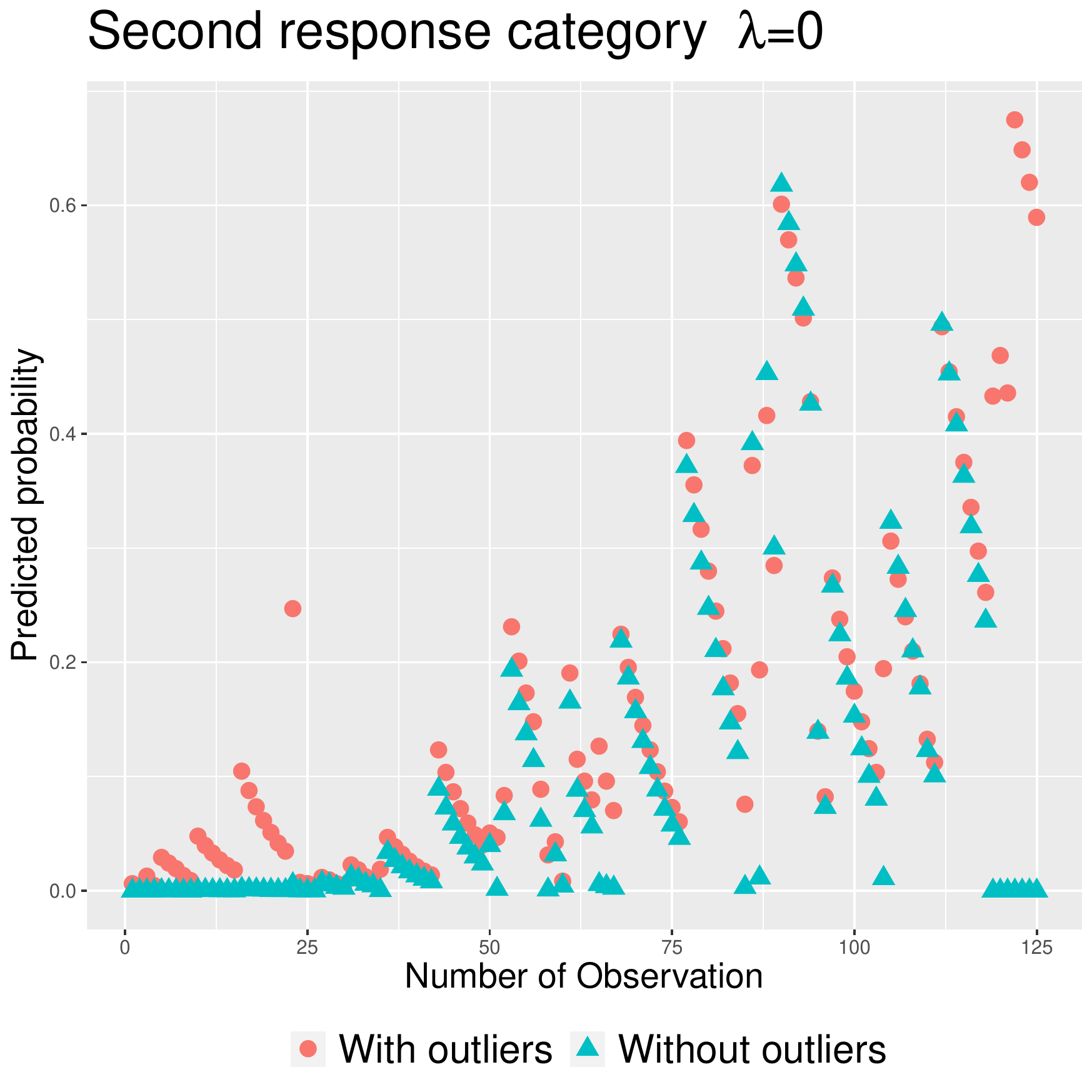}\\
 \includegraphics[scale=0.4]{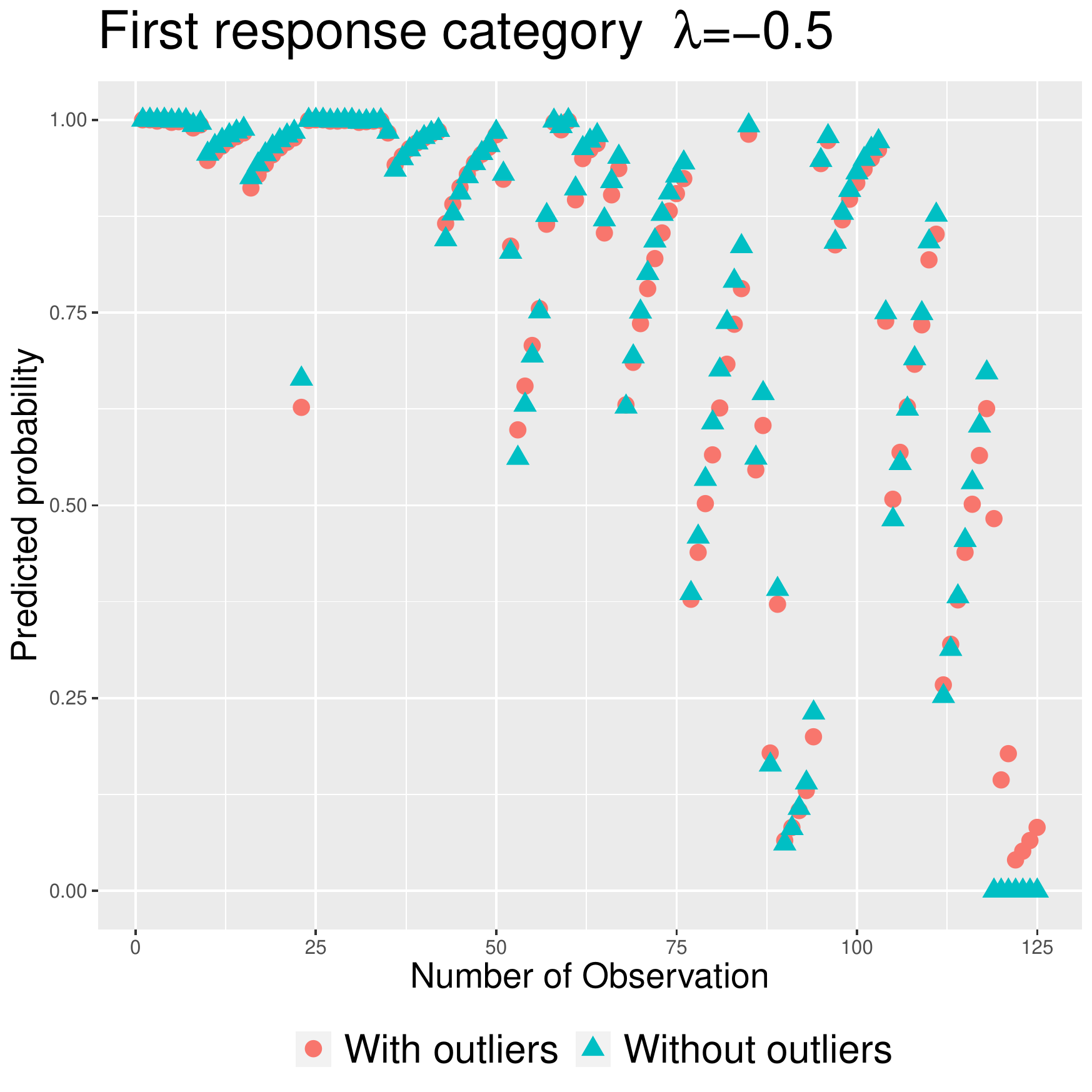}&
\includegraphics[scale=0.4]{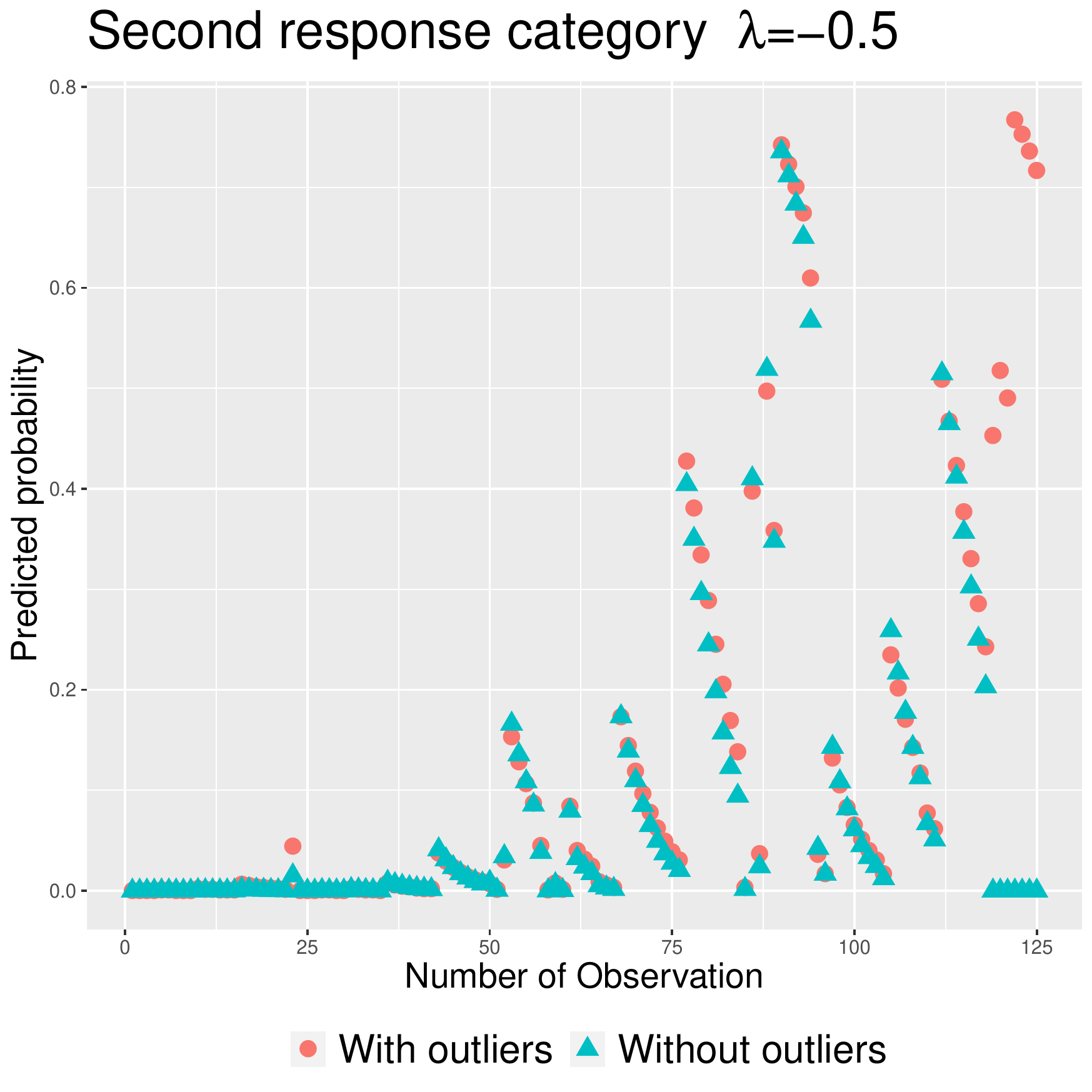}
\end{tabular}
\caption{Mamography example: predicted category probabilities of the response variable for the MLE ($\lambda=0$) and M$\phi$Es with $\lambda=2/3$ and $\lambda=-0.5$.  \label{fig:mamo_scatter}}
\end{figure}

\subsection{Mammography experience data \label{sec:mammo}}
As noted in Remark \ref{remark:simple_design}, multinomial logistic regression model under complex sample design is an extension of the classical one, evaluated under a  simple sample design. Therefore, the tools developed in this paper, can be also applied to these cases, in which the data design may be much simpler. In this section, we study the Mammography  experience data, a subset of a study by the University of Massachusetts Medical School,  introduced in \cite{hosmer2000}  and recently studied by \cite{martin2015} and \cite{castilla2018_BIOMETRICS}.   This study, which assess factors associated with women's knowledge, attitude and behavior towards mammography, involves $412$ individuals, grouped in $125$ distinct covariates values (which, somehow correspond to the ``clusters'' in a more complex survey) and $8$ explanatory variables, detailed in the cited bibliography. The response variable ME (Mammography experience) is a categorical factor with three levels: “Never”, “Within a Year” and “Over a Year”. As suggested by \cite{martin2015}, the groups of observations associated with covariate values $\boldsymbol{x}_i$ for $i\in \{1, 3, 17, 35,75,81,102\}$ can be treated as outliers. So this data set is a perfect candidate to show the robustness performance of the proposed estimators. \\

We compute the minimum $\phi$-divergence estimators (M$\phi$Es) of $\boldsymbol{\beta}$ for $\lambda \in \{-0.5,0,2/3 \}$, for the full dataset and also for the outliers deleted dataset. Moreover, we plot the corresponding (estimated) category probabilities for each available distinct covariate values. The left panel of Figure \ref{fig:mamo_scatter} presents these category probabilities for the first category, while the right panel presents these category probabilities for the second category.  Results clearly indicate the significant variation of the MLE and M$\phi$E with $\lambda=2/3$ in the presence  or absence of the outliers (red circles and blue triangles, respectively). However, the M$\phi$E with $\lambda=-0.5$ is shown to be much more stable, which is in concordance with its theoretical robustness. \\

\begin{figure}[h!!!!]
\centering
\includegraphics[scale=0.4]{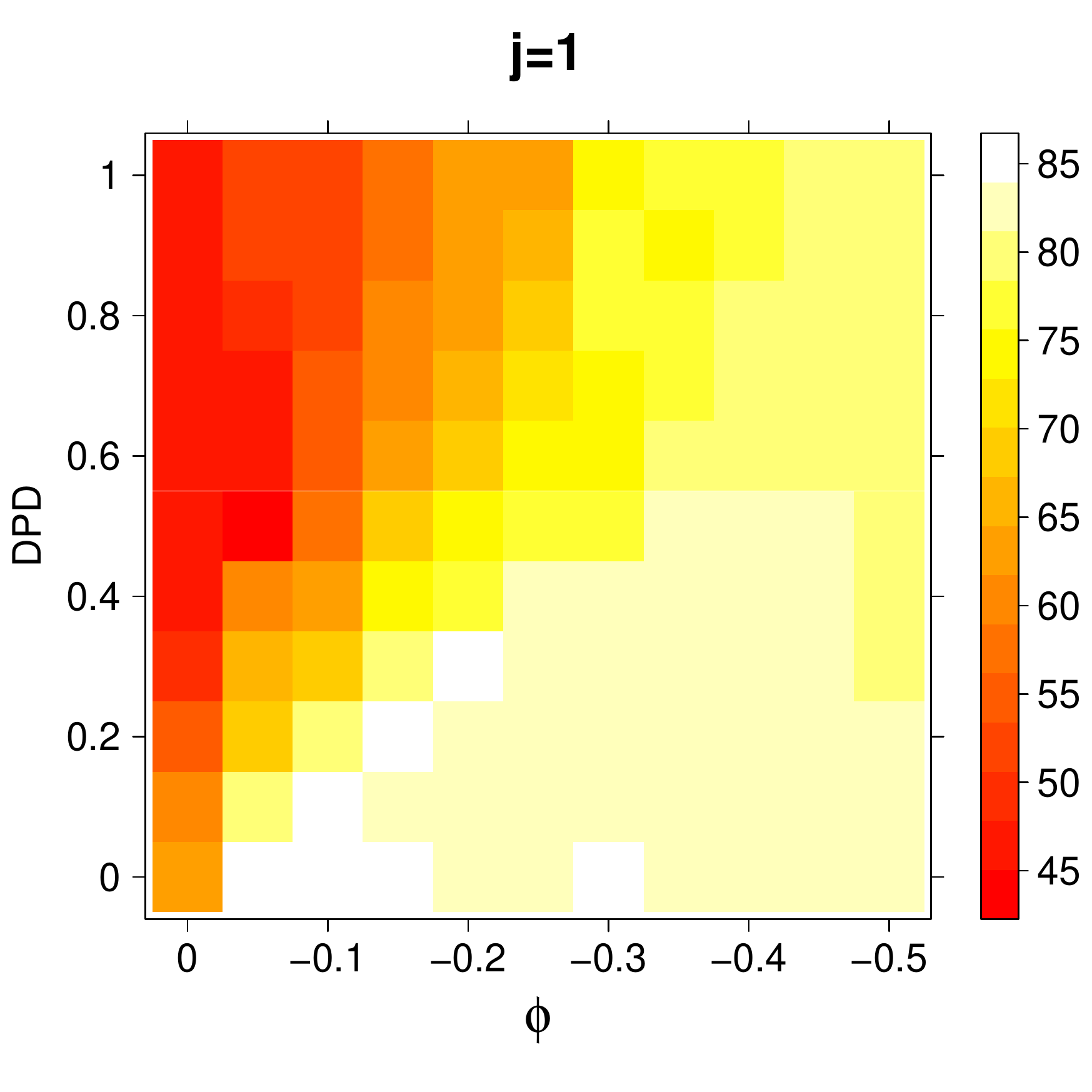} \
\includegraphics[scale=0.4]{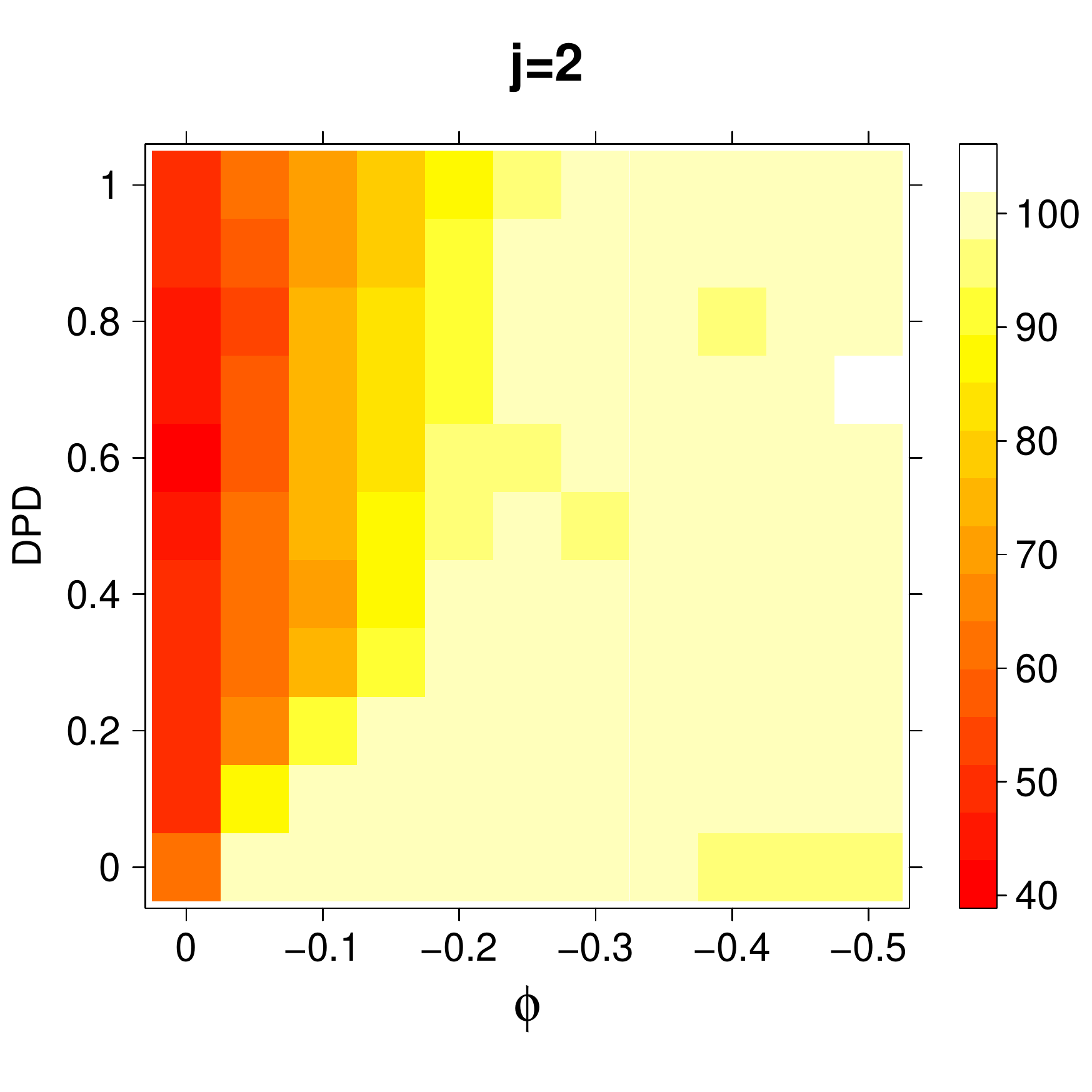} \
\caption{Mamography example: comparation of M$\phi$Es and MDPDEs. \label{fig:mamo_heats}}
\end{figure}

In \cite{castilla2018_BIOMETRICS} this dataset was also analyzed in order to illustrate the robustness of other family of estimators, those based on DPD divergences, the MDPDEs. This family is also parametrized by a tuning parameter, let say, $\lambda^* \geq 0$, and contains the MLE as a particular case for $\lambda^*=0$.   The question that could arise here is the point on using PM$\phi$Es instead of  MDPDEs.   In this regard, the efficiency of both family of estimators is compared  in the following way: for each pair of tuning parameters $(\lambda,\lambda^*)$ in a grid on $[-0.5,0]\times [0,1]$, the M$\phi$Es and MDPDEs are computed and the estimated probabilities for each of the categories of the response variable are obtained for each of the $I=125$ observations.  Then, we count the number of times, in these $125$ observations, that the M$\phi$E presents a lower error (the estimated probability is closer to the observed probability) than the MDPDE. The higher this value is, the better is the M$\phi$E with respect to the MDPDE. If the value is under $[I/2]= 63$, then the MDPDE is preferable. These results are illustrated in the two heat plots (for the first and second category, the third is omitted since it is similar) on Figure \ref{fig:mamo_heats}. We can observe how M$\phi$Es with a low value of $\lambda$ improves any MDPDE, while MDPDEs with a large value of $\lambda^*$ only improves PM$\phi$Es with tuning parameters near to $0$. The efficiency of the MLE is not comparable to any other option.
\bigskip

\begin{figure}[h!]
\centering
\begin{tabular}{rc}
 \includegraphics[scale=0.4]{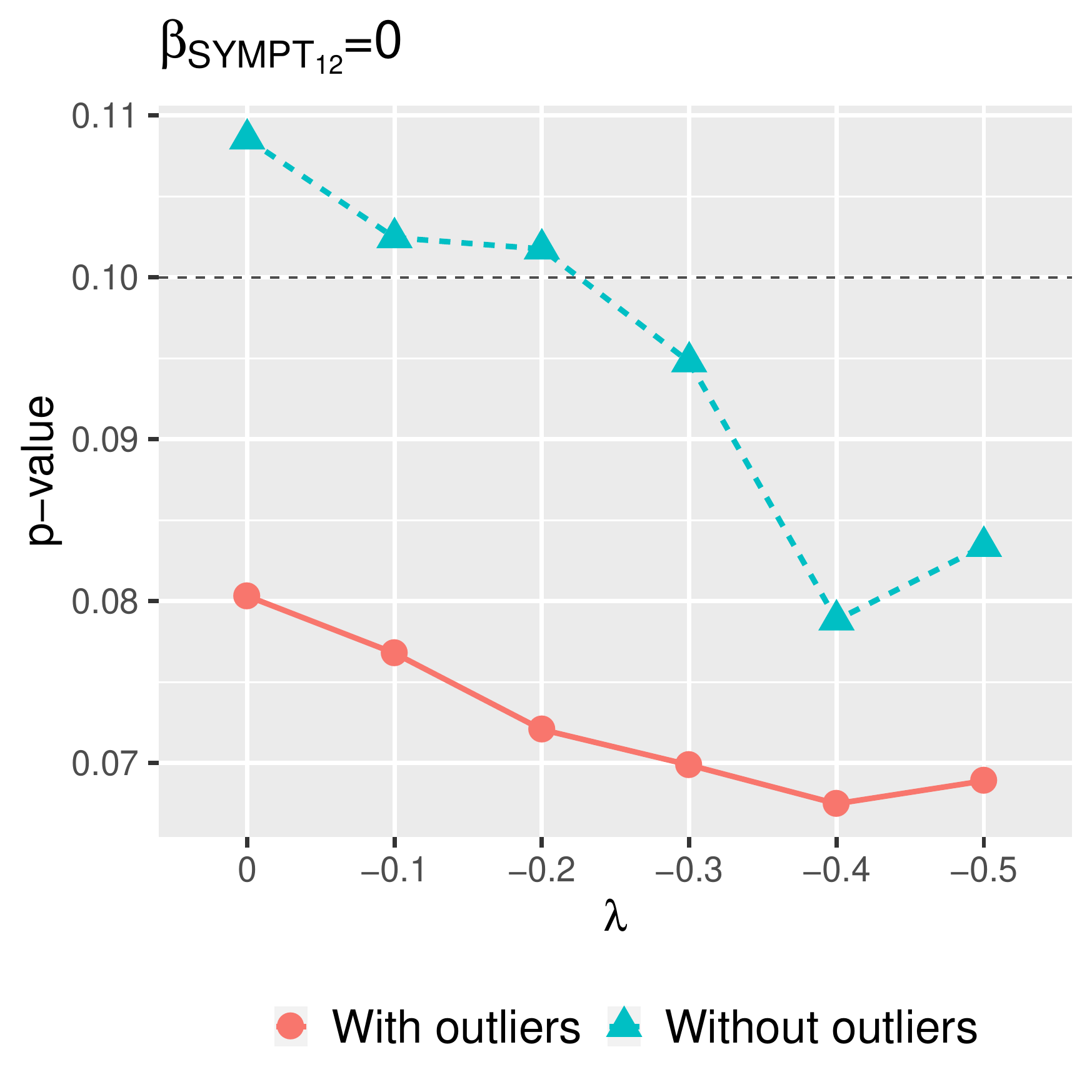}&
\includegraphics[scale=0.4]{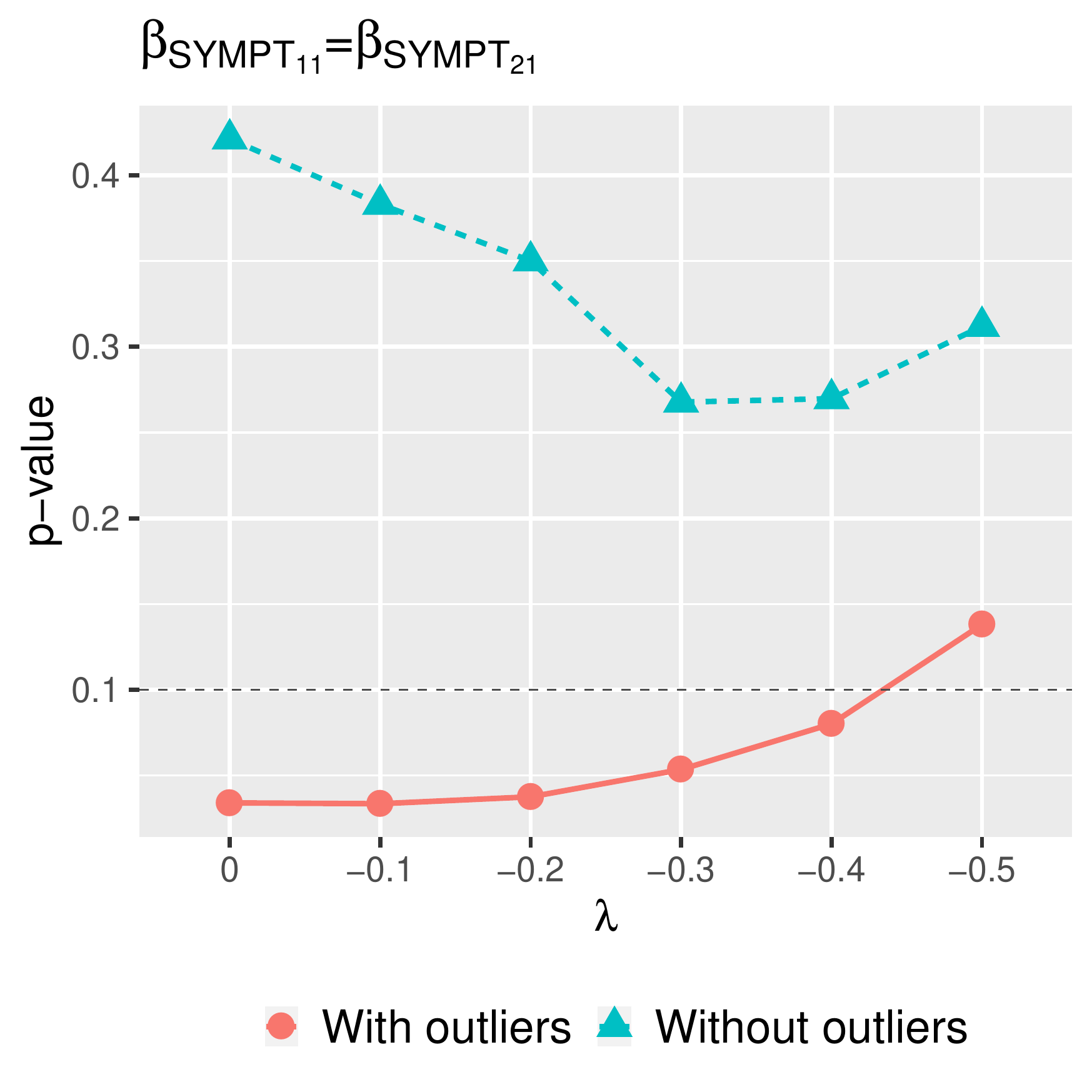}\\
\end{tabular}
\caption{Mamography example: p-values of the proposed M$\phi$Es based Wald-type tests. \label{fig:mamo_pvalues}}
\end{figure}

Now, we want to evaluate the robustness of the proposed Wald-type tests.  We consider the problem of testing 
\begin{align*}
&H_0: \beta_{SYMPT_{12}}=0,\\
&H_0: \beta_{SYMPT_{11}}=\beta_{SYMPT_{21}},
\end{align*}
for the variable SYMPT (``You do not need a mammogram unless you develop symptoms: 1, strongly agree; 2, agree; 3, disagree; 4, strongly disagree). The p-values obtained based on the proposed test are plotted over $\lambda$ in Figure \ref{fig:mamo_pvalues} for both the full and the outlier deleted data.  Clearly, the test decision at the significance level $\alpha=0.1$ changes completely in the presence of outliers for $\lambda$ near to $0$.\\

\section{Concluding Remarks and Future Work}
In this paper, we present robust estimators (PM$\phi$E) and Wald-type tests based on them for the multinomial logistic regression under complex survey, by means of $\phi$-divergence measures. In particular, we focus our study in the Cressie-Read subfamily of divergences, which are modelized by a tuning parameter $\lambda$. It is theoretically proved and empirically illustrated how PM$\phi$Es and Wald-type tests with $-1<\lambda<0$ are more robust than the classical  PMLE and Wald-test, presenting an interesting alternative to them. We believe that this method may be of special interest for analyzing demographic and health surveys, such as the one presented in Section \ref{sec:mammo}, as well as overall complex surveys for developing countries.

One of the problems that arises here is, given any data set, the choice of the tuning parameter $\lambda$. Robustness is usually accompanied with a loss of efficiency and other factors, such as sample size, can be also determinant in this decision. One possible way to make this choice is as follows: in a grid of possible tuning parameters, apply a measure of discrepancy to the data. Then, the tuning parameter that leads to the minimum discrepancy-statistic can be chosen as the “optimal” one. This is, somehow, the idea followed in the examples (see Figure \ref{fig:malawi_SMAE} and \ref{fig:mamo_heats}). Another alternative is the one proposed by  \cite{warwick2005_JSCS}, which consists on minimizing the estimated mean square error, computed as the sum of the squared (estimated) bias and variance. One of the main drawbacks of this method is the fact that it depends on a pilot estimator to estimate the bias. This problem was also highlighted recently in \cite{doi:10.1080/02664763.2020.1736524}, where an ``iterative Warwick and Jones algorithm'' (IJW algorithm) is proposed. Application of these methods and a development of new ones will be a challenging and interesting problem for further consideration.

\bigskip
\noindent \textbf{Acknowledgments:} This research is partially supported by Grant PGC2018-095194-B-I00, Grant FPU16/03104 and Grant BES-2016- 076669 from \textit{Ministerio de Ciencia, Innovación y Universidades}  and  \textit{Ministerio de Economía, Industria y Competitividad} (Spain). E. Castilla is member of the \textit{Instituto de Matemática Interdisciplinar}, Complutense University of Madrid.

\clearpage

\appendix
\section{Study of the Robustness of the proposed estimators and Wald-type tests \label{sec:modelrobust}}

An important concept in robustness theory is the influence function (\cite{hampel1986}). For any estimator defined in terms of a statistical functional $\boldsymbol{U}(F)$ from the true distribution $F$, its influence function (IF) is defined as 

\begin{equation}\label{eq:robust_1}
{IF}(t,\boldsymbol{U},F)=\lim_{\varepsilon \downarrow 0}\frac{\boldsymbol{U}(F_{\varepsilon })-\boldsymbol{U}(F)}{\varepsilon }=\left.  \frac{\partial \boldsymbol{U}(F_{\varepsilon })}{\partial \varepsilon }\right\vert _{\varepsilon =0^{+}},
\end{equation}
where $F_{\varepsilon }=(1-\varepsilon )F+\varepsilon \Delta _{t}$, with $\varepsilon $ being the contamination proportion and $\Delta _{t}$ being the degenerate distribution at the contamination point $t$. Thus, the (first-order) IF, as a function of $t$, measures the standardized asymptotic bias (in its first-order approximation) caused by the infinitesimal contamination at the point $t$. The maximum of this IF over $t$ indicates the extent of bias due to contamination and so smaller its value, the more robust the estimator is.

The classical theory deals with the case of independent and identically distributed (iid) observations (\cite{lindsay1994}). However, the present set-up of complex survey is not as simple as the iid set-up; in fact the observations within a cluster of a stratum are iid but the observations in different cluster and stratum are independent non-homogeneous. So we need to modify the definition of the influence function accordingly. Recently, \cite{ghosh2013,ghosh2015,ghosh2018} have discussed the extended definition of the influence function for the independent but non-homogeneous observations; however, all these works have been developed for  DPD-based estimators. Here, we will extend their definition for the case of multinomial logistic regression and PM$\phi$DEs.

We first need to define the statistical functional $\boldsymbol{U}_{\phi}(\boldsymbol{\beta})$ corresponding to the PM$\phi$Es as the minimizer of the phi-divergence between  the true and the model densities.  This is defined as the minimizer of

\begin{align*}
\boldsymbol{H}_{\phi}(\boldsymbol{\beta})=\frac{1}{\tau}\sum\limits_{h=1}^{H}\sum\limits_{i=1}^{n_{h}}w_{hi}m_{hi}\sum\limits_{s=1}^{d+1}\pi_{his}\left(\boldsymbol{\beta}\right) \phi\left( \frac{g_{his}}{\pi_{his}\left( \boldsymbol{\beta}\right) }\right). 
\end{align*}  
Under appropriate differentiability conditions as the solution of the estimating equations

\begin{align*}
\frac{\partial \boldsymbol{H}_{\phi}(\boldsymbol{\beta})}{\partial \boldsymbol{\beta}}=\sum\limits_{h=1}^{H}\sum\limits_{i=1}^{n_{h}}\frac{\omega_{hi}m_{hi}}{\phi''(1)}\frac{\partial \boldsymbol{\pi}^T_{hi}(\boldsymbol{\beta})}{\partial \boldsymbol{\beta}}\boldsymbol{f}_{\phi,hi}\left(\boldsymbol{g}_{hi},\boldsymbol{\beta} \right)=\boldsymbol{0}.
\end{align*}
In particular, for the Cressie-Read subfamily

\begin{align}\label{eq:robust_H}
\frac{\partial \boldsymbol{H}_{\phi_{\lambda}}(\boldsymbol{\beta})}{\partial \boldsymbol{\beta}}=\sum\limits_{h=1}^{H}\sum\limits_{i=1}^{n_{h}}\frac{\omega_{hi}m_{hi}}{\lambda+1}\frac{\partial \boldsymbol{\pi}^T_{hi}(\boldsymbol{\beta})}{\partial \boldsymbol{\beta}}\text{diag}^{-(\lambda+1)}(\boldsymbol{\pi}_{hi}(\boldsymbol{\beta}))\boldsymbol{g}_{hi}^{\lambda+1}=\boldsymbol{0}.
\end{align}
 For simplicity, let us first assume that the contamination is only in one
cluster probability $g_{h_0i_0}$ for some fix $h_0$ and $i_0$.  Consider the contaminated probability vector

\begin{align*}
\boldsymbol{g}_{h_0i_0,\varepsilon}=(1-\varepsilon)\boldsymbol{\pi}_{h_0i_0}(\boldsymbol{\beta}^0)+\varepsilon \delta_{\boldsymbol{t}},
\end{align*}
where $\varepsilon$ is the contamination proportion and $\delta_{\boldsymbol{t}}$ is the degenerate probability at the outlier point $\boldsymbol{t}=(t_1,\dots,t_{d+1})^T\in \{0,1 \}^{d+1}$ with $\sum_{s=1}^{d+1}t_s=1$ and 

\begin{align*}
\boldsymbol{g}_{hi}= \left\{ \begin{array}{ll}
         \boldsymbol{\pi}_{hi}(\boldsymbol{\beta}^0) & \mbox{if $(i,h) \neq (i_0,h_0)$};\\
        \boldsymbol{g}_{h_0i_0,\varepsilon} & \mbox{if $(i,h) =(i_0,h_0)$}.\end{array} \right.
\end{align*}

 Denote the corresponding contaminated full probability vector as $\boldsymbol{g}_{\varepsilon}$ which is the same as $\boldsymbol{g}$ except $g_{h_0i_0}$ being replaced by $g_{h_0i_0,\varepsilon}$ and let the corresponding contaminated distribution vector be $\boldsymbol{G}_{\varepsilon}$. We replace $\boldsymbol{\beta}$ in (\ref{eq:robust_H}) by \ $\boldsymbol{\beta}_{\varepsilon}^{h_0i_0}=\boldsymbol{U}_{\phi_{\lambda}}(\boldsymbol{G}_{\varepsilon})$. Then, we have

\begin{align}\label{eq:robust_H_mod}
&\left.\frac{\partial \boldsymbol{H}_{\phi_{\lambda}}(\boldsymbol{\beta})}{\partial \boldsymbol{\beta}}\right|_{\boldsymbol{\beta}=\boldsymbol{\beta}_{\varepsilon}^{h_0i_0}}  \\
&=\sum_{\substack{h=1 \\ (i,h)\neq (i_0,h_0)}}^{H}\sum\limits_{i=1}^{n_{h}}\left\{\frac{\omega_{hi}m_{hi}}{\lambda+1}\left.\frac{\partial \boldsymbol{\pi}^T_{hi}(\boldsymbol{\beta})}{\partial \boldsymbol{\beta}}\right|_{\boldsymbol{\beta}=\boldsymbol{\beta}_{\varepsilon}^{h_0i_0}}\text{diag}^{-(\lambda+1)}(\boldsymbol{\pi}_{hi}(\boldsymbol{\beta}_{\varepsilon}^{h_0i_0}))\boldsymbol{\pi}_{hi}^{\lambda+1}(\boldsymbol{\beta}^0) \right\}\notag \\
& \quad + \frac{\omega_{h_0i_0}m_{h_0i_0}}{\lambda+1}\left.\frac{\partial \boldsymbol{\pi}^T_{hi}(\boldsymbol{\beta})}{\partial \boldsymbol{\beta}}\right|_{\boldsymbol{\beta}=\boldsymbol{\beta}_{\varepsilon}^{h_0i_0}}\text{diag}^{-(\lambda+1)}(\boldsymbol{\pi}_{h_0i_0}(\boldsymbol{\beta}_{\varepsilon}^{h_0i_0}))\left[ (1-\varepsilon)\boldsymbol{\pi}_{h_0i_0}^{\lambda+1}(\boldsymbol{\beta}^0)+\varepsilon \delta_{\boldsymbol{t}}^{\lambda+1}\right]. \notag
\end{align}
Now, we are going to get the derivative of (\ref{eq:robust_H_mod}) with respect to $\varepsilon$.

\begin{align*}
&-(\lambda+1)\sum_{\substack{h=1 \\ (i,h)\neq (i_0,h_0)}}^{H}\sum\limits_{i=1}^{n_{h}}\left\{\frac{\omega_{hi}m_{hi}}{\lambda+1}\left.\frac{\partial \boldsymbol{\pi}^T_{hi}(\boldsymbol{\beta})}{\partial \boldsymbol{\beta}}\right|_{\boldsymbol{\beta}=\boldsymbol{\beta}_{\varepsilon}^{h_0i_0}}\text{diag}^{-(\lambda+2)}(\boldsymbol{\pi}_{hi}(\boldsymbol{\beta}_{\varepsilon}^{h_0i_0})) \right.\\
& \left. \quad \quad  \times \left.\frac{\partial \boldsymbol{\pi}^T_{hi}(\boldsymbol{\beta})}{\partial \boldsymbol{\beta}}\right|_{\boldsymbol{\beta}=\boldsymbol{\beta}_{\varepsilon}^{h_0i_0}} \boldsymbol{\pi}_{hi}^{\lambda+1}(\boldsymbol{\beta}^0)\frac{\partial \boldsymbol{\beta}_{\varepsilon}^{h_0i_0}}{\partial \varepsilon} \right\}\notag \\
&+\sum_{\substack{h=1 \\ (i,h)\neq (i_0,h_0)}}^{H}\sum\limits_{i=1}^{n_{h}}\left\{\frac{\omega_{hi}m_{hi}}{\lambda+1}\left.\frac{\partial^2 \boldsymbol{\pi}^T_{hi}(\boldsymbol{\beta})}{\partial \boldsymbol{\beta}\partial \boldsymbol{\beta}^T}\right|_{\boldsymbol{\beta}=\boldsymbol{\beta}_{\varepsilon}^{h_0i_0}}\text{diag}^{-(\lambda+1)}(\boldsymbol{\pi}_{hi}(\boldsymbol{\beta}_{\varepsilon}^{h_0i_0}))\boldsymbol{\pi}_{hi}^{\lambda+1}(\boldsymbol{\beta}^0)\frac{\partial \boldsymbol{\beta}_{\varepsilon}^{h_0i_0}}{\partial \varepsilon} \right\}\\
&  + \left\{ \frac{\omega_{h_0i_0}m_{h_0i_0}}{\lambda+1}\left.\frac{\partial \boldsymbol{\pi}^T_{hi}(\boldsymbol{\beta})}{\partial \boldsymbol{\beta}}\right|_{\boldsymbol{\beta}=\boldsymbol{\beta}_{\varepsilon}^{h_0i_0}}\text{diag}^{-(\lambda+2)}(\boldsymbol{\pi}_{h_0i_0}(\boldsymbol{\beta}_{\varepsilon}^{h_0i_0}))\left.\frac{\partial \boldsymbol{\pi}^T_{hi}(\boldsymbol{\beta})}{\partial \boldsymbol{\beta}}\right|_{\boldsymbol{\beta}=\boldsymbol{\beta}_{\varepsilon}^{h_0i_0}} \right.\\
& \left. \quad \quad  \times \left[ (1-\varepsilon)\boldsymbol{\pi}^{\lambda+1}_{h_0i_0}(\boldsymbol{\beta}^0)+\varepsilon \delta^{\lambda+1}_{\boldsymbol{t}}\right]\frac{\partial \boldsymbol{\beta}_{\varepsilon}^{h_0i_0}}{\partial \varepsilon} \right\} \notag\\
& + \left\{ \frac{\omega_{h_0i_0}m_{h_0i_0}}{\lambda+1}\left.\frac{\partial^2 \boldsymbol{\pi}^T_{hi}(\boldsymbol{\beta})}{\partial \boldsymbol{\beta}\partial \boldsymbol{\beta}^T}\right|_{\boldsymbol{\beta}=\boldsymbol{\beta}_{\varepsilon}^{h_0i_0}}\text{diag}^{-(\lambda+1)}(\boldsymbol{\pi}_{h_0i_0}(\boldsymbol{\beta}_{\varepsilon}^{h_0i_0})) \right. \\
& \left. \quad \quad \times \left[ (1-\varepsilon)\boldsymbol{\pi}^{\lambda+1}_{h_0i_0}(\boldsymbol{\beta}^0)+\varepsilon \delta^{\lambda+1}_{\boldsymbol{t}}\right] \frac{\partial \boldsymbol{\beta}_{\varepsilon}^{h_0i_0}}{\partial \varepsilon} \right\} \\
&  + \left\{ \omega_{h_0i_0}m_{h_0i_0}\left.\frac{\partial \boldsymbol{\pi}^T_{hi}(\boldsymbol{\beta})}{\partial \boldsymbol{\beta}}\right|_{\boldsymbol{\beta}=\boldsymbol{\beta}_{\varepsilon}^{h_0i_0}}\text{diag}^{-(\lambda+1)}(\boldsymbol{\pi}_{h_0i_0}(\boldsymbol{\beta}_{\varepsilon}^{h_0i_0}))  \left[ -\boldsymbol{\pi}^{\lambda+1}_{h_0i_0}(\boldsymbol{\beta}^0)+ \delta^{\lambda+1}_{\boldsymbol{t}}\right] \right\}=\boldsymbol{0}. 
\end{align*}
Now, evaluating the previous expression in $\varepsilon=0$ and simplifying, we have

\begin{align*}
&IF(t_{h_0i_0},\boldsymbol{U}_{\phi_{\lambda}}(\boldsymbol{\beta}),F_{\boldsymbol{\beta}^0})\sum_{\substack{h=1}}^{H}\sum\limits_{i=1}^{n_{h}}\left\{\frac{\omega_{hi}m_{hi}}{\lambda+1}\left.\frac{\partial^2 \boldsymbol{\pi}^T_{hi}(\boldsymbol{\beta})}{\partial \boldsymbol{\beta}\partial \boldsymbol{\beta}^T}\right|_{\boldsymbol{\beta}=\boldsymbol{\beta}^{0}}\text{diag}^{-(\lambda+1)}(\boldsymbol{\pi}_{hi}(\boldsymbol{\beta}^{0}))\boldsymbol{\pi}_{hi}^{\lambda+1}(\boldsymbol{\beta}^0) \right.\\
& \left. \quad \quad \quad  \quad  - \omega_{hi}m_{hi}\left.\frac{\partial \boldsymbol{\pi}^T_{hi}(\boldsymbol{\beta})}{\partial \boldsymbol{\beta}}\right|_{\boldsymbol{\beta}=\boldsymbol{\beta}^{0}}\text{diag}^{-(\lambda+2)}(\boldsymbol{\pi}_{hi}(\boldsymbol{\beta}^{0})) \left.\frac{\partial \boldsymbol{\pi}^T_{hi}(\boldsymbol{\beta})}{\partial \boldsymbol{\beta}}\right|_{\boldsymbol{\beta}=\boldsymbol{\beta}^{0}} \boldsymbol{\pi}_{hi}^{\lambda+1}(\boldsymbol{\beta}^0)\right\}\\
& + \left\{\omega_{h_0i_0}m_{h_0i_0}\left.\frac{\partial \boldsymbol{\pi}^T_{hi}(\boldsymbol{\beta})}{\partial \boldsymbol{\beta}}\right|_{\boldsymbol{\beta}=\boldsymbol{\beta}^{0}}\text{diag}^{-(\lambda+1)}(\boldsymbol{\pi}_{h_0i_0}(\boldsymbol{\beta}^{0}))\left[ -\boldsymbol{\pi}^{\lambda+1}_{h_0i_0}(\boldsymbol{\beta}^0)+ \delta^{\lambda+1}_{\boldsymbol{t}}\right]\right\}=\boldsymbol{0}.
\end{align*}
Proposition \ref{th: IF} follows straightforward.

\begin{theorem}\label{th: IF} Let us consider the multinomial logistic regression model under complex design given in (\ref{eq:pi}). The IF of the PM$\phi$Es with respect to the $i_0$ cluster in the $h_0$ stratum is given by

\begin{align}\label{eq:IF}
&IF(t_{h_0i_0},\boldsymbol{U}_{\phi_{\lambda}}(\boldsymbol{\beta}),F_{\boldsymbol{\beta}^0})= \boldsymbol{\Psi}_{n,\lambda}^{-1}(\boldsymbol{\beta}) \boldsymbol{u}^*_{\phi_{\lambda},h_0i_0}(\boldsymbol{\beta}),
\end{align}
where

\begin{align}\label{eq:IF_u}
\boldsymbol{u}^*_{\phi_{\lambda},h_0i_0}(\boldsymbol{\beta})=\left[\omega_{h_0i_0}m_{h_0i_0}\left.\frac{\partial \boldsymbol{\pi}^T_{hi}(\boldsymbol{\beta})}{\partial \boldsymbol{\beta}}\right|_{\boldsymbol{\beta}=\boldsymbol{\beta}^{0}}\text{diag}^{-(\lambda+1)}(\boldsymbol{\pi}_{h_0i_0}(\boldsymbol{\beta}^{0}))\delta^{\lambda+1}_{\boldsymbol{t}}\right]
\end{align}
and

\begin{align*}
\boldsymbol{\Psi}_{n,\phi_{\lambda}}(\boldsymbol{\beta})=&-\sum_{\substack{h=1}}^{H}\sum\limits_{i=1}^{n_{h}}\left\{\frac{\omega_{hi}m_{hi}}{\lambda+1}\left.\frac{\partial^2 \boldsymbol{\pi}^T_{hi}(\boldsymbol{\beta})}{\partial \boldsymbol{\beta}\partial \boldsymbol{\beta}^T}\right|_{\boldsymbol{\beta}=\boldsymbol{\beta}^{0}}\boldsymbol{1}_{d+1}^T \right\}\\
& +\sum_{\substack{h=1}}^{H}\sum\limits_{i=1}^{n_{h}}\left\{\omega_{hi}m_{hi}\left.\frac{\partial \boldsymbol{\pi}^T_{hi}(\boldsymbol{\beta})}{\partial \boldsymbol{\beta}}\right|_{\boldsymbol{\beta}=\boldsymbol{\beta}^{0}}\text{diag}^{-(\lambda+2)}(\boldsymbol{\pi}_{hi}(\boldsymbol{\beta}^{0})) \left.\frac{\partial \boldsymbol{\pi}^T_{hi}(\boldsymbol{\beta})}{\partial \boldsymbol{\beta}}\right|_{\boldsymbol{\beta}=\boldsymbol{\beta}^{0}} \boldsymbol{\pi}_{hi}^{\lambda+1}(\boldsymbol{\beta}^0)\right\}.
\end{align*}
\end{theorem}

\smallskip

\begin{remark}
Let us consider the right part of the IF (\ref{eq:IF}). Equation (\ref{eq:IF_u}) can be expressed as

\begin{align*}
\boldsymbol{u}^*_{\phi_{\lambda},h_0i_0}(\boldsymbol{\beta})
&=\omega_{h_0i_0}m_{h_0i_0}\boldsymbol{\Delta}^*(\boldsymbol{\pi}_{h_0i_0}(\boldsymbol{\beta}^0))\text{diag}^{-(\lambda+1)}(\boldsymbol{\pi}_{h_0i_0}(\boldsymbol{\beta}^{0}))\delta^{\lambda+1}_{\boldsymbol{t}}\otimes \boldsymbol{x}_{h_0i_0}
\end{align*}
and, in the particular case $\lambda=0$ (PMLE)

\begin{align}\label{eq:IF_PMLE}
\boldsymbol{u}^*_{h_0i_0}(\boldsymbol{\beta})
&=\omega_{h_0i_0}m_{h_0i_0}(\boldsymbol{\pi}^*_{h_0i_0}(\boldsymbol{\beta}^0)-\delta^*_{\boldsymbol{t}})\otimes \boldsymbol{x}_{h_0i_0}.
\end{align}
While the influence of vertical outliers on the PM$\phi$Es is bounded as $\boldsymbol{t}$ changes its indicative category only; by the assumed form of the model probability, the IF is bounded for all $-1< \lambda<0$ but unbounded at $\lambda=0$ against ``bad'' leverage points. Effectively,  when $\boldsymbol{x}_{h_0i_0}$ increases in (\ref{eq:IF_PMLE}), the residual $(\boldsymbol{\pi}^*_{h_0i_0}(\boldsymbol{\beta}^0)-\delta^*_{\boldsymbol{t}})$ will typically tend much faster to zero than $\boldsymbol{x}_{h_0i_0}$ to infinity, resulting in a small influence. But a ``bad'' leverage point associated to a misclassified  observation will lead to an infinite value of the IF (see \cite{croux2003a}).
\end{remark}

Similarly, one can show that, in the case there is contamination in some of the clusters within some stratum, the boundedness and robustness implications for the IF are exactly the same as before.

\bigskip

Let us now study the robustness of proposed Wald-type tests through the IF of the corresponding Wald-type test statistics defined in Section \ref{Sec:Wald}. In our context, the functional associated with the Wald-type test, evaluated at $\boldsymbol{U}_{\beta}(\boldsymbol{G})$ is given by

\begin{equation*}
W_{n}(\boldsymbol{U}_{\beta}(\boldsymbol{G}))=n\left(\boldsymbol{M}^T \boldsymbol{U}_{\beta}(\boldsymbol{G})-\boldsymbol{m}\right)^T\left( \boldsymbol{M}^{T}\boldsymbol{V}_n(\boldsymbol{U}_{\beta}(\boldsymbol{G}))\boldsymbol{M}\right) ^{-1}\left(\boldsymbol{M}^T \boldsymbol{U}_{\beta}(\boldsymbol{G})-\boldsymbol{m}\right).
\label{eq:IF_Wald}
\end{equation*}%

The IF of general Wald-type tests under such non-homogeneous set-up has been extensively studied in \cite{basu2018_METRIKA}, for the case of DPD estimators. Here, we can follow that the first-order IF of $W_{n}$, defined as the first order derivative of
its value at the contaminated distribution with respect to $\varepsilon$ at $\varepsilon=0$, become null at the null distribution.  Therefore, the first order IF is not informative in this case of Wald-type tests, and we need to investigate the second-order IF, let say $IF^{(2)}$, of $W_{n}$.  Through some computations we obtain the following proposition.

\begin{theorem}\label{th: IF} Let us consider the multinomial logistic regression model under complex design given in (\ref{eq:pi}). The second-order IF of the functional associated with the Wald-type tests with respect to the $i_0$ cluster in the $h_0$ stratum is given by

\begin{align}\label{eq:IF2_Wald}
&IF^{(2)}(t_{h_0i_0},W_n,F_{\boldsymbol{\beta}^0}) \notag \\
&= 2 IF^{T}(t_{h_0i_0},W_n,F_{\boldsymbol{\beta}^0}) \boldsymbol{M}\left( \boldsymbol{M}^{T}\boldsymbol{V}_n(\boldsymbol{U}_{\beta}(\boldsymbol{G}))\boldsymbol{M}\right) ^{-1} \boldsymbol{M}^TIF(t_{h_0i_0},W_n,F_{\boldsymbol{\beta}^0}), 
\end{align}
where $IF(t_{h_0i_0},W_n,F_{\boldsymbol{\beta}^0})$ is the IF  of the PM$\phi$Es given in (\ref{eq:IF}).
\end{theorem}

\begin{remark}
Note that, the second-order IF of the proposed Wald-type tests is a quadratic function of the corresponding IF of the PM$\phi$Es. Therefore, the boundedness of the influence functions of PM$\phi$E at $-1<\lambda<0$ also indicates the boundedness of the IFs of the Wald-type tests functional $W_n$.
\end{remark}

\bigskip
\section{Proof of Results \label{app:proofs}}

\subsection{Proof of Theorem \ref{th:wald_test_phi}}
\begin{proof}
We have by (\ref{eq:asymEst}) 
\begin{equation*}
\sqrt{n}(\widehat{\boldsymbol{\beta }}_{\phi ,P}-\boldsymbol{\beta }^0)\underset{n\rightarrow \infty }{\overset{\mathcal{L}}{\longrightarrow }}\mathcal{N}\left( \mathbf{0}_{d+1},\mathbf{V}\left( \boldsymbol{\beta }^0\right) \right) ,
\end{equation*}%
with $\mathbf{V}\left( \boldsymbol{\beta }^0\right)=\mathbf{J}^{-1}(\boldsymbol{\beta }^0)\mathbf{G}(\boldsymbol{\beta }^0)\mathbf{J}^{-1}(\boldsymbol{\beta }^0)$. Therefore,
\begin{equation*}
\sqrt{n}(\boldsymbol{M}^{T}\widehat{\boldsymbol{\beta }}_{\phi ,P}-\boldsymbol{m})\underset{n\rightarrow \infty }{\overset{\mathcal{L}}{\longrightarrow }}\mathcal{N}\left( \mathbf{0}_{d+1},\boldsymbol{M}^{T}\mathbf{V}\left( \boldsymbol{\beta }^0\right) \boldsymbol{M}\right) .
\end{equation*}%
As $\mathrm{rank}(\boldsymbol{M})=r$, we have that 
\begin{equation*}
n(\boldsymbol{M}^{T}\widehat{\boldsymbol{\beta }}_{\phi ,P}-\boldsymbol{m})^{T}\left( \boldsymbol{M}^{T}\mathbf{V}\left( \boldsymbol{\beta }^0\right)\boldsymbol{M}\right) ^{-1}(\boldsymbol{M}^{T}\widehat{\boldsymbol{\beta }}_{\phi ,P}-\boldsymbol{m})
\end{equation*}%
converge in law to a chi-square distribution with $r$ degrees of freedom. But $\mathbf{V}_n\left( \widehat{\boldsymbol{\beta }}_{\phi ,P}\right)$ is a consistent estimator of $\mathbf{V}\left(\boldsymbol{\beta }^0\right)$.  Therefore, we have that under $H_{0}$, $\boldsymbol{W}_{n}(\widehat{\boldsymbol{\beta}}_{\phi ,P})$ defined in (\ref{eq:Wald_phi}), converges in law to a chi-square distribution with $r$ degrees of freedom.
\end{proof}

\subsection{Proof of Theorem \ref{th:power1}}

\begin{proof}
A first order Taylor expansion of $\ell ^{\ast }(\widehat{\boldsymbol{\beta }}_{\phi ,P},\boldsymbol{\beta} ^{0})$ at $\widehat{\boldsymbol{\beta }}_{\phi ,P}$ around $\boldsymbol{\beta }^{0}$ gives
\begin{equation*}
\ell ^{\ast }(\widehat{\boldsymbol{\beta }}_{\phi ,P},\boldsymbol{\beta}^{0})-\ell ^{\ast }\left( \boldsymbol{\beta} ^{0},\boldsymbol{\beta}^{0}\right) =\left. \frac{\partial \ell ^{\ast }\left( \boldsymbol{\beta },\boldsymbol{\beta ^{0}}\right) }{\partial \boldsymbol{\beta }^{T}}\right \vert _{\boldsymbol{\beta} =\boldsymbol{\beta }^{0}}(\widehat{\boldsymbol{\beta }}_{\phi ,P}-\boldsymbol{\beta} ^{0})+o_{p}\left( \left \Vert 
\widehat{\boldsymbol{\beta }}_{\phi ,P}-\boldsymbol{\beta }^{}\right
\Vert \right) .
\end{equation*}
The asymptotic distribution of $\sqrt{n}\left(\ell ^{\ast }(\widehat{\boldsymbol{\beta }}_{\phi ,P},\boldsymbol{\beta} ^{0})-\ell^{\ast }\left( \boldsymbol{\beta} ^{0},\boldsymbol{\beta} ^{0}\right)\right) $ coincides with the asymptotic distribution of%
\begin{equation*}
\sqrt{n}\left(\left. \frac{\partial \ell ^{\ast }\left( \boldsymbol{\beta },\boldsymbol{\beta} ^{0}\right) }{\partial \boldsymbol{\beta }^{T}}\right \vert _{\boldsymbol{\beta =\beta }^{0}}(\widehat{\boldsymbol{\beta }}_{\phi ,P}-\boldsymbol{\beta} ^{0})\right),
\end{equation*}
but 
$$
\left. \frac{\partial \ell ^{\ast }\left( \boldsymbol{\beta },\boldsymbol{\beta} ^{0}\right) }{\partial \boldsymbol{\beta }^{T}}\right \vert _{\boldsymbol{\beta} =\boldsymbol{\beta }^0}=2\left( \boldsymbol{M}^{T}\boldsymbol{\beta }^{0}-\boldsymbol{m}\right) ^{T}\left( \boldsymbol{M}^{T}\mathbf{V}\left( \boldsymbol{\beta }^{0}\right) \boldsymbol{M}\right) ^{-1} \boldsymbol{M}^{T}.
$$
Now the result follows.
\end{proof}

\subsection{Proof of Theorem \ref{Th:PowerApprox_Wald}}

\begin{proof}

\begin{small}
\begin{align*}
\Pi_{W_{n}(\widehat{\boldsymbol{\beta }}_{\phi ,P})}\left(\boldsymbol{\beta}^{0}\right) & \simeq \textrm{P}_{\boldsymbol{\beta ^{0}}}(W_{n}(\widehat{\boldsymbol{\beta }}_{\phi ,P})>\chi_{r,\alpha }^{2})=\textrm{P}_{\boldsymbol{\beta ^{0}}}(n\ell^{\ast }(\widehat{\boldsymbol{\beta }}_{\phi ,P},\widehat{\boldsymbol{\beta }}_{\phi ,P})>\chi _{r,\alpha }^{2}) \\
& =\textrm{P}_{\boldsymbol{\beta ^{0}}}(\sqrt{n}\ell ^{\ast }(\widehat{\boldsymbol{\beta }}_{\phi ,P},\widehat{\boldsymbol{\beta }}_{\phi,P})>\frac{\chi _{r,\alpha }^{2}}{\sqrt{n}}) \\
& =\text{\textrm{P}}_{\boldsymbol{\beta ^{0}}}\left( \sqrt{n}\left(\ell ^{\ast }(\widehat{\boldsymbol{\beta }}_{\phi ,P},\boldsymbol{\beta}^{0})-\ell ^{\ast }\left( \boldsymbol{\beta}^{0},\boldsymbol{\beta}^{0}\right) \right) >\frac{\chi _{r,\alpha }^{2}}{\sqrt{n}}-\sqrt{n}\ell ^{\ast }\left( \boldsymbol{\beta}^{0},\boldsymbol{\beta}^{0}\right) \right) \\
& =1-\textrm{P}_{\boldsymbol{\beta}^{0}}\left( \frac{\sqrt{n}\left( \ell ^{\ast }(\widehat{\boldsymbol{\beta }}_{\phi ,P},\boldsymbol{\beta}^{0})-\ell ^{\ast }\left( \boldsymbol{\beta}^{0},\boldsymbol{\beta}^{0}\right) \right) }{\sigma _{W}\left( \boldsymbol{\beta }^{0}\right) }\leq \frac{1}{\sigma _{W}\left( \boldsymbol{\beta }^{0}\right) }\left( \frac{\chi _{r,\alpha }^{2}}{\sqrt{n}}-\sqrt{n}\ell ^{\ast}\left( \boldsymbol{\beta }^{0},\boldsymbol{\beta }^{0}\right)\right) \right) \\
& =1-\Phi _{n}\left( \frac{\sqrt{n}}{\sigma _{W}\left( \boldsymbol{\beta}^{\ast}\right) }\left( \frac{\chi _{r,\alpha }^{2}}{n}-\ell ^{\ast }\left( \boldsymbol{\beta}^{0},\boldsymbol{\beta}^{0}\right) \right)\right) ,
\end{align*}%
\end{small}
\noindent where $\Phi _{n}\left( x\right) $ uniformly tends to the standard normal distribution $\phi \left( x\right) $ as $n\rightarrow \infty$. 
\end{proof}

\subsection{Proof of Theorem \ref{th:respower_ult}}

\begin{proof}
It is easy to see that 
\begin{align*}
\boldsymbol{M}^{T}\widehat{\boldsymbol{\beta }}_{\phi ,P}-\boldsymbol{m} =n^{-1/2}\boldsymbol{M}^{T}\boldsymbol{\boldsymbol{d+M}}^{T}(\widehat{\boldsymbol{\beta }}_{\phi ,P}-\boldsymbol{\beta }_{n}).
\end{align*}%

We know, under $H_{1,n}$, that $\sqrt{n}(\widehat{\boldsymbol{\beta }}_{\phi ,P}-\boldsymbol{\beta }_{n})\underset{n\rightarrow \infty }{\overset{\mathcal{L}}{\longrightarrow }}\mathcal{N}\left( \mathbf{0}_{d+1},\mathbf{V}(\boldsymbol{\beta }_{n})\right)$ and $\boldsymbol{\beta }_{n}\underset{n\rightarrow \infty }{\longrightarrow }\boldsymbol{\beta }_{0}$. Therefore, as $\boldsymbol{M}^{T}\boldsymbol{\beta }_{0}=\boldsymbol{m}$,

\begin{equation*}
\sqrt{n}(\boldsymbol{M}^{T}\widehat{\boldsymbol{\beta }}_{\phi ,P}-\boldsymbol{m})\underset{n\rightarrow \infty }{\overset{\mathcal{L}}{\longrightarrow }}\mathcal{N}\left( \boldsymbol{M}^{T}\boldsymbol{\boldsymbol{d}},\boldsymbol{M}^{T}\mathbf{V}(\boldsymbol{\beta }_{0})\boldsymbol{M}\right) .
\end{equation*}%
But it is known that if $\boldsymbol{Z}\in \mathcal{N}\left( \boldsymbol{\mu ,\Sigma }\right)$, $\boldsymbol{\Sigma }$ is a symmetric projection of rank $k$ and $\boldsymbol{\Sigma \mu }=\boldsymbol{\mu }$, then $\boldsymbol{Z}^{T}\boldsymbol{Z}$ is a chi-square distribution with $k$ degrees of freedom and non-centrality parameter $\boldsymbol{\mu }^{T}\boldsymbol{\mu}$. So considering the  quadratic form 

\begin{equation*}
W_{n}(\widehat{\boldsymbol{\beta }}_{\phi ,P})=\boldsymbol{Z}^{T}\boldsymbol{Z},
\end{equation*}
with \ $\boldsymbol{Z}=\sqrt{n}\left[ \boldsymbol{M}^T\mathbf{V}\left(\widehat{\boldsymbol{\beta}}_{\phi,P}\right)^{-1} \boldsymbol{M}\right]^{-1/2}(\boldsymbol{M}^{T}\widehat{\boldsymbol{\beta}}_{\phi,P}-\boldsymbol{m})$ \ and 

\begin{equation*}
\boldsymbol{Z}\underset{n\rightarrow \infty }{\overset{\mathcal{L}}{\longrightarrow }}\mathcal{N}\left( \left[ \boldsymbol{M}^{T}\mathbf{V}(\boldsymbol{\beta }_{0})\boldsymbol{M}\right] ^{-1/2}\boldsymbol{M}^{T}\boldsymbol{d},\boldsymbol{I}_{r\times r}\right) ,
\end{equation*}%
the application  is immediate with the non-centrality parameter being

$$\boldsymbol{d}^{T}\boldsymbol{M}\left[ \boldsymbol{M}^{T}\mathbf{V}(\boldsymbol{\beta }_{0})\boldsymbol{M}\right] ^{-1}\boldsymbol{\boldsymbol{M}}^{T}\boldsymbol{d}.$$

The point (b) is straightforward taking into account that  the equivalence between the hypotheses (\ref{eq:alt1}) and (\ref{eq:alt2}) is given by  $\boldsymbol{M}^{T}\boldsymbol{\boldsymbol{d}}=\boldsymbol{\delta }$.
\end{proof}

\bigskip
\bigskip

\section{Some extensions of the Simulation Study \label{sec:app_sim}}

\begin{figure}[p]
\centering
\begin{tabular}{rc}
 \includegraphics[scale=0.4]{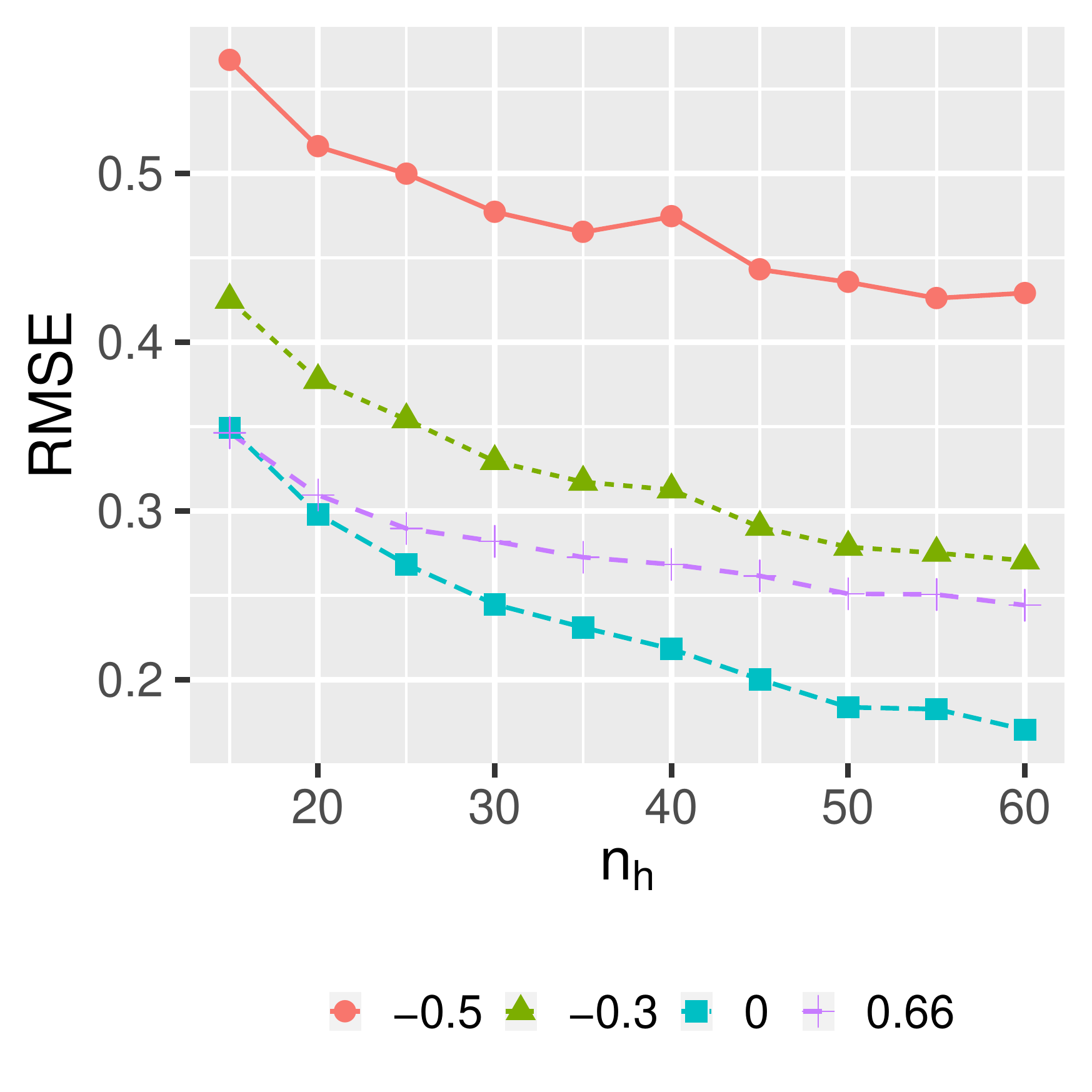}&
\includegraphics[scale=0.4]{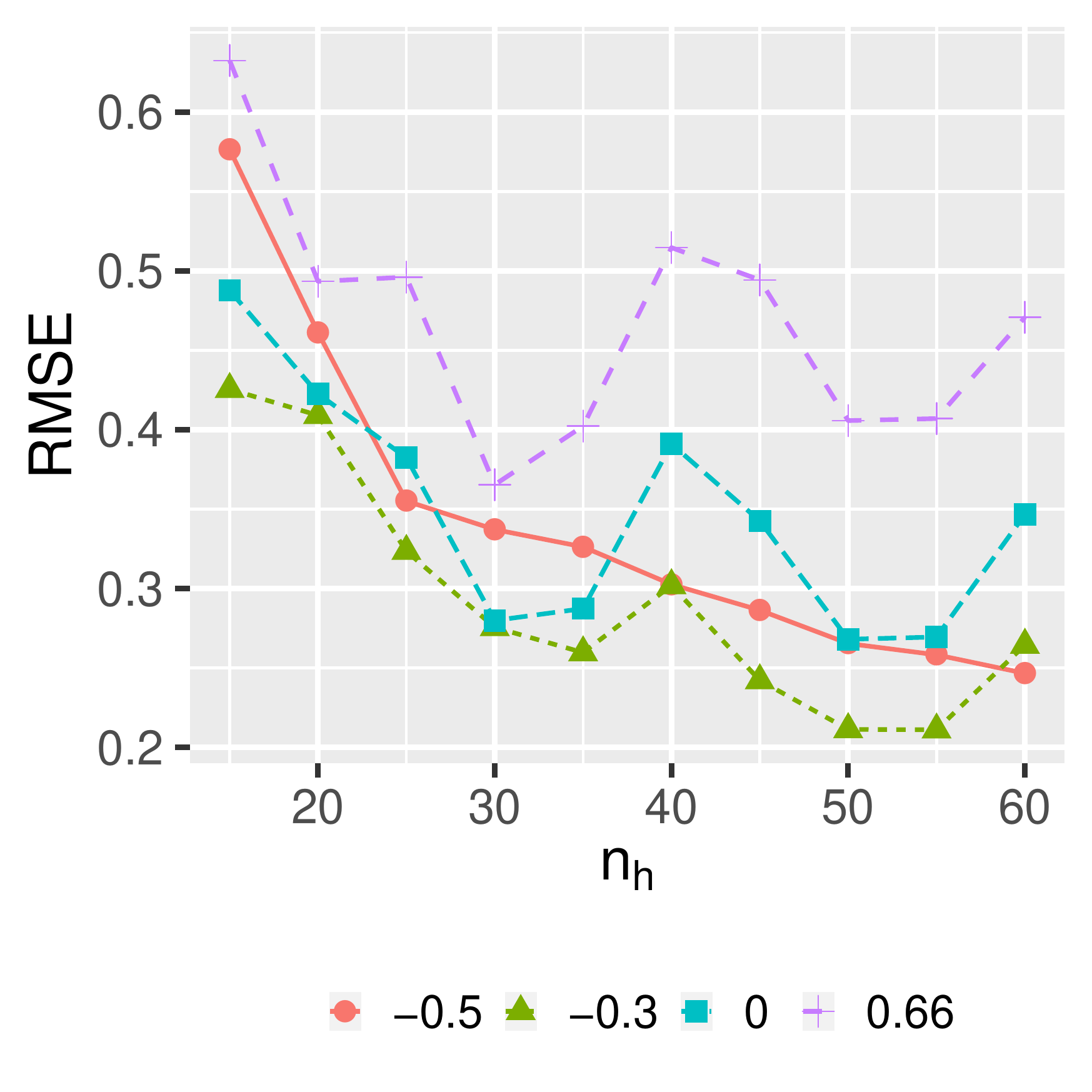}\\
 \includegraphics[scale=0.4]{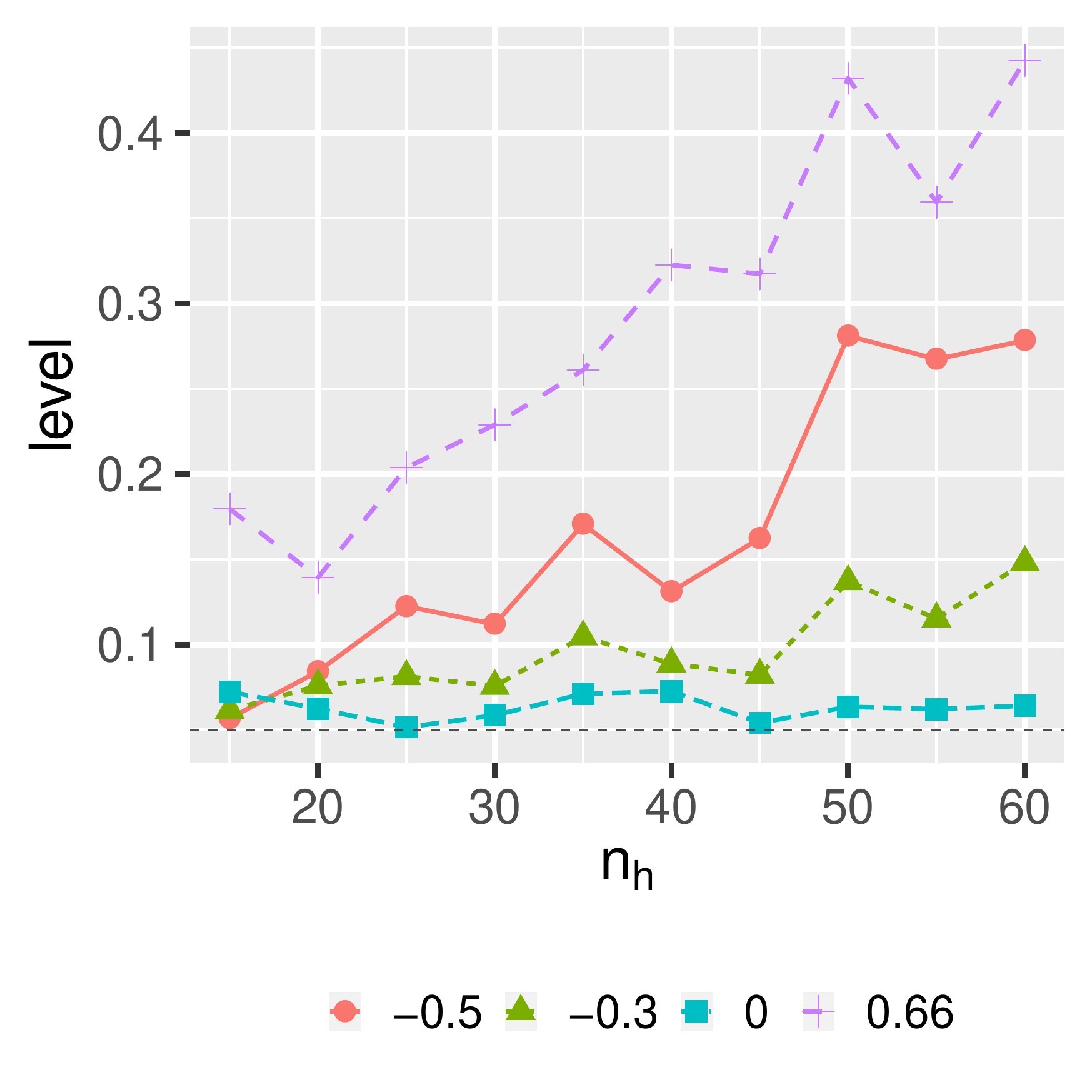}&
\includegraphics[scale=0.4]{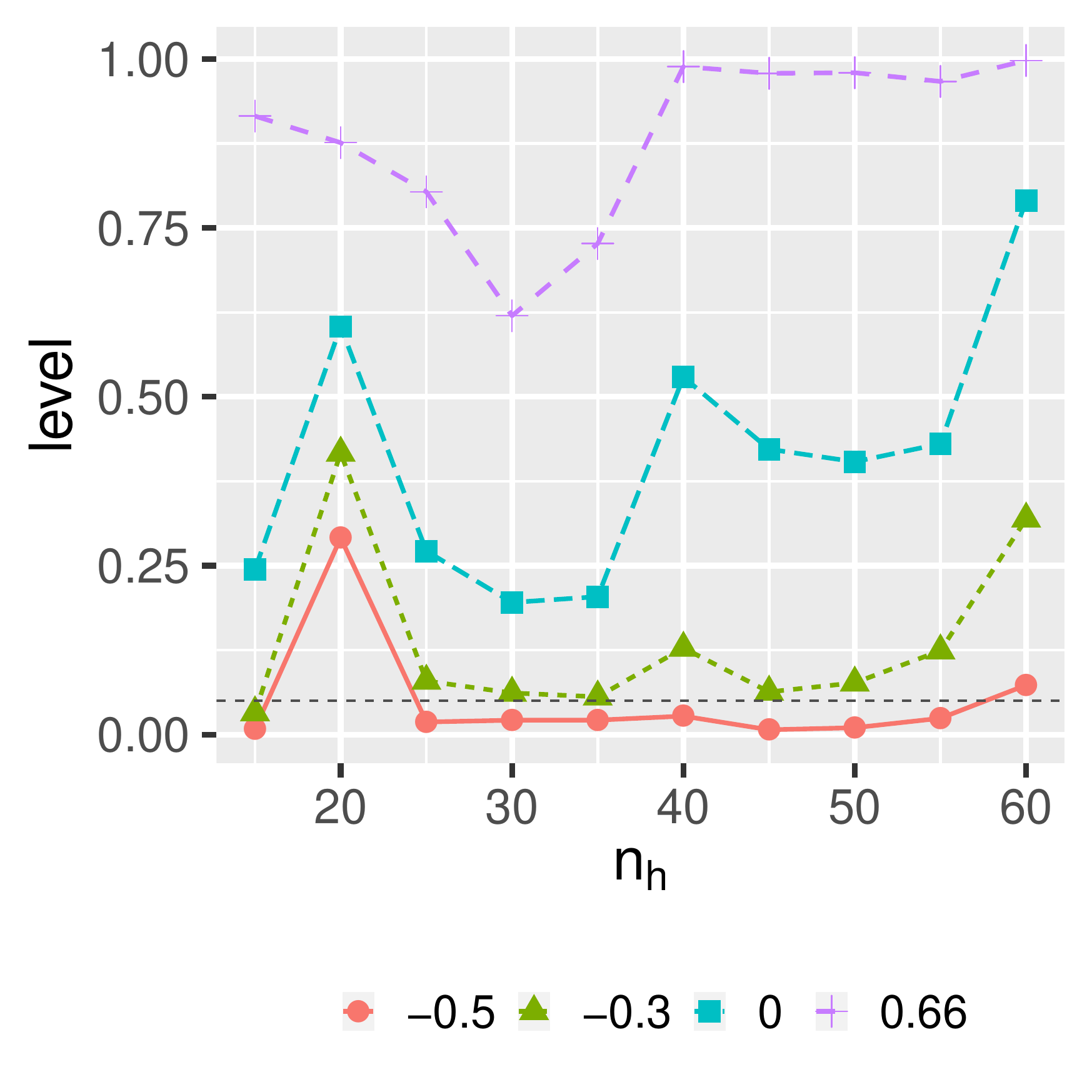}\\
 \includegraphics[scale=0.4]{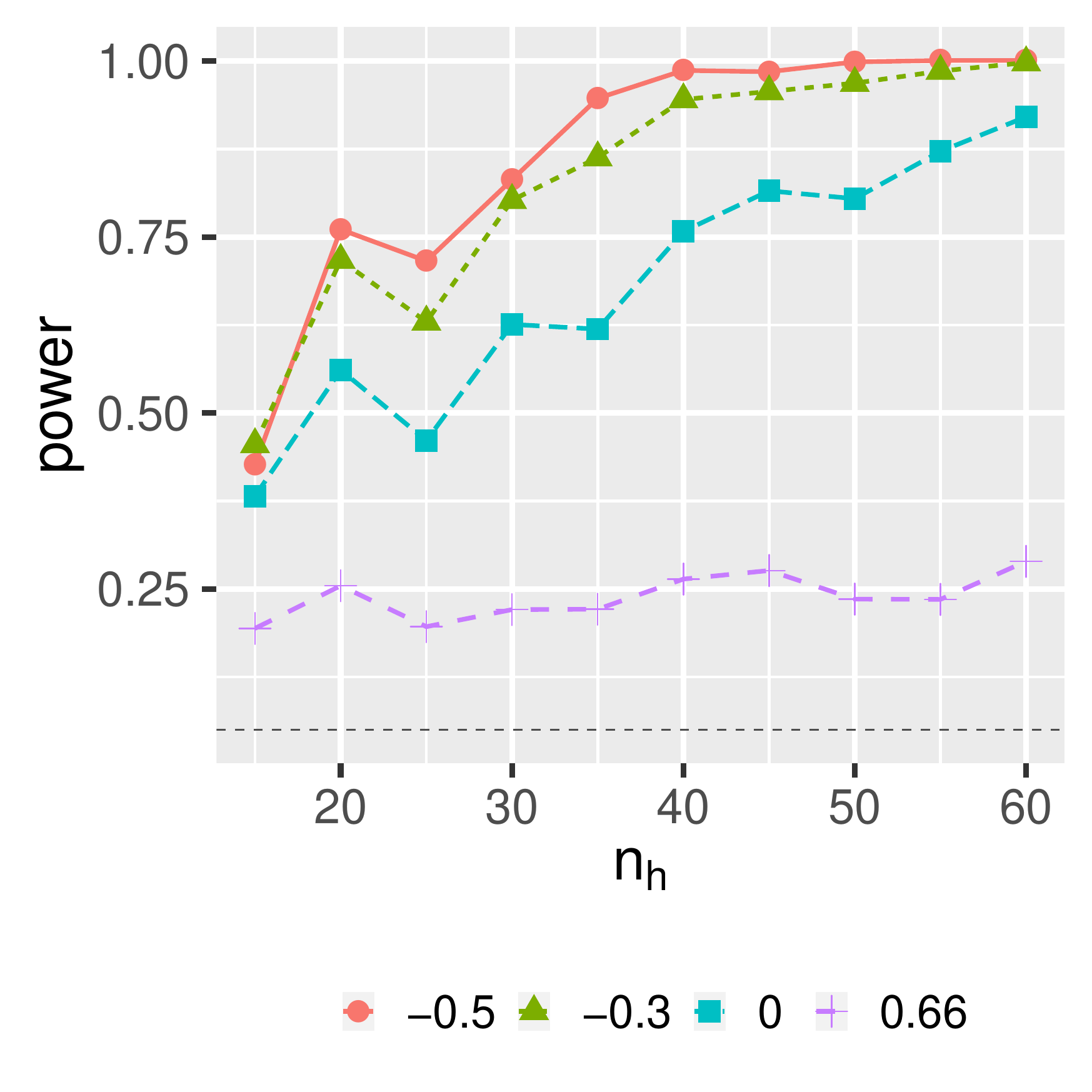}&
\includegraphics[scale=0.4]{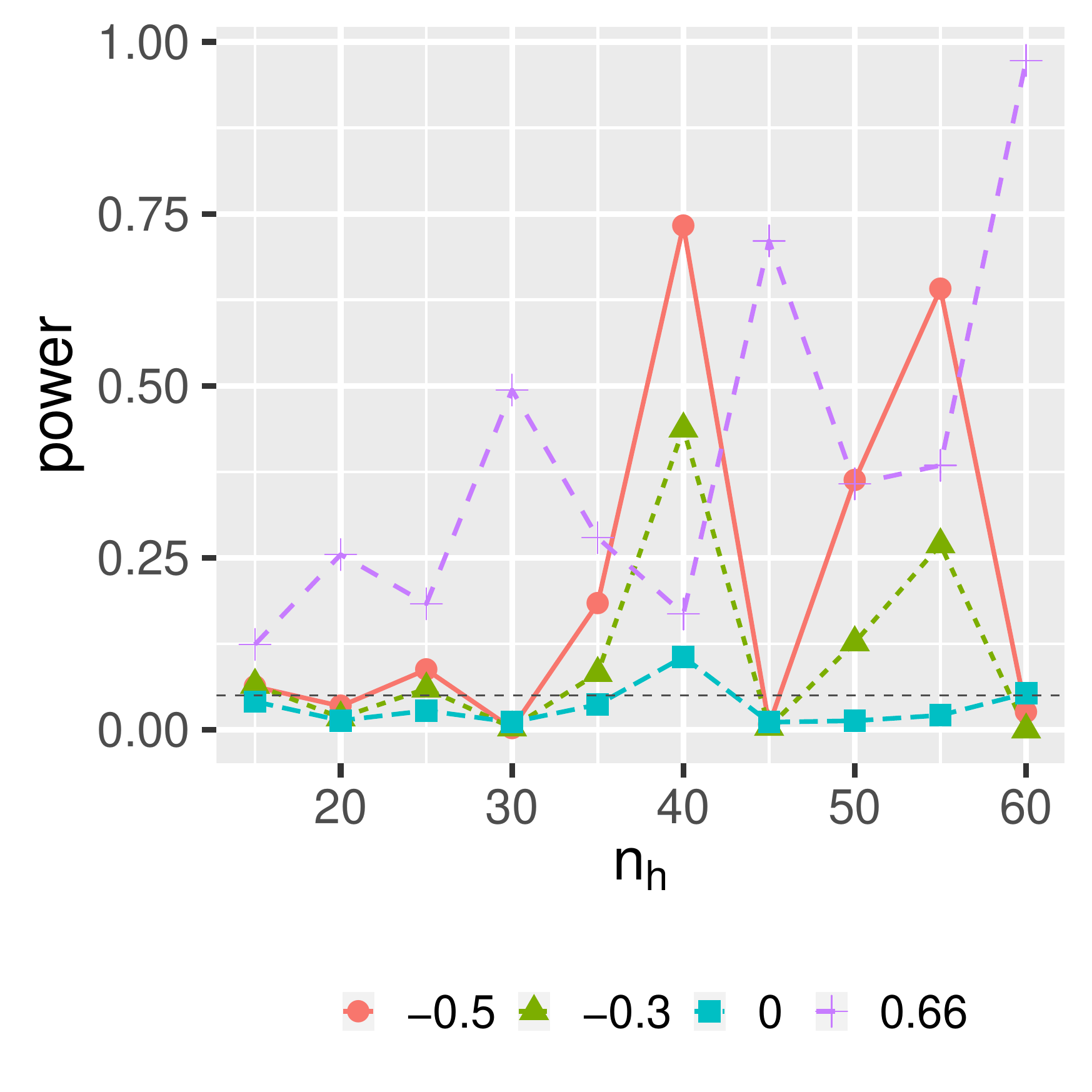}
\end{tabular}
\caption{RMSEs (top), emprirical levels (middle) and empirical powers (bottom).  Non-contaminated and contaminated settings (left and right, respectively). Random Clumped distribution  \label{fig:MC_RC}}
\end{figure}

\begin{figure}[p]
\centering
\begin{tabular}{rc}
 \includegraphics[scale=0.4]{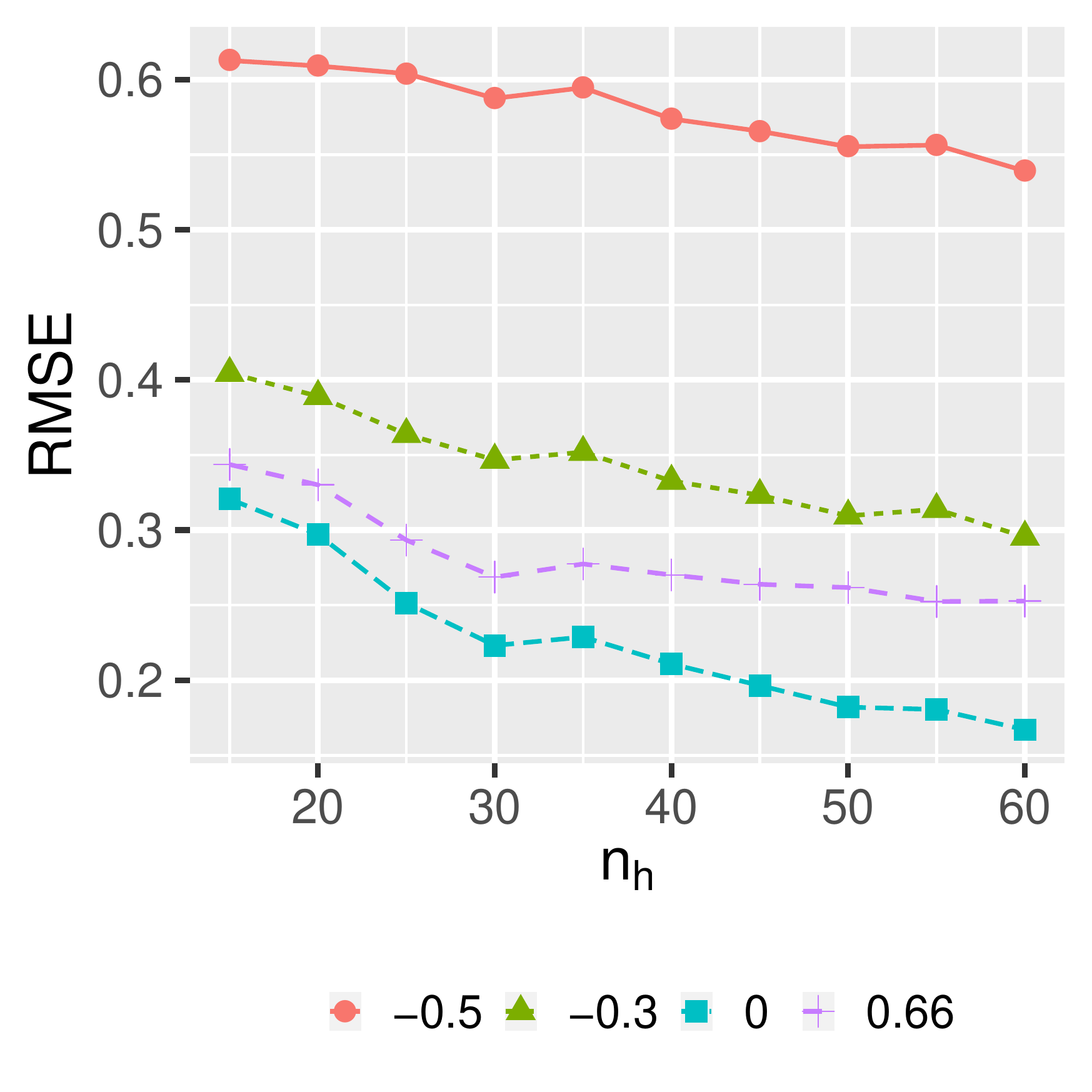}&
\includegraphics[scale=0.4]{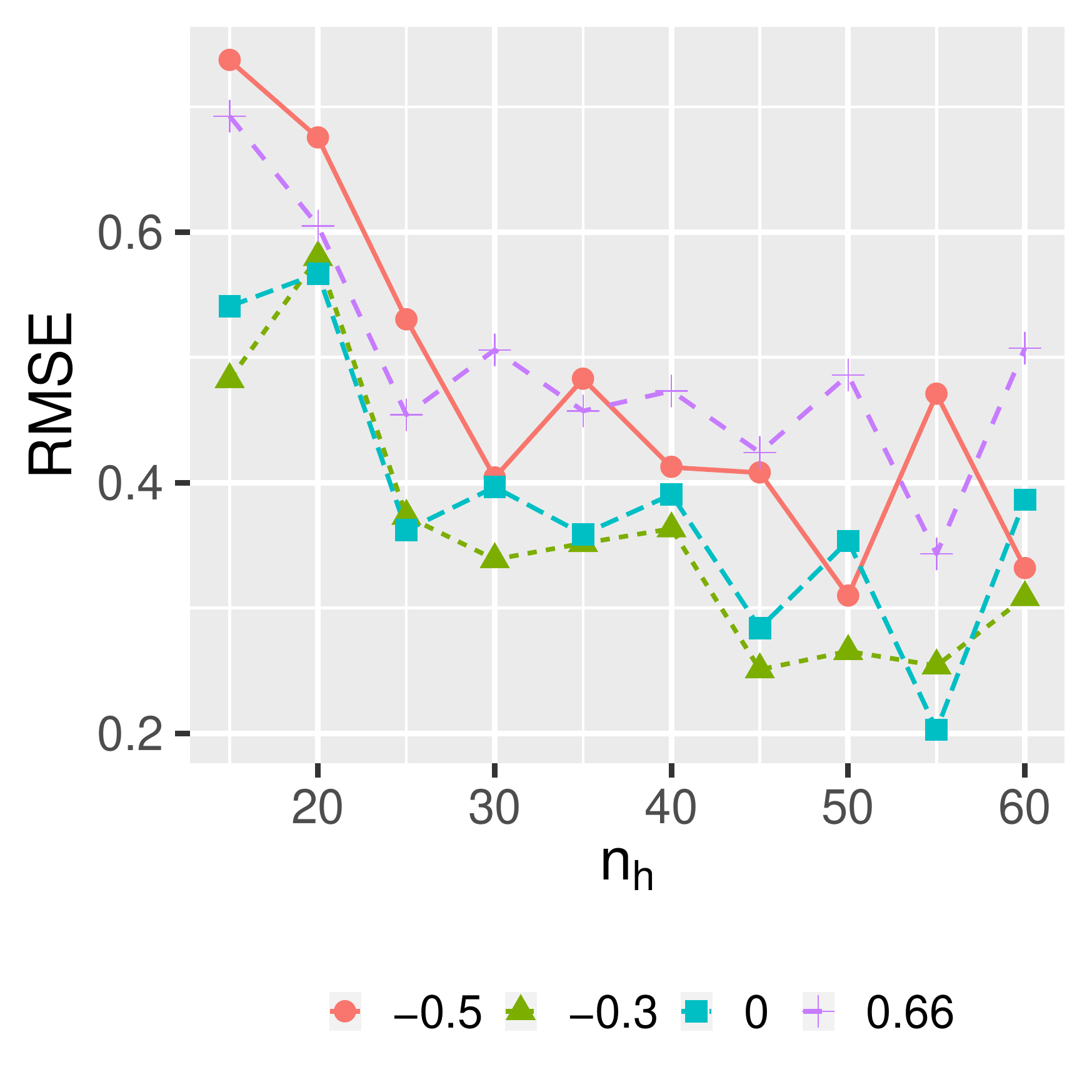}\\
 \includegraphics[scale=0.4]{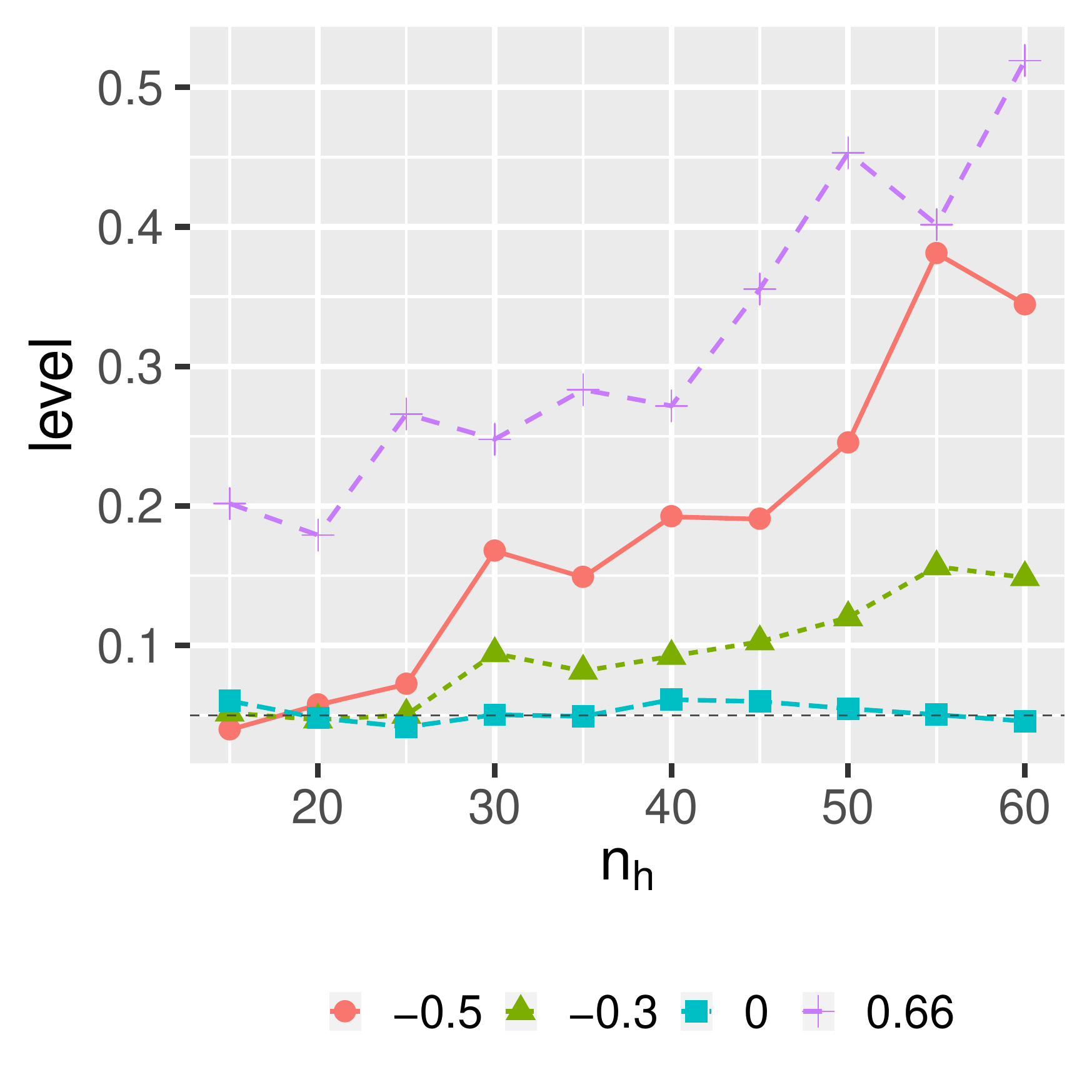}&
\includegraphics[scale=0.4]{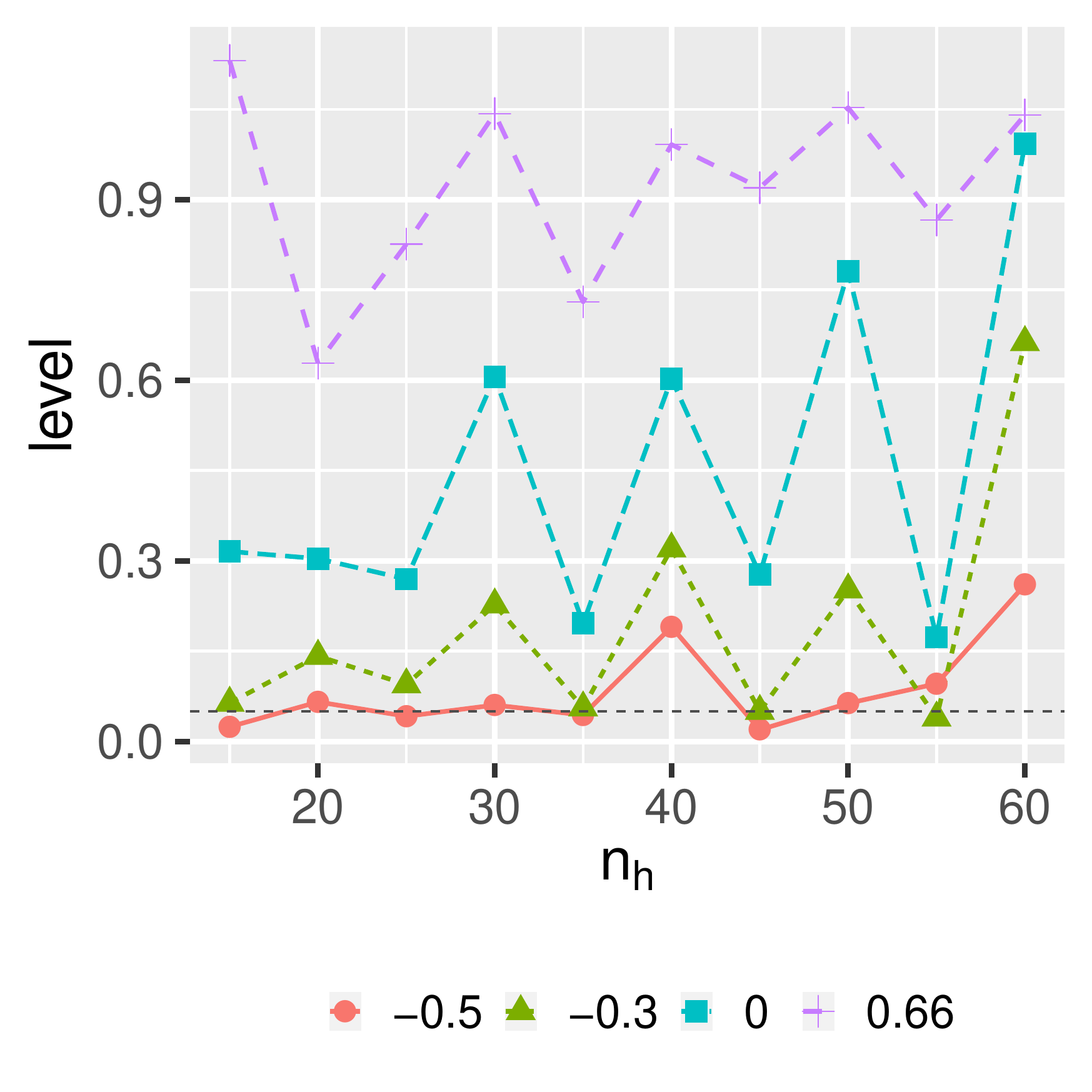}\\
 \includegraphics[scale=0.4]{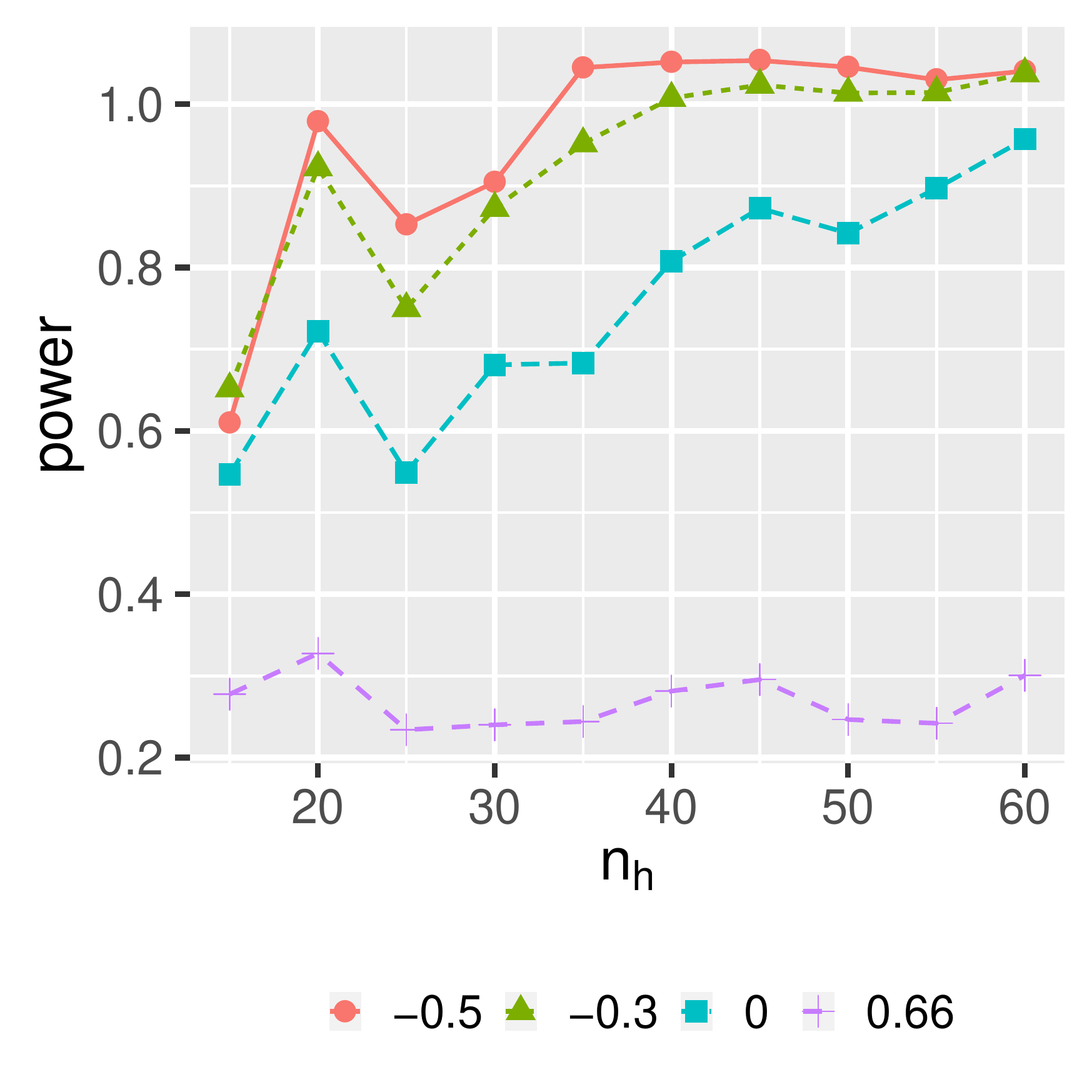}&
\includegraphics[scale=0.4]{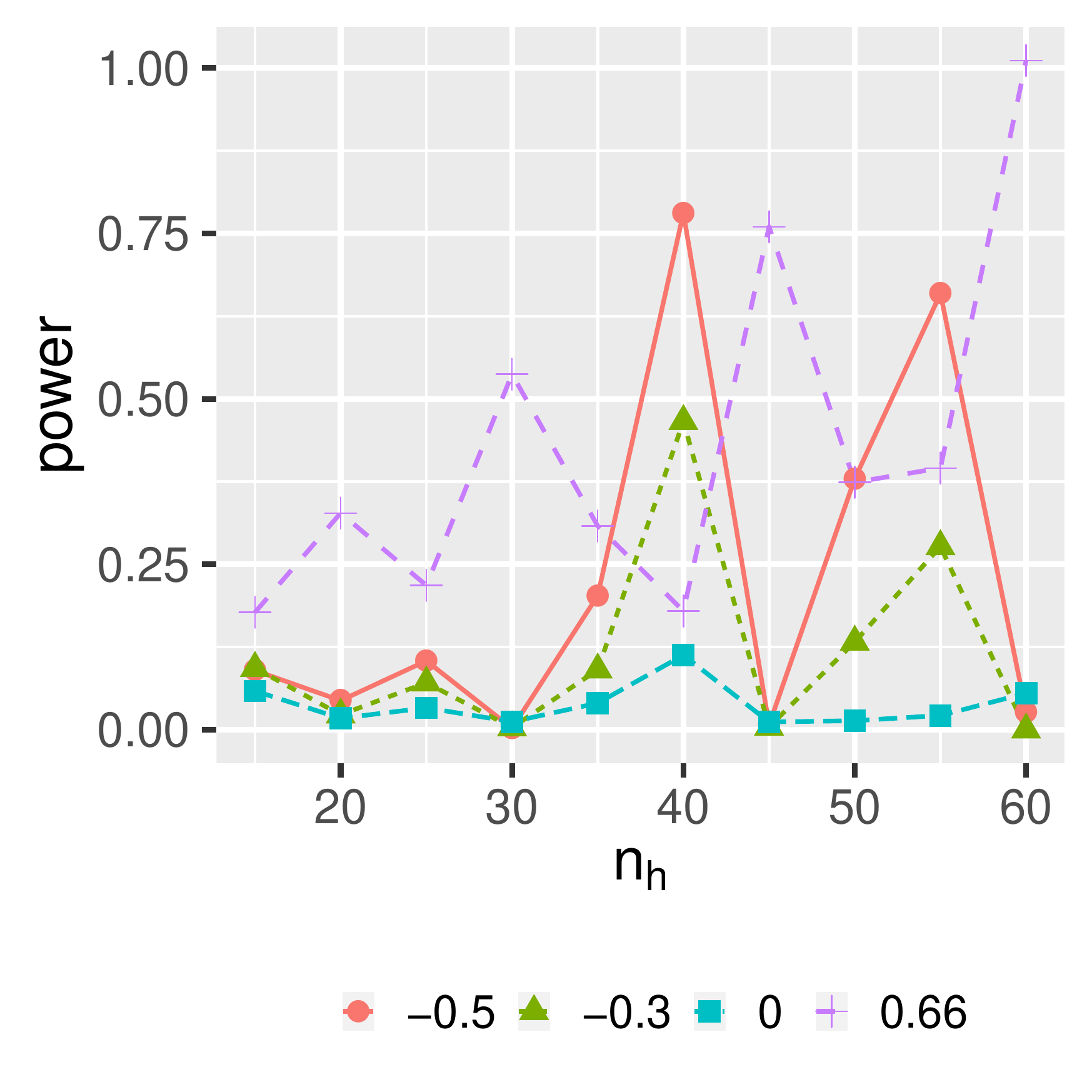}
\end{tabular}
\caption{RMSEs (top), emprirical levels (middle) and empirical powers (bottom).  Non-contaminated and contaminated settings (left and right, respectively). Dirichlet Multinomial distribution. \label{fig:MC_DM}}
\end{figure}

We extend the simulation study presented in Section \ref{sec:mc} to other scenarios. In particular, we first study the same scheme as in Section \ref{sec:mc} but considering two different overdispersed distributions for the response variable: the Random-Clumped and the Dirichlet Multinomial    distributions. Results are presented in Figure \ref{fig:MC_RC} and Figure \ref{fig:MC_DM}.  

Same conclusions as in Section \ref{sec:mc} are obtained, illustrating again the robustness of proposed estimators and Wald-type tests against classical PMLE.

\subsection{Algorithms for m-Inflated, Random Clumped and Dirichlet multinomial   distributions in the context of PLR models with complex design \label{sec:algorithms}}

We present the algorithms that are needed in order to compute the  Random-Clumped  (\cite{Morel1993}), Dirichlet-multinomial (\cite{Mosimann1962}) and m-Inflated (\cite{Cohen1976}) multinomial distributions in the context of PLR models with complex design. We consider, without loss of generality, that the intra-cluster correlation parameter is equal in all the clusters and strata ($\rho_{hi}\equiv \rho$, $h=1,\dots, H$, $i=1,\dots,n_h$).

\subsubsection*{m-Inflated distribution}

 \begin{algorithmic}[1]
  \small
 \renewcommand{\algorithmicrequire}{\textbf{Goal:}}
\REQUIRE Generation of response variable with m-Inflated distribution of parameters $\rho$ and $\boldsymbol{\pi}(\boldsymbol{\beta})$, in a scenario with $H$ strata and $n_h$ clusters in the stratum $h$, $h=1,\dots,H$.
  \FOR  {$h=1,\dots ,H$}
    \FOR  {$i=1,\dots ,n_h$}
		\STATE $k_1  \gets Ber(\rho^2)$
		\IF {$k_1=0$}
			\STATE  $\boldsymbol{y}_{hi}  \gets \mathcal{M}(m_{hi},\boldsymbol{\pi}_{hi}(\boldsymbol{\beta}))$
		\ELSE 
			\STATE $\boldsymbol{y}_{hi}  \gets m_{hi}\times\mathcal{M}(1,\boldsymbol{\pi}_{hi}(\boldsymbol{\beta}))$
		\ENDIF
  	\ENDFOR
  \ENDFOR
  \RETURN $\boldsymbol{y}=(\boldsymbol{y}_{11},\dots, \boldsymbol{y}_{Hn_{H}})^T$

 \end{algorithmic}

\subsubsection*{Random-Clumped  distribution}

 \begin{algorithmic}[1]
 \small
 \renewcommand{\algorithmicrequire}{\textbf{Goal:}}
\REQUIRE Generation of response variable with Random Clumped distribution of parameters $\rho$ and $\boldsymbol{\pi}(\boldsymbol{\beta})$, in a scenario with $H$ strata and $n_h$ clusters in the stratum $h$, $h=1,\dots,H$.
  \FOR  {$h=1,\dots ,H$}
    \FOR  {$i=1,\dots ,n_h$}
		\STATE $\boldsymbol{y}^{(0)} \gets \mathcal{M}(1,\boldsymbol{\pi}_{hi}(\boldsymbol{\beta}))$
		\STATE $k_1  \gets Bin(m_{hi},\rho)$
		\STATE $\boldsymbol{y}^{(1)}  \gets \mathcal{M}(m_{hi}-k_1,\boldsymbol{\pi}_{hi}(\boldsymbol{\beta}))$
		\STATE $\boldsymbol{y}_{hi}  \gets \boldsymbol{y}^{(0)}\times k_1+\boldsymbol{y}^{(1)}$

  	\ENDFOR
  \ENDFOR
\RETURN $\boldsymbol{y}=(\boldsymbol{y}_{11},\dots, \boldsymbol{y}_{Hn_{H}})^T$
 \end{algorithmic}

\subsubsection*{Dirichlet-Multinomial distribution}

 \begin{algorithmic}[1]
  \small
 \renewcommand{\algorithmicrequire}{\textbf{Goal:}}
\REQUIRE Generation of response variable with Dirichlet Multinomial distribution of parameters $\rho$ and $\boldsymbol{\pi}(\boldsymbol{\beta})$, in a scenario with $H$ strata and $n_h$ clusters in the stratum $h$, $h=1,\dots,H$.
  \FOR  {$h=1,\dots ,H$}
    \FOR  {$i=1,\dots ,n_h$}
		\STATE $\alpha_1 \gets \frac{1-\rho^2}{\rho^2}\pi_{hi1}(\boldsymbol{\beta})$
		\STATE $\alpha_2 \gets \frac{1-\rho^2}{\rho^2}(1-\pi_{hi1}(\boldsymbol{\beta}))$
		\STATE $y_{hi1}\gets Bin(m_{hi},Beta(\alpha_1,\alpha_2))$
			\FOR{$r=1,\dots,d$}
				\STATE $\alpha_1 \gets \frac{1-\rho^2}{\rho^2}\pi_{hir}(\boldsymbol{\beta})$
		\STATE $\alpha_2 \gets \frac{1-\rho^2}{\rho^2}(1-\sum_{l=1}^r\pi_{hil}(\boldsymbol{\beta}))$
		\STATE $y_{hir}\gets Bin(m_{hi}-\sum_{l=1}^{r-1}y_{hil},Beta(\alpha_1,\alpha_2))$
		\ENDFOR
		\STATE $y_{hi,d+1} \gets m_{hi}-\sum_{r=1}^{d}y_{hir}$
		\STATE $\boldsymbol{y}_{hi}\gets(y_{hi,1},\dots, y_{hi,d+1})$
  	\ENDFOR
  \ENDFOR
  \RETURN $\boldsymbol{y}=(\boldsymbol{y}_{11},\dots, \boldsymbol{y}_{Hn_{H}})^T$

 \end{algorithmic} 

\clearpage
\bibliography{bibliography2}
\bibliographystyle{abbrvnat}

\end{document}